%% file: 2-14.tex
\documentclass[a4,10pt]{article}
\usepackage{amsmath,amscd,amssymb}
\usepackage[margin=.8in]{geometry}

\input{notation}
\input{new_theorem}

\begin{document}

\title{Specialized Simpson's main estimates
for cyclic harmonic $G$-bundles}

\author{Takuro Mochizuki}
\date{}
\maketitle

\begin{abstract}
We study a generalization of
specialized Simpson's main estimate
in the context of cyclic harmonic $G$-bundles
induced by split automorphisms.
We apply it to the classification of 
Toda type $G$-harmonic bundles.
 
\vspace{.1in}
\end{abstract}

\section{Introduction}

\subsection{Harmonic bundles on Riemann surfaces}

Let $(E,\delbar_E)$ denote a holomorphic vector bundle
on a Riemann surface $X$,
and let $\theta$ be
an $\End(E)$-valued holomorphic $1$-form.
Such a pair $(E,\delbar_E,\theta)$ is called a Higgs bundle.
For a Hermitian metric $h$ of $E$,
we obtain the Chern connection $\nabla_h$
which is a unitary connection whose $(0,1)$-part equals $\delbar_E$.
We obtain the curvature $R(h)$ of the Chern connection.
We also obtain the adjoint $\theta^{\dagger}_h$
of the Higgs field $\theta$ with respect to $h$.
The metric $h$ is called a harmonic metric of
the Higgs bundle $(E,\delbar_E,\theta)$
if $R(h)+[\theta,\theta^{\dagger}_h]=0$.
A Higgs bundle with a harmonic metric
is called a harmonic bundle.
This condition is equivalent to
the integrability of the connection
$\DD^1_h=\nabla_h+\theta+\theta^{\dagger}_h$.

By definition,
a harmonic bundle $(E,\delbar_E,\theta,h)$
has the underlying Higgs bundle $(E,\delbar_E,\theta)$
and the underlying flat bundle $(E,\DD^1_h)$.
By the following theorem proved in
\cite{corlette, don2, Hitchin-self-duality, s1},
the concept of harmonic bundle
has been one of useful bridges to connect
topology and differential geometry,
and complex geometry and algebraic geometry of Riemann surfaces.

\begin{thm}[Corlette-Donaldson-Hitchin-Simpson]
If $X$ is compact,
the above correspondences induce
equivalences of harmonic bundles,
polystable Higgs bundles of degree $0$,
and semisimple flat bundles.
\hfill\qed
\end{thm}

There are generalizations in various directions.
For example, see
\cite{Biquard-Boalch,Biquard-Garcia-Prada-Mundet,
Bradlow-Garcia-Prada-Mundet,
corlette, JZ2,
Mochizuki-wild, Mochizuki-KH-Higgs, s2,s5}.

\subsubsection{Simpson's main estimate}

Let us recall another basic result
called Simpson's main estimate,
which is useful in analysis of harmonic bundles.
Simpson \cite{s2} pioneered such estimates
in his study of tame harmonic bundles.
It has been further developed in
\cite{Mochizuki-wild, Decouple}.
We briefly recall the estimates given in \cite{Decouple}.

Let $\Harm(E,\theta)$
denote the set of harmonic metrics of
a Higgs bundle $(E,\delbar_E,\theta)$.
Let $\Sigma_{\theta}$ be the spectral curve of $\theta$.
We fix a conformal metric $g_X$ of $X$.
Let $W$ denote any compact subset of $X$.

\begin{thm}
\label{thm;25.9.26.40}
There exists a positive constant $A_0$
depending only on $(\Sigma_{\theta},g_X,W,\rank(E))$
such that the following holds.
\begin{itemize}
 \item For any $t>0$ and for any $h\in\Harm(E,t\theta)$,
       we obtain
       $|t\theta|_{h,g_X}\leq A_0(1+t)$ on $W$.
       \hfill\qed
\end{itemize}
\end{thm}

Let us consider the case where
there exists a decomposition of the Higgs bundle
$(E,\theta)=(E_1,\theta_1)\oplus (E_2,\theta_2)$
such that
$\Sigma_{\theta_1}\cap \Sigma_{\theta_2}=\emptyset$.

\begin{thm}
\label{thm;25.9.26.31}
There exist positive constants $A_1,A_2$
depending only on $(\Sigma_{\theta},g_X,W,\rank(E))$
such that the following holds.
 \begin{itemize}
  \item For any $t\geq 1$,
	$h\in\Harm(E,t\theta)$,
       $P\in W$,
       and $u_i\in E_{i|P}$ $(i=1,2)$, we obtain
\[
	\bigl|h(u_1,u_2)\bigr|
       \leq
       A_1\exp(-A_2t)
	|u_1|_h\cdot |u_2|_h.
\]
\hfill\qed
\end{itemize}
\end{thm}

According to Theorem \ref{thm;25.9.26.31},
if the spectral curves are far enough apart,
then the decomposition is almost orthogonal.
Variants of these theorems are
fundamental in the study of singular and higher dimensional
harmonic bundles
in \cite{s2} and \cite{mochi4,Mochizuki-wild,Decouple}.
It was also applied in \cite{Sagman-Smillie}
to prove some harmonic bundles are counterexamples to 
Labourie's conjecture.

\subsubsection{Specialized Simpson's main estimate
in the case of symmetric Higgs bundles}

Under some additional assumptions,
we obtain some useful estimates
for comparison of harmonic metrics of $(E,\delbar_E,\theta)$
from Theorem \ref{thm;25.9.26.40}
and Theorem \ref{thm;25.9.26.31},
which we call specialized Simpson's main estimate.
Let us explain it in the symmetric case.

Let $C$ be a non-degenerate holomorphic symmetric pairing of
the Higgs bundle $(E,\delbar_E,\theta)$.
Namely,
$C:E\otimes E\to\nbigo_X$ denotes a symmetric product
such that
the induced map $\Psi_C:E\to E^{\lor}$ is an isomorphism,
where $E^{\lor}$ denotes the dual of $E$,
and $\theta$ is self-adjoint with respect to $C$
in the sense
$C(\theta u_1\otimes u_2)=C(u_1\otimes\theta u_2)$
for local sections $u_i$.
A harmonic metric $h$ of $(E,\delbar_E,\theta)$
is called compatible with $C$
if $\Psi_C$ is isometric with respect to $h$
and the induced metric $h^{\lor}$ of $E^{\lor}$.
Let $\Harm(E,\theta;C)$
denote the set of harmonic metrics
compatible with $C$.

Moreover, we assume that $\theta$ is regular semisimple,
i.e.,
$\Sigma_{\theta}\to X$ is a covering map
of degree $\rank(E)$.
The following lemma is explained in \cite{Li-Mochizuki3}.

\begin{lem}
There exists a uniquely determined harmonic metric $h^{\can}$
of $(E,\theta)$ compatible with $C$
such that it is decoupled in the sense
$R(h^{\can})=[\theta,\theta^{\dagger}_{h^{\can}}]=0$. 
In particular,
$h^{\can}\in \Harm(E,\theta;C)$.
\hfill\qed
\end{lem}

For any $h\in \Harm(E,\theta;C)$,
we obtain the endomorphism
$v(h^{\can},h)$ of $E$
such that
(i) $v(h^{\can},h)$ 
is
self-adjoint with respect to both $h^{\can}$ and $h$,
(ii) $v(h^{\can},h)$ satisfies
\[
 h(u_1,u_2)=h^{\can}\Bigl(\exp\bigl(2v(h^{\can},h)\bigr)u_1,u_2\Bigr).
\]
As a specialization of Simpson's main estimate,
we obtained the following in
\cite{Li-Mochizuki3, Mochizuki-Szabo}.
\begin{thm}
\label{thm;25.9.26.50}
There exist $A_i$ $(i=1,2)$
depending only on 
$(\Sigma_{\theta},g_X,W)$
such that the following holds.
\begin{itemize}
 \item For any $t\geq 1$ and $h\in \Harm(E,t\theta;C)$, we obtain
 $\bigl|
v(h^{\can},h)
\bigr|_{h^{\can}}
\leq
A_1
\exp\bigl(-tA_2\bigr)$
on $W$.
\hfill\qed
\end{itemize}
\end{thm}

As an application,
we applied 
Theorem \ref{thm;25.9.26.50}
in \cite{Li-Mochizuki3}
to prove $\Harm(E,\theta,C)\neq\emptyset$
if $\theta$ is generically regular semisimple,
i.e.,
there exists a discrete subset
$D\subset X$ such that
$(E,\theta)_{|X\setminus D}$
is regular semisimple.
Let us emphasize that Theorem \ref{thm;25.9.26.50}
is also useful even in the non-symmetric case.
For example,
suppose that (i) $X$ is compact,
(ii) $\Sigma_{\theta}$ is smooth, (iii) $\deg(E)=0$.
By fixing a flat metric $h_{\det(E)}$ of $\det(E)$,
we obtain harmonic metrics $h_t$ of $(E,t\theta)$
such that $\det(h_t)=h_{\det(E)}$.
In \cite{Mochizuki-Szabo},
we applied Theorem \ref{thm;25.9.26.50}
to study the limiting behaviour of $h_t$ as $t\to\infty$.
In \cite{Mochizuki-Hitchin-metric},
we also applied Theorem \ref{thm;25.9.26.50}
in the study of the asymptotic behaviour of
the Hitchin metric of the moduli space of
Higgs bundles on compact Riemann surfaces.

\subsubsection{Similar estimates for Toda equations}

We set
$\hyperk_{X,r}=\bigoplus_{j=1}^r K_X^{(r+1-2j)/2}$.
For any $r$-differential $q$,
let $\theta(q)$ denote the Higgs field of $\hyperk_{X,r}$
induced by
$K_X^{(r+1-2j)/2}\simeq
K_X^{(r+1-2j-2)/2}\otimes K_X$ $(j=1,\ldots,r-1)$
and
$q:K_X^{-(r-1)/2}\simeq K_X^{(r-1)/2}\otimes K_X$.
A harmonic metric $h$ of $(\hyperk_{X,r},\theta(q))$
are closely related with the Toda equation associated with $q$
if the decomposition 
$\hyperk_{X,r}=\bigoplus_{j=1}^r K_X^{(r+1-2j)/2}$
is orthogonal with respect to $h$.
In \cite{Li-Mochizuki2},
a similar estimate is obtained
for such harmonic metrics in the case where
$q$ is nowhere vanishing.
We applied it to study
the classification of such harmonic metrics.
(See \cite{Li-Mochizuki2} for more details.
See \cite{GuestItsLin, GuestLin,Li,Mochizuki-Toda}
for previous works for the classification of solutions.)

\subsubsection{Goal of this paper}

In this paper,
we shall study a generalization of
specialized Simpson's main estimates
to the context of cyclic $G$-Higgs bundles
induced by split automorphisms in a unified way.
We shall deduce the estimates from 
Theorem \ref{thm;25.9.26.40}
and Theorem \ref{thm;25.9.26.31}.

\begin{rem}
See {\rm\cite{Garcia-Prada-Gonzalez-cyclic, Sagman-Tosic}}
for recent studies of cyclic $G$-Higgs bundles
from different perspectives. 
\hfill\qed
\end{rem}

\subsection{Split automorphisms for simple complex algebraic groups}

\subsubsection{Split automorphism of finite order}

Let $G$ be a semisimple complex Lie group,
i.e.,
$G$ is a connected complex Lie group
whose Lie algebra $\gminig$ is semisimple.
An element $u\in\gminig$ is called semisimple
if $\ad(u)\in\End_{\cnum}(\gminig)$ is semisimple.
We set $C_{\gminig}(u)=\Ker\ad(u)$.
An element $u\in\gminig$ is called regular
if $\dim C_{\gminig}(u)=\dim\gminit$
for a Cartan subalgebra $\gminit$ of $\gminig$.
If $u$ is regular and semisimple,
then $C_{\gminig}(u)$ is a Cartan subalgebra.

Let $\sigma$ be a holomorphic
automorphism of $G$ of finite order $m>0$.
We set
$G^{\sigma}=\{g\in G\,|\,\sigma(g)=g\}$.
Let $\omega$ be a primitive $m$-th root of $1$.
We obtain the eigen decomposition
$(\gminig,\sigma)
=\bigoplus_{\ell=0}^{m-1}
(\gminig_{\ell},\omega^{\ell})$,
where $\sigma=\omega^{\ell}$ on $\gminig_{\ell}$.
Let $\gminig_1^{\rs}$ denote the set of
regular semisimple elements contained in $\gminig_1$.

\begin{df}
$(\sigma,\omega)$ is called split
if the following two conditions are satisfied. 
\begin{itemize}
 \item $\gminig_1^{\rs}\neq\emptyset$.
 \item For any $u\in\gminig_1^{\rs}$,
       we have
       $\gminig_0\cap C_{\gminig}(u)=0$.
\hfill\qed
\end{itemize}
\end{df}
       
\subsubsection{Canonical maximal compact groups}

Let $\Herm(G)$ denote the set of maximal compact subgroups.
For any $K\in\Herm(G)$,
let $\rho_K$ denote the involution of $\gminig$
defined by
$\rho_K(u_1+\sqrt{-1}u_2)=u_1-\sqrt{-1}u_2$
for $u_i\in \gminik$.
We set
$\Herm(G,\sigma)=\bigl\{
K\in\Herm(G)\,\big|\,\sigma(K)=K
\bigr\}$.

\begin{prop}[Corollary
\ref{cor;25.6.21.10} and Corollary \ref{cor;25.9.27.2}]
For any $u\in\gminig_1^{\rs}$,
there exists a unique
$K^{\can}(u)\in\Herm(G,\sigma)$
such that
$\rho_{K^{\can}(u)}(C_{\gminig}(u))=C_{\gminig}(u)$.
This induces a $C^{\infty}$-map
$K^{\can}:\gminig_1^{\rs}\to \Herm(G,\sigma)$.
If $u_1,u_2\in\gminig_1^{\rs}$ are commuting,
then $K^{\can}(u_1)=K^{\can}(u_2)$.
 \end{prop}

\subsubsection{Example 1}

A typical example of split automorphism is the involution
which induces a split real form.
\begin{prop}[Proposition
\ref{prop;25.9.27.3}
and Proposition \ref{prop;25.9.27.4}]
Let $\sigma$ be a holomorphic involution of $G$.
Then, $(\sigma,-1)$ is split
if and only if
$\sigma$  induces a split real form of $\gminig$
in the following sense.
\begin{itemize}
 \item Let $K\in\Herm(G)$ such that $\sigma(K)=K$.
       We have the eigen decomposition
       $(\gminik,\sigma)=
       (\gminik^{\sigma=1},1)
       \oplus
       (\gminik^{\sigma=-1},-1)$.
       Then,
       $\gminik^{\sigma=1}\oplus
       \sqrt{-1}\gminik^{\sigma=-1}$
       is a split real form of $\gminig$.
\end{itemize} 
 \end{prop}

\begin{rem}
For $G=\SL(n,\cnum)$,
we have the involution
given by
$g\longmapsto \lefttop{t}g^{-1}$.
\hfill\qed
\end{rem}

\subsubsection{Example 2}
\label{subsection;25.9.26.101}

Another typical example
appears in the classical work of Kostant \cite{Kostant-TDS}.
Suppose $\gminig$ is simple.
Let $G$ be of adjoint type.
Let $\gminit$ be a Cartan subalgebra.
We have the corresponding subgroup $T\subset G$.
We obtain the root decomposition
$\gminig=\gminit\oplus
 \bigoplus_{\phi\in\Delta}
 \gminig_{\phi}$.
We fix a Weyl chamber.
Let $\alpha_1,\ldots,\alpha_{\ell}\in \gminit^{\lor}$
denote the set of positive simple roots.
There exists a unique element
$\ttx_0\in\gminit$
determined by
$\alpha_i(\ttx_0)=1$ $(i=1,\ldots,\ell)$.
Let $\psi\in\Delta$
denote the highest root.
We obtain the positive integer
$\tth=\psi(\ttx_0)$.

We obtain the element
$\ttw=\exp\Bigl(
 \frac{2\pi\sqrt{-1}}{1+\tth}\ttx_0
 \Bigr)\in T$.
We consider the automorphism
$\sigma:G\to G$
defined by
$\sigma(g)=\ttw g \ttw^{-1}$.
We set $\omega=\exp(\frac{2\pi\sqrt{-1}}{\tth+1})$.
The following proposition is
essentially contained in the classical work of Kostant
\cite{Kostant-TDS}.

\begin{prop}[Proposition
\ref{prop;25.9.27.10}]
$(\sigma,\omega)$ is split.
\hfill\qed
 \end{prop}

\subsection{Specialized Simpson's main estimate for cyclic harmonic $G$-bundles}

\subsubsection{Harmonic $G$-bundles}

Let $G$ be a semisimple complex algebraic group
with the Lie algebra $\gminig$.
Let $\gbigp_G$ be a holomorphic principal $G$-bundle.
Let $\gbigp_G(\gminig)$ denote the bundle of Lie algebras
associated with the adjoint representation.

Let $K\subset G$ be a maximal compact subgroup.
Because $G$ is semisimple,
we have the natural isomorphism
$\gbigp_G(G/K)\simeq \gbigp_G(\Herm(G))$.
A section $h$ of $\gbigp_G(\Herm(G))$
corresponds to
a $K$-reduction $\gbigp^h_K\subset\gbigp_G$,
which is called a Hermitian metric of $\gbigp_G$ in this paper.
As in the case of vector bundles,
there exists the connection of $\gbigp^h_{K}$
associated with $(\gbigp_G,h)$.
We obtain the curvature $R(h)$ of the Chern connection
as a section of
$\gbigp_G(\gminig)\otimes\Omega^{1,1}$.

A $\gbigp_G(\gminig)$-valued
holomorphic $1$-form $\theta$
is called a Higgs field of $\gbigp_G$.
For a Hermitian metric $h$ of $\gbigp_G$,
we obtain the adjoint
$-\rho_{h}:\gbigp^{h}_{K}(\gminig)\otimes\Omega^{1,0}
\to \gbigp^{h}_{K}(\gminig)\otimes\Omega^{0,1}$,
and hence we obtain
a $C^{\infty}$-section
$-\rho_h(\theta)$ of
$\gbigp_G(\gminig)\otimes\Omega^{0,1}$.
Then, $h$ is called a harmonic metric of
the $G$-Higgs bundle $(\gbigp_G,\theta)$
if $R(h)+[\theta,-\rho_h(\theta)]=0$.

\subsubsection{Compatibility with an automorphism $\sigma$}

Let $\sigma$ be a complex algebraic automorphism of
$G$ of finite order $m>0$.
Let $G_0^{\sigma}$ denote the connected component of 
$G^{\sigma}$ containing $1$.
Let $H\subset G^{\sigma}$ be a closed subgroup
such that $G_0^{\sigma}\subset H$.

Let $\gbigp_{H}$
be a holomorphic principal $H$-bundle.
Let $\theta$ be a section of
$\gbigp_{H}(\gminig_1)\otimes\Omega^1$.
We obtain the holomorphic principal $G$-bundle
$\gbigp_{G}=\gbigp_{H}(G)$.
Because
$\gbigp_{H}(\gminig_1)
\subset
\gbigp_G(\gminig)$,
we obtain
the section $\theta$ of $\gbigp_G(\gminig)\otimes K_X$.
In this way,
we obtain a $G$-Higgs bundle
$(\gbigp_G,\theta)$.
We say a harmonic metric $h$ of $(\gbigp_G,\theta)$
is compatible with $\sigma$
if $\sigma(h)=h$
as a section of $\gbigp_G(\Herm(G))$,
i.e.,
$h$ is a section of 
$\gbigp_{H}(\Herm(G,\sigma))$.

Let
$\Harm(\gbigp,\theta,\sigma)$
denote the set of
harmonic metrics of $(\gbigp_G,\theta)$
compatible with $\sigma$.

\subsubsection{Specialized Simpson's main estimates}

We assume that $(\sigma,\omega)$ is split.
Assume that 
$\theta\in H^0(X,\gbigp_{H}(\gminig_1)\otimes \Omega_X^1)$
is regular semisimple in the following sense.
\begin{itemize}
 \item For a local expression $\theta=f\,dz$,
$f$ is a holomorphic section of
$\gbigp_{H}(\gminig_1^{\rs})$.
\end{itemize} 

\begin{prop}[Lemma
\ref{lem;25.5.13.30}]
There exists a unique
$h^{\can}\in\Harm(\gbigp_{G},\theta,\sigma)$ 
such that
it is decoupled
in the sense
$R(h^{\can})=[\theta,-\rho_{h^{\can}}(\theta)]=0$.
\end{prop}

As a result, we obtain the $K$-reduction
$\gbigp^{\can}_{K}:=
\gbigp^{h^{\can}}_{K}\subset \gbigp_G$
corresponding to $h^{\can}$.
The induced metric of
$\gbigp_G(\gminig)=
\gbigp^{\can}_{K}(\gminig)$
is also denoted by $h^{\can}$.

Let $W(\gminit)$ denote the Weyl group of $\gminit$.
We have the natural action of $W(\gminit)$ on
$\gminit\otimes\Omega_X^1$.
We obtain the quotient space
$(\gminit\otimes\Omega_X^1)/W(\gminit)$.
We obtain the section
$\Psi_{\theta}$
of $(\gminit\otimes\Omega_X^1)/W(\gminit)$.

For any  $h\in\Harm(\gbigp_G,\theta,\sigma)$,
we obtain
the section
$v(h^{\can},h)$ of
$\gbigp^{\can}_{K}(\sqrt{-1}\gminik)
\subset
 \gbigp_{G}(\gminig)$
determined by
$h^{\can}$ and $h$.
(See Lemma \ref{lem;25.5.2.3}.)
Let $W$ be any compact subset of $X$.
We obtain the following generalization
of Theorem \ref{thm;25.9.26.50}.

\begin{thm}[Theorem
\ref{thm;25.6.20.20}]
\label{thm;25.9.26.100}
There exist $A_{11},A_{12}>0$
depending only on
$(\Psi_{\theta},W,g_X)$
such that the following holds.
\begin{itemize}
 \item For any $t\geq 1$
       and $h\in \Harm(\gbigp_G,t\theta,\sigma)$,
       we obtain
       $\bigl|
       v(h^{\can},h)
       \bigr|_{h^{\can}}
       \leq
       A_{11}\exp(-A_{12}t)$
       on $W$.
\end{itemize}
\end{thm}

As a direct application,
we obtain the following theorem from 
Theorem \ref{thm;25.9.26.100}
by using an argument in \cite{Li-Mochizuki3}.

\begin{thm}[Theorem
\ref{thm;25.5.13.32}]
If $\theta\in H^0(X,\gbigp_{H}(\gminig_1\otimes\Omega_X^1))$
is generically regular semisimple,
then we have
$\Harm(\gbigp_G,\theta,\sigma)\neq\emptyset$.
 \end{thm}

\subsubsection{Harmonic metrics for generic cyclic $G$-Higgs bundles}

Let us consider the more detailed classification problem
in the setting of \S\ref{subsection;25.9.26.101}.
We have
$\gminig_1=\bigoplus_{i=1}^{\ell}\gminig_{\alpha_i}
\oplus\gminig_{-\psi}$.
Let $\gbigp_T$ be a holomorphic principal $T$-bundle.
We have the decomposition
$\gbigp_T(\gminig)
 =\bigoplus_{\phi\in\Delta}
 \gbigp_T(\gminig_{\phi})$.
Let
$\psi_i$ be the positive integers such that
$\psi=\sum_{i=1}^{\ell} \psi_i\alpha_i$.
We set
$\Upsilon(\gminig,\gminit):=
\bigotimes_{i=1}^{\ell}
\gminig_{\alpha_i}^{\otimes\psi_i}
\otimes\gminig_{-\psi}$.
Note that the induced $T$-representation on
$\Upsilon(\gminig,\gminit)$ is trivial.

Let $\theta\in H^0(X,\gbigp_T(\gminig_1)\otimes\Omega^1_X)$.
We obtain the decomposition
$\theta=\sum\theta_{\alpha_i}+\theta_{-\psi}$.
We obtain the $(\tth+1)$-differential
\[
 \gminio(\theta)
 = \prod_{i=1}^{\ell}
 \theta_{\alpha_i}^{\otimes\psi_i}
 \otimes
 \theta_{-\psi}
 \in
 H^0\Bigl(X,
 \Upsilon(\gminig,\gminit)
 \otimes K_X^{\otimes(\tth+1)}
 \Bigr).
\]

Let us consider the case $X=\Xbar\setminus D$,
where $\Xbar$ is a compact Riemann surface,
and $D$ is a finite subset of $\Xbar$.
We assume that
$\gminio(\theta)
\in H^0\bigl(\Xbar,\Upsilon(\gminig,
\gminit)\otimes K_{\Xbar}^{\otimes(\tth+1)}(\ast D)\bigr)$,
and that $\gminio(\theta)$ is not constantly $0$.
For each $P\in D$,
let $\ord_P(\gminio(\theta))$
denote the zero order of $\gminio(\theta)$ at $P$.

Let $D^{>0}$ denote the set of $P\in D$
such that
$\ord_P\gminio(\theta)+\tth+1>0$.
For $P\in D^{>0}$,
let 
$\nbigs_P(\theta)$ denote the set of
$\veca=(a_i)\in \real^{\ell}_{\leq 0}$
such that
$\sum_{i=1}^{\ell} a_{i}+\tth+1+\ord_P(\gminio(\theta))\geq 0$.
\begin{thm}[Theorem
\ref{thm;25.5.10.20}]
There exists a natural bijection
$\Psi:\Harm(\gbigp_G,\theta,\sigma)
\simeq
 \prod_{P\in D^{>0}}
 \nbigs_P(\theta)$.
Here, 
if $D_{>0}=\emptyset$,
then
$\prod_{P\in D^{>0}}
\nbigs_P(\theta)$
is defined to be the set of one point.
\end{thm}

\paragraph{Acknowledgement}

I thank Qiongling Li for discussions on many occasions.
This study partially grew out of our studies
on Toda equations, real harmonic bundles
and harmonic principal $G$-bundles.
I thank Oscar Garc\'{\i}a-Prada and Martin Guest
for discussions on cyclic Higgs bundles
and Toda equations.
I thank Nan-Kuo Ho for informing the work with Guest.

I gave talks on this study
in the workshops
``Higher rank geometric structures, 
Higgs bundles and physics''
held at the Institut Henri Poincar\'{e} in May 2025,
and
``Higgs bundles and harmonic maps''
held at the University of Kiel in July 2025.
The preparation of the talks at these conferences
were useful in this research.
I thank the organizers for the opportunities.

I am partially supported by
the Grant-in-Aid for Scientific Research (A) (No. 21H04429),
the Grant-in-Aid for Scientific Research (A) (No. 22H00094),
the Grant-in-Aid for Scientific Research (A) (No. 23H00083),
and the Grant-in-Aid for Scientific Research (C) (No. 20K03609),
and the Grant-in-Aid for Scientific Research (C) (No. 25K06973)
Japan Society for the Promotion of Science.
I am also partially supported by the Research Institute for Mathematical
Sciences, an International Joint Usage/Research Center located in Kyoto
University.

\section{Preliminary for complex semisimple Lie groups}

\subsection{Algebraricity}

Let $\gminig$ be a semisimple complex Lie algebra.
Let $G$ be a connected complex Lie group 
whose Lie algebra is $\gminig$.
\begin{prop}[\mbox{\cite[Corollary 7.6]{Knapp-Book},
\cite[Theorem VIII.10]{Serre-Lie-algebra}}]
\mbox{{}}\label{prop;25.6.17.20}
\begin{itemize}
 \item There is a unique complex algebraic structure on $G$
       compatible with its complex Lie group structure.
 \item $G$ is realized as a closed algebraic subgroup of
       $\GL(n,\cnum)$ for some $n$.
 \item If $G'$ is a complex algebraic group,
       then any homomorphism of complex Lie groups $G\to G'$ is algebraic.
       \hfill\qed
\end{itemize}
\end{prop}

\subsection{Maximal compact subgroups}

Let $K\subset G$ be a maximal compact subgroup.
Let $\gminik\subset\gminig$ denote the Lie algebra of $K$.
\begin{prop}[\mbox{\cite[Theorem 12]{Serre-Lie-algebra}}]
\mbox{{}}
\begin{itemize}
 \item $\gminik$ is a compact real form of $\gminig$.
 \item We obtain the global Cartan decomposition.
       Namely, 
\begin{equation}
\label{eq;25.5.2.2}
 K\times(\sqrt{-1}\gminik)
\lrarr
 G,
 \quad
 (k,v)
 \longmapsto
 k\cdot \exp(v)
\end{equation}
is an isomorphism of real analytic manifolds.       
In particular,
$K$ is connected because $G$ is assumed to be connected.
\hfill\qed
\end{itemize}
\end{prop}

Let $\rho_{\gminik}$ denote
the anti-$\cnum$-linear involution with respect to $\gminik$,
i.e.,
\[
 \rho_{\gminik}=\id_{\gminik}\oplus
 (-\id_{\sqrt{-1}\gminik}).
\]
It is an automorphism of real Lie algebras.
It induces an anti-holomorphic involution $\rho_K$ of $G$.
We have
$\rho_K(k\exp(v))=k\cdot \exp(-v)$
for $k\in K$ and $v\in\sqrt{-1}\gminik$.

\begin{lem}
\label{lem;25.4.29.20}
The normalizer of $K$ in $G$ equals $K$.
\end{lem}
\pf
Though this is well known,
we explain an outline of a proof.
Let $g\in G$ such that $gKg^{-1}=K$.
There exists $k\in K$ and $v\in \sqrt{-1}\gminik$
such that $g=\exp(v)k$.
We obtain $\exp(v)K\exp(-v)=K$.
For any $k_1\in K$,
there exists $k_2\in K$
such that
$\exp(v)k_1=k_2\exp(v)=\exp(\Ad(k_2)v)k_2$.
It implies
$k_1=k_2$ and $v=\Ad(k_1)v$.
We obtain $[\gminik,v]=0$.
It implies that $v$ is contained in
the center of $\gminig$.
Because $\gminig$ is semisimple,
we obtain $v=0$.
\hfill\qed

\vspace{.1in}
Let $\Herm(G)$ denote the set of maximal compact subgroups.
It is equipped with the left $G$-action defined by
$(g,K_1)\longmapsto gK_1g^{-1}$.
We have the map $G\to \Herm(G)$
given by $g\longmapsto gKg^{-1}$.
It induces a $G$-equivariant map
\begin{equation}
\label{eq;25.5.2.1}
 G/K\lrarr \Herm(G).
\end{equation}
\begin{lem}
The map {\rm(\ref{eq;25.5.2.1})} is a bijection.
\hfill\qed 
\end{lem}

\subsubsection{Difference between two compact real forms}

Let $K_i\subset G$ $(i=1,2)$ be maximal compact subgroups.
Let $\gminik_i\subset\gminig$ denote the corresponding Lie subalgebras.
There exists $g\in G$
such that $K_2=gK_1g^{-1}$.
\begin{lem}
\label{lem;25.5.7.1}
$s(K_1,K_2):=\rho_{K_1}(g)g^{-1}
 \in \exp(\sqrt{-1}\gminik_1)$
depends only on $(K_1,K_2)$.
Moreover,
$s(aK_1a^{-1},aK_2a^{-1})=
as(K_1,K_2)a^{-1}$ holds
for any $a\in G$.
\end{lem}
\pf
For another choice $g'$
such that $K_2=g'K_1(g')^{-1}$,
there exists $k\in K_1$ such that $g'=gk$.
We obtain
$\rho_{K_1}(g')(g')^{-1}
=\rho_{K_1}(g)kk^{-1}g^{-1}
=\rho_{K_1}(g)g^{-1}$.
For $a\in G$,
we have
$aK_2a^{-1}=(aga^{-1})(aK_1a^{-1})(aga^{-1})^{-1}$
and
$\rho_{aK_1a^{-1}}(aga^{-1})=a\rho_{K_1}(g)a^{-1}$.
Hence, we can check
$s(aK_1a^{-1},aK_2a^{-1})=as(K_1,K_2)a^{-1}$
directly.
\hfill\qed

\vspace{.1in}

There exists the Cartan decomposition $g=ab$,
where $a\in \exp(\sqrt{-1}\gminik_1)$
and $b\in K_1$.
Then, $s(K_1,K_2)=a^{-2}$.

\begin{lem}
We have
$\rho_{K_2}(x)=
\rho_{K_1}\bigl(s(K_1,K_2)xs(K_1,K_2)^{-1}\bigr)$
for any $x\in G$,
and
$\rho_{\gminik_2}(u)=\rho_{\gminik_1}\circ \Ad(s(K_1,K_2))(u)$
for any $u\in\gminig$.
\end{lem}
\pf
Because $K_2=gK_1g^{-1}$,
we have
\[
 \rho_{K_2}(x)
 =g\rho_{K_1}(g^{-1}xg)g^{-1}
 =\rho_{K_1}(\rho_{K_1}(g)g^{-1}x g\rho_{K_1}(g^{-1}))
=\rho_{K_1}\bigl(s(K_1,K_2)xs(K_1,K_2)^{-1}\bigr).
\]
Because
$\gminik_2=\Ad(g)(\gminik_1)$,
we obtain
$\rho_{\gminik_2}
=\Ad(g)\circ\rho_{\gminik_1}\circ\Ad(g^{-1})
=\rho_{\gminik_1}\circ
\Ad(\rho_{K_1}(g)g^{-1})
=\rho_{\gminik_1}\circ\Ad(s(K_1,K_2))$.
\hfill\qed

\begin{lem}
$s(K_2,K_1)=s(K_1,K_2)^{-1}$. 
\end{lem}
\pf
We have
\begin{multline}
s(K_2,K_1)
=\rho_{K_2}(g^{-1})g
=\rho_{K_1}
\bigl(s(K_1,K_2)g^{-1}s(K_1,K_2)^{-1}\bigr)g
=\rho_{K_1}\bigl(
 \rho_{K_1}(g)g^{-1}g^{-1}g\rho_{K_1}(g^{-1})\rho_{K_1}(g)
\bigr)
\\
 = g\rho_{K_1}(g^{-1})
=s(K_1,K_2)^{-1}.
\end{multline}
\hfill\qed

\begin{lem}
We have
$\rho_{K_j}(s(K_1,K_2)^{-1})=s(K_1,K_2)$ $(j=1,2)$.
\end{lem}
\pf
We have
\[
 \rho_{K_1}(s(K_1,K_2)^{-1})
=\rho_{K_1}\bigl(
g\rho_{K_1}(g)^{-1}
\bigr)
=\rho_{K_1}(g)g^{-1}=s(K_1,K_2).
\]
Because $s(K_1,K_2)^{-1}=s(K_2,K_1)$,
we also obtain
$\rho_{K_2}(s(K_1,K_2)^{-1})=s(K_1,K_2)$.
\hfill\qed

\begin{cor}
$s(K_1,K_2)\in\exp(\sqrt{-1}\gminik_1\cap\sqrt{-1}\gminik_2)$.
\hfill\qed
\end{cor}

\begin{lem}
\label{lem;25.5.2.3}
There exists a unique element
$v(K_1,K_2)\in \sqrt{-1}\gminik_1\cap\sqrt{-1}\gminik_2$
such that
\[
K_2=\exp(v(K_1,K_2))K_1\exp(-v(K_1,K_2)).
\]
\end{lem}
\pf
It is determined by
the condition $\exp(2v(K_1,K_2))=s(K_1,K_2)$.
\hfill\qed

\subsubsection{The induced Hermitian metrics on $\gminig$}

Let $B_{\gminig}$ denote the Killing form on $\gminig$,
i.e.,
$B_{\gminig}(u,v)=\Tr(\ad(u)\ad(v))$
for any $u,v\in\gminig$.
It is a $\cnum$-valued $G$-invariant non-degenerate
symmetric complex bilinear form.

For $K\in\Herm(G)$,
we obtain the $K$-invariant
positive definite Hermitian metric $h_{\gminig,K}$ of $\gminig$
as follows:
\[
 h_{\gminig,K}(u,v)
 =-B_{\gminig}(u,\rho_{\gminik}(v))
 \quad
 (u,v\in\gminig).
\]

\begin{lem}
\label{lem;25.4.29.2}
For $u,v,w\in\gminig$,
we have 
$h_{\gminig,K}(\ad(w)u,v)+h_{\gminig,K}(u,\ad(\rho_{\gminik}(w))v)=0$.
In other words,
the adjoint of $\ad(w)$ with respect to $h_{\gminig,K}$
equals 
$\ad(-\rho_{\gminik}(w))$,
i.e.,
$\ad(w)^{\dagger}_{h_{\gminig,K}}=\ad(-\rho_{\gminik}(w))$.
\hfill\qed
\end{lem}

\begin{lem}
For $u,v\in\gminig$ and $g\in G$,
we have
$h_{\gminig,K}(\Ad(g)u,v)=h_{\gminig,K}(u,\Ad(\rho_K(g)^{-1})v)$.
In other words,
$\Ad(g)^{\dagger}_{h_{\gminig,K}}
=\Ad(\rho_K(g)^{-1})$. 
\hfill\qed
\end{lem}

\begin{lem}
Let $K_1,K_2\in\Herm(G)$.
For any $u,v\in\gminig$,
we have
\begin{equation}
\label{eq;25.5.2.4}
 h_{\gminig,K_2}(u,v)
=h_{\gminig,K_1}\bigl(\Ad(s(K_1,K_2))u,v\bigr).
\end{equation}
\end{lem}
\pf
We have
\[
 h_{\gminig,K_2}(u,v)
=-B_{\gminig} (u,\rho_{\gminik_2}(v))
=-B_{\gminig}\Bigl(u,\rho_{\gminik_1}(\Ad(s(K_1,K_2))v)\Bigr)
=h_{\gminig,K_1}(u,\Ad(s(K_1,K_2))v).
\]
We also have
$\Ad(s(K_1,K_2))^{\dagger}_{h_{\gminig,K_1}}
=\Ad(\rho_{K_1}(s(K_1,K_2)^{-1}))
=\Ad(s(K_1,K_2))$.
Hence, we obtain (\ref{eq;25.5.2.4}).
\hfill\qed

\subsubsection{Norms}

Let $K$ be a maximal compact subgroup of $G$.

\begin{lem}
\label{lem;25.4.27.20}
\label{lem;25.4.27.31}
Let $u\in\gminig$.
\begin{itemize}
\item 
Let $\alpha_1,\ldots,\alpha_{\dim \gminig}$
be the eigenvalues of $\ad(u)\in\End_{\cnum}(\gminig)$
with multiplicity.
Then, the following holds:
\begin{equation}
\label{eq;25.5.2.10}
 h_{\gminig,K}(u,u)
\geq \sum_{i=1}^{\dim \gminig}|\alpha_i|^2.
\end{equation}
\item The equality holds
       if and only if
       $[u,\rho_{\gminik}(u)]=0$.
\end{itemize}
\end{lem}
\pf
We obtain (\ref{eq;25.5.2.10})
by the definition of $h_{\gminig,K}$.
The equality holds
if and only if
$[\ad(u),\ad(u)^{\dagger}_{h_{\gminig,K}}]=0$.
Because $\ad:\gminig\to\End_{\cnum}(\gminig)$ is injective,
we obtain the second claim.
\hfill\qed

\subsection{Cartan subalgebras and regular semisimple elements}

\subsubsection{Cartan subalgebras}

We recall the following general theorem for Cartan subalgebras.
\begin{prop}[\mbox{\cite[Theorem 13.1]{Nishiyama-Ohta}, \cite[Theorem VIII.1]{Serre-Lie-algebra}}]
\label{prop;23.2.16.1}
Let $\gminit\subset\gminig$ be any Cartan subalgebra.
Let $T$ denote the corresponding Lie subgroup of $G$.
\begin{itemize}
 \item $T$ is a closed algebraic subgroup of $G$,
       and it is a maximal complex torus in $G$,
       i.e.,
       $T\simeq (\cnum^{\ast})^{\dim_{\cnum}\gminit}$.
       In particular,
       $T$ has a unique maximal compact subgroup $T_c$
       which is a real torus.
       The real dimension of $T_c$
       equals
       $\dim_{\cnum}\gminit$.
 \item $T$ is the centralizer of $\gminit$ in $G$,
       i.e.,
       $a\in G$ is contained in $T$
       if and only if
       $\Ad(a)u=u$ for any $u\in \gminit$.
 \item $T$ is the centralizer of $T$ in $G$,
       i.e.,
       $a\in G$ is contained in $T$
       if and only if
       $aga^{-1}=g$ for any $g\in T$.
       \hfill\qed
\end{itemize}
\end{prop}

For any Cartan subalgebra $\gminit$,
we set
$N_G(\gminit):=\bigl\{
b\in G\,\big|\,
\Ad(b)\gminit=\gminit
\bigr\}$.
Let $W(\gminit)$ denote the Weyl group of $(\gminig,\gminit)$.
According to \cite[Lemma 8.8]{Matsushima-Lie-algebra},
the natural morphism
$N_{G}(\gminit)\to\GL(\gminit)$ induces
a surjection $N_G(\gminit)\to W(\gminit)$.
The kernel is $T$
by Proposition \ref{prop;23.2.16.1}.

\subsubsection{The real part and the purely imaginary part}
\label{subsection;25.6.15.20}

For any Cartan subalgebra $\gminit$,
we obtain the root system
$\Delta(\gminit,\gminig)\subset\gminit^{\lor}$.
Let $\gminit_{\real}\subset\gminit$
be the real subspace determined by the condition
$\phi(\gminit_{\real})\subset\real$
for any $\phi\in\Delta(\gminit,\gminig)$.
The Lie algebra of $T_c$ equals
$\sqrt{-1}\gminit_{\real}$.

\subsubsection{Regular semisimple elements}

For any $a\in\gminig$,
let $C_{\gminig}(a)$ denote the centralizer of $a$ in $\gminig$,
i.e.,
$C_{\gminig}(a)=\{b\in \gminig\,|\,[b,a]=0\}$.
The element $a$ is called regular
if
$\dim_{\cnum} C_{\gminig}(a)
=\min_{b\in\gminig}\dim_{\cnum} C_{\gminig}(b)$.
If $a$ is regular, $\dim_{\cnum} C_{\gminig}(a)=\dim_{\cnum}\gminit'$ holds
for any Cartan subalgebra $\gminit'$ of $\gminig$.

An element $a\in\gminig$ is called semisimple
if $\ad(a)\in\End_{\cnum}(\gminig)$ is semisimple.
It is equivalent to the condition that
$F(a)$ is semisimple for any representation $F$ of
the Lie algebra $\gminig$.
Any semisimple element is contained in
a Cartan subalgebra.
(See \cite[Corollary 20.6.3]{Tauvel-Yu}.)
\begin{lem}
For any regular semisimple element $a$,
$C_{\gminig}(a)$ is a Cartan subalgebra of $\gminig$. 
It is characterized as a unique Cartan subalgebra
containing $a$.
\end{lem}
\pf
Because $a$ is semisimple,
there exists a Cartan subalgebra $\gminit$
such that $a\in\gminit$.
By definition, we have $\gminit\subset C_{\gminig}(a)$.
By comparing the dimensions,
we obtain $\gminit=C_{\gminig}(a)$.
\hfill\qed

\vspace{.1in}
Let $\gminit$ be any Cartan subalgebra of $\gminig$.
Any element of $\gminit$ is semisimple by definition.
Let $\gminit^{\rs}$ denote the set of
regular semisimple elements in $\gminit$.
For each $\alpha\in\Delta(\gminit,\gminig)$,
we obtain the hypersurface $P_{\alpha}$ of $\gminit$
as the zero of $\alpha$.
Then, $H(\gminit)=\gminit\setminus \gminit^{\rs}$
equals
$\bigcup_{\alpha\in\Delta(\gminit,\gminig)}P_{\alpha}$.

\begin{lem}
\label{lem;25.6.17.1}
Let $u\in\gminit$.
If $g\cdot u=u$ for some $g\in W(\gminit)\setminus\{1\}$,
then $u$ is contained in $\gminit\setminus\gminit^{\rs}$.
As a result, the action of $W(\gminit)$ on $\gminit^{\rs}$ is free.
\end{lem}
\pf
For any $g\in W(\gminig)\setminus\{1\}$,
we set
$\EE_1(\gminit,g)
=\bigl\{v\in\gminit\,\big|\,
g\cdot v=v\bigr\}\subset \gminit$.
Because the action of $g$ is defined on $\gminit_{\real}$,
we have
$\EE_1(\gminit,g)=(\EE_1(\gminit,g)\cap\gminit_{\real})\otimes\cnum$.
By \cite[\S10.3 Lemma B]{Humphreys},
$\EE_1(\gminit,g)\cap\gminit_{\real}$
is contained in
$\bigcup_{\alpha\in\Delta(\gminit,\gminig)}
(P_{\alpha}\cap\gminit_{\real})$.
Then, we obtain the claim of the lemma.
\hfill\qed

\vspace{.1in}

Let $C_G(a)$ denote the centralizer of $a$ in $G$,
i.e.,
$C_G(a)=\bigl\{
b\in G\,\big|\,\Ad(b)(a)=a
\bigr\}$.

\begin{lem}
If $a$ is regular semisimple,
$C_G(a)$ equals the centralizer of
the Cartan subalgebra $C_{\gminig}(a)$,
i.e.,
\[
 C_G(a)
 =\bigl\{
 g\in G\,\big|\,
 \Ad(g)u=u\,\,\forall u\in C_{\gminig}(a)
 \bigr\}.
\]
In particular, $C_G(a)$ is a maximal torus.
\end{lem}
\pf
Let $g\in C_G(a)$.
Because $\gminit:=C_{\gminig}(a)$ is the unique Cartan subalgebra
containing $a$,
we obtain $\Ad(g)\gminit=\gminit$.
It implies $g\in N_G(\gminit)$.
By Lemma \ref{lem;25.6.17.1},
we obtain $g\in C_G(\gminit)$.
\hfill\qed

\subsection{Maximal compact subgroups and Cartan subalgebras}

\subsubsection{Compatibility of
maximal compact subgroups and Cartan subalgebras}

A maximal compact subgroup $K$ of $G$
is called compatible with a Cartan subalgebra $\gminit$
if
$\gminit\cap\gminik=\sqrt{-1}\gminit_{\real}$.
It is equivalent to the condition
$\rho_{\gminik}(\gminit)=\gminit$.
It is also equivalent to the condition that
$T_c\subset K$.
Let $\Herm(G,\gminit)$
denote the set of maximal compact subgroups of $G$
compatible with $\gminit$.
The following lemma is clear.
\begin{lem}
 $\Herm(G,\gminit)=
\{K\in \Herm(G)\,|\,T_c\subset K\}$ is not empty.
\hfill\qed
\end{lem}

\begin{lem}
\label{lem;25.6.14.10}
For $K\in\Herm(G)$,
the following conditions are equivalent.
\begin{itemize}
 \item $K$ is compatible with $\gminit$.
 \item $[u,\rho_{\gminik}(u)]=0$ holds
       for any $u\in \gminit$.
 \item $[u,\rho_{\gminik}(u)]=0$ holds
       for some $u\in \gminit^{\rs}$.
\end{itemize}
\end{lem}
\pf
Suppose $K\in\Herm(G,\gminit)$.
Because
$\gminit=(\gminik\cap\gminit)\oplus
 \sqrt{-1}(\gminik\cap\gminit)$,
any $u\in\gminit$
is expressed as
$u=u_1+\sqrt{-1}u_2$,
where $u_i\in \sqrt{-1}\gminit_{\real}\subset\gminik$.
We obtain
$\rho_{\gminik}(u)=u_1-\sqrt{-1}u_2$,
and
$[u,\rho_{\gminik}(u)]=0$.
Therefore, the first condition 
implies the second condition.
The second condition clearly implies the third.
Suppose $[u,\rho_{\gminik}(u)]=0$ holds
for some $u\in \gminit^{\rs}$.
We obtain $\rho_{\gminik}(u)\in\gminit$.
It implies $\rho_{\gminik}(\gminit)=\gminit$,
i.e., the first condition holds.
\hfill\qed

\subsubsection{Normalizers and the Weyl group}

Let $K\in\Herm(G,\gminit)$.
We set
\[
 N_{K}(\gminit)
=\bigl\{
 b\in K\,\big|\,
 \Ad(b)(\gminit)=\gminit
\bigr\}
=\bigl\{
 b\in K\,\big|\,
 \Ad(b)(\sqrt{-1}\gminit_{\real})=\sqrt{-1}\gminit_{\real}
\bigr\}.
\]
The natural morphism
$N_{K}(\gminit)
\to
\GL(\gminit)$
induces the surjection
$N_{K}(\gminit)\to W(\gminit)$,
and the kernel is $T_c$.
(See \cite[Theorem 4.54]{Knapp-Book}.)
We obtain the following lemma.

\begin{lem}
$N_K(\gminit)$ is a normal subgroup of
$N_G(\gminit)$,
and  
we have $N_G(\gminit)/N_K(\gminit)=
 T/T_c=\exp(\gminit_{\real})$.
\hfill\qed
\end{lem}

\begin{lem}
\label{lem;25.6.14.11}
Let $K_1,K_2\in \Herm(G,\gminit)$.
\begin{itemize}
 \item $s(K_1,K_2)\in\exp(\gminit_{\real})$.
 \item There exists $v\in\gminit_{\real}$
       such that
       $K_2=\exp(v)K_1\exp(-v)$.
\end{itemize}
\end{lem}
\pf
There exists $g\in G$
such that $K_2=gK_1g^{-1}$
and $T_c=gT_cg^{-1}$.
It implies that $g\in N_G(\gminit)$.
There exists $v\in \gminit_{\real}$
and $a\in N_{K_1}(\gminit)$
such that
$g=\exp(v)a$.
Then, we obtain the claim of the lemma.
\hfill\qed

\begin{cor}
\label{cor;25.5.3.1}
For any $K\in \Herm(G,\gminit)$,
we obtain the bijection
\[
  \gminit_{\real}\simeq \Herm(G,\gminit)
\]
induced by
$v\mapsto \exp(v)K\exp(-v)$.
\hfill\qed
\end{cor}

\section{Split automorphisms}

\label{section;25.5.13.10}

\subsection{Automorphisms of finite order}
\subsubsection{The associated periodically graded Lie algebras}
\label{subsection;25.6.18.10}

Let $G$ be a complex semisimple Lie group.
Let $\sigma$ be a complex algebraic automorphism of $G$
of finite order $m>0$,
i.e., $\sigma^m=\id_G$, $\sigma^{\ell}\neq \id_G$ $(\ell=1,\ldots,m-1)$.
We set $G^{\sigma}:=\{g\in G\,|\,\sigma(g)=g\}$.
Let $G_0^{\sigma}$ denote
the connected component of $G^{\sigma}$
containing the unit element.
We set
$G_Z^{\sigma}
:=\bigl\{
 g\in G\,\big|\,
 \Ad(\sigma(g))=\Ad(g)
 \bigr\}$.
We clearly have
$G_0^{\sigma}\subset G^{\sigma}\subset G^{\sigma}_Z$.
We obtain the induced action $\sigma$ on $\gminig$.
We have
$G_Z^{\sigma}
=\bigl\{
 g\in G\,\big|\,
 \Ad(g)\circ\sigma=\sigma\circ\Ad(g)
\bigr\}$.

Let $\omega$ be a primitive $m$-th root of $1$.
We obtain the eigen decomposition
\begin{equation}
\label{eq;25.5.13.11}
 (\gminig,\sigma)
 =\bigoplus_{\ell=0}^{m-1}
 (\gminig_{\ell},\omega^{\ell}\id_{\gminig_{\ell}}).
\end{equation}
We have
$[\gminig_{k},\gminig_{\ell}]
\subset\gminig_{k+\ell}$,
where the indexes $k,\ell$ are considered in $\seisuu/m\seisuu$.
Thus,
we obtain a periodically graded Lie algebra.

\subsubsection{Cartan subspaces and the restricted Weyl groups}

We recall some general results from \cite{Vinberg}.
The author also referenced
\cite[\S14]{Nishiyama-Ohta} which is a Japanese textbook.
A subspace $\gminia\subset\gminig_1$ 
is called a Cartan subspace if
$\gminia$ is maximally abelian subspace
and if any element of $\gminia$ is semisimple.

\begin{thm}[\mbox{\cite[Theorem 1]{Vinberg}}]
\label{thm;25.6.13.2}
For any two Cartan subspaces $\gminia_i$ $(i=1,2)$
in $\gminig_1$,
there exists $g\in G_0^{\sigma}$
such that $\Ad(g)\gminia_1=\gminia_2$.
\hfill\qed
 \end{thm}

Let $\gminia\subset\gminig_1$ 
be any Cartan subspace.
Let $H$ be a complex algebraic closed subgroup of $G$
such that
$G^{\sigma}_0\subset H\subset G^{\sigma}_Z$.
We consider
$N_{H}(\gminia)=\bigl\{
g\in H\,\big|\,\Ad(g)\gminia=\gminia
\bigr\}$
and
$C_H(\gminia)=\bigl\{
g\in N_H(\gminia)\,\big|\,
\Ad(g)_{|\gminia}=\id_{\gminia}
\bigr\}$.
We set $W_H(\gminia)=N_H(\gminia)/C_H(\gminia)$,
which naturally acts on $\gminia$.

\begin{prop}[\mbox{\cite[Proposition 8]{Vinberg}}]
$W_H(\gminia)$ is a finite group.
\hfill\qed
\end{prop}

\begin{thm}[\mbox{\cite[Theorem 7]{Vinberg}}]
\label{thm;25.5.13.3}
The natural morphism
$\cnum[\gminig_1]^{H}
 \to
 \cnum[\gminia]^{W_H(\gminia)}$ 
is an isomorphism of algebras. 
As a result, 
we obtain the isomorphism of algebraic varieties
$\Spec\cnum[\gminig_1]^{H}
\simeq
\Spec\cnum[\gminia]^{W_H(\gminia)}$.
\hfill\qed
\end{thm}

\subsubsection{Maximal compact subgroups compatible with $\sigma$}

We say that a maximal compact subgroup $K$ of $G$
is compatible with $\sigma$
if $\sigma(K)=K$.
Because $G$ is connected,
it is equivalent to the condition
$\sigma(\gminik)=\gminik$.
The $\sigma$-invariant part of $\gminik$
is denoted by $\gminik^{\sigma}$,
i.e., $\gminik^{\sigma}=\gminik\cap\gminig_0$.
Let $\Herm(G,\sigma)$ denote the set of
maximal compact subgroups compatible with $\sigma$.
If $K\in\Herm(G,\sigma)$,
then $gKg^{-1}\in\Herm(G,\sigma)$ for any $g\in G_Z^{\sigma}$.
Hence, $G_Z^{\sigma}$ naturally acts on $\Herm(G,\sigma)$.

\begin{lem}
Let $K_1,K_2\in\Herm(G,\sigma)$.
Then,
 $v(K_1,K_2)\in
 \sqrt{-1}(\gminik_1^{\sigma}\cap\gminik_2^{\sigma})$. 
\end{lem}
\pf
Because $\sigma(K_i)=K_i$,
we obtain
$\sigma(v(K_1,K_2))=
v(\sigma(K_1),\sigma(K_2))
=v(K_1,K_2)$.
\hfill\qed

\begin{cor}
\label{cor;25.6.17.10}
If $\Herm(G,\sigma)\neq\emptyset$,
then for any $K\in \Herm(G,\sigma)$
we obtain
the diffeomorphism
$\sqrt{-1}\gminik^{\sigma}
\simeq
\Herm(G,\sigma)$
given by
$v\mapsto \exp(v)K\exp(-v)$. 
\hfill\qed
\end{cor}

Let $H$ be a complex algebraic closed subgroup of $G$
such that $G_0^{\sigma}\subset H\subset G_Z^{\sigma}$.
For any $K\in\Herm(G,\sigma)$, we set $K_H=K\cap H$.

\begin{lem}
If $K\in\Herm(G,\sigma)$,
then the Cartan decomposition
$K\times \sqrt{-1}\gminik\simeq G$
induces
$K_H\times\sqrt{-1}\gminik^{\sigma}\simeq H$.
\end{lem}
\pf
Let $g=k\exp(v)\in G$
for $k\in K$ and $v\in\sqrt{-1}\gminik$.
First, let us consider the case $H=G_Z$.
Suppose $g\in G_Z$.
We have
$\Ad(g)=\Ad(k)\cdot \Ad(\exp(v))$.
Note that $\Ad(k)$ is anti-self-adjoint
and $\Ad\exp(v)$ is self-adjoint
with respect to $h_{\gminig,K}$.
We also note that $\sigma$ on $\gminig$ is unitary
with respect to $h_{\gminig,K}$.
Hence, we obtain that
$\sigma\circ\Ad(k)\circ\sigma^{-1}$
and
$\sigma\circ\Ad(\exp(v))\circ\sigma^{-1}$
are anti-self-adjoint
and self-adjoint with respect to $h_{\gminig,K}$.
By the uniqueness of Cartan decomposition,
we obtain
$\sigma\circ\Ad(k)\circ\sigma^{-1}=\Ad(k)$
and 
$\sigma\circ\Ad(\exp(v))\circ\sigma^{-1}=\Ad\exp(v)$
from $\sigma\circ \Ad(g)\circ\sigma^{-1}=\Ad(g)$.
We obtain $v\in\sqrt{-1}\gminik^{\sigma}$
and $k\in K\cap G_Z$.

Let $g\in H$.
We obtain the decomposition
$g=k\cdot \exp(v)$,
where $k\in K\cap G_Z$
and $v\in\sqrt{-1}\gminik^{\sigma}$.
Because $\exp(tv)\in G_0^{\sigma}$ for any $t\in\cnum$,
we obtain $k\in K_H$.
\hfill\qed

\subsubsection{Regular semisimple elements in $\gminig_1$}

Let $\gminig_1^{\rs}$ denote
the set of regular semisimple elements contained in $\gminig_1$.
Suppose that $\gminig_1^{\rs}\neq\emptyset$.

Let $u\in\gminig_1^{\rs}$.
Because
$\sigma(u)=\omega u$,
we obtain
$\sigma(C_{\gminig}(u))=C_{\gminig}(u)$.
By setting
$C_{\gminig}(u)_{\ell}=
C_{\gminig}(u)\cap \gminig_{\ell}$,
we obtain the decomposition
\[
 C_{\gminig}(u)
 =\bigoplus_{\ell=0}^{m-1}
 C_{\gminig}(u)_{\ell}.
\]

\begin{lem}
\label{lem;25.6.13.1}
$C_{\gminig}(u)_1$
is the unique Cartan subspace containing $u$.
\end{lem}
\pf
Because $C_{\gminig}(u)$ is a Cartan subalgebra of $\gminig$,
$C_{\gminig}(u)_1$
is an abelian subspace,
and any element of $C_{\gminig}(u)_1$ is semisimple.
Let $\gminia'\subset\gminig_1$ be an abelian subalgebra
such that $C_{\gminig}(u)_1\subset\gminia'$.
Because any element of $\gminia'$ is commuting with $u$,
we obtain $\gminia'\subset C_{\gminig}(u)$,
and hence $\gminia'=C_{\gminig}(u)_1$,
i.e., $C_{\gminig}(u)_1$ is a Cartan subspace.
If $\gminia\subset\gminig_1$ is a Cartan subspace
containing $u$,
we obtain
$\gminia\subset C_{\gminig}(u)$,
and hence
$\gminia\subset C_{\gminig}(u)_1$.
\hfill\qed

\begin{cor}
Under the assumption $\gminig_1^{\rs}\neq\emptyset$,
for any Cartan subspace $\gminia\subset\gminig_1$,
$\gminia^{\rs}:=\gminia\cap\gminig_1^{\rs}$ is not empty.
\end{cor}
\pf
By Lemma \ref{lem;25.6.13.1},
there exists a Cartan subspace $\gminia$
such that $\gminia^{\rs}\neq\emptyset$.
By Theorem \ref{thm;25.6.13.2},
we obtain $\gminia^{\rs}\neq\emptyset$
for any Cartan subspace $\gminia$.
\hfill\qed

\vspace{.1in}
Let $H$ be a complex algebraic closed subgroup of $G$
such that $G_0^{\sigma}\subset H\subset G_Z^{\sigma}$.
For any $K\in\Herm(G,\sigma)$, we set $K_H=K\cap H$.
Let $\gminia$ be a Cartan subspace.

\begin{lem}
\label{lem;25.6.17.2}
For any $u\in \gminia^{\rs}$,
we have
$C_H(\gminia)=C_H(u)$.
\end{lem}
\pf
Clearly, we have $C_H(\gminia)\subset C_H(u)$.
Let $g\in C_H(u)$.
Because $u$ is regular semisimple,
we have $g\in C_G(C_{\gminia}(u))$.
It implies $g\in C_H(\gminia)$.
\hfill\qed

\begin{lem}
For $a\in \gminig_1^{\rs}$,
we set
$M_H(a,\gminia)=
\bigl\{
 g\in H\,\big|\,
 \Ad(g)a\in\gminia
\bigr\}$.
\begin{itemize}
 \item $M_H(a,\gminia)=\bigl\{
 g\in H\,\big|\,
 \Ad(g)C_{\gminig}(a)_1=\gminia
 \bigr\}$.
 \item $M_H(a,\gminia)$ is
       non-empty.
 \item The natural left action of
       $N_H(\gminia)$ on $M_H(a,\gminia)$
       is free and transitive.
\end{itemize}
\end{lem}
\pf
Let us study the first claim.
The implication $\supset$ is clear.
Let $g\in M_H(a,\gminia)$.
Because $\Ad(g)a\in \gminia^{\rs}$,
we have
$\Ad(g)C_{\gminig}(a)_1
=C_{\gminig}(\Ad(g)a)_1
=\gminia$.
Thus, we obtain the first claim.
We obtain the second claim by Theorem \ref{thm;25.6.13.2}.
Let $g_1,g_2\in M_H(a,\gminia)$.
We have $\Ad(g_2)a\in\gminia^{\rs}$
and 
$\Ad(g_1\cdot g_2^{-1})(\Ad(g_2)a)\in\gminia^{rs}$.
We obtain
$g_1\cdot g_2^{-1}\in N_H(\gminia)$.
Hence, we obtain the third claim.
\hfill\qed

\begin{lem}
The $W_H(\gminia)$-action on $\gminia^{\rs}$ is free.
\end{lem}
\pf
Let $u\in \gminia^{\rs}$.
Suppose $\Ad(g)u=u$ for some $g\in H$.
Because $u$ is regular semisimple,
we obtain
$\Ad(g)v=v$ for any $v\in C_{\gminig}(u)$.
Because $\gminia=C_{\gminig}(u)_1$,
we obtain $g\in C_H(\gminia)$.
\hfill\qed

\subsubsection{Normality and diagonalizability}

We continue to assume $\gminig_1^{\rs}\neq\emptyset$.
Let $\gminia\subset\gminig_1$ be a Cartan subspace.
Let $K$ be a maximal compact subgroup of $G$
such that
$K\in\Herm(G,\sigma)\cap\Herm(G,C_{\gminig}(\gminia))$.
We set
$K_0^{\sigma}=K\cap G_0^{\sigma}$.
\begin{prop}
\label{prop;25.6.15.1}
The following conditions are equivalent for $u\in\gminig_1^{\rs}$.
\begin{itemize}
 \item $[u,\rho_{\gminik}(u)]=0$.
 \item There exists $g\in K_0^{\sigma}$
       such that $\Ad(g)u\in\gminia$.
\end{itemize} 
\end{prop}
\pf
If there exists $g\in K^{\sigma}_0$
such that $\Ad(g)(u)\in\gminia$,
we obtain
\[
 [u,\rho_{\gminik}(u)]
 =
  \Ad(g)^{-1}\Bigl[
 \Ad(g)(u),
 \rho_{\gminik}(\Ad(g)(u))
 \Bigr]=0.
\]

Suppose that $[u,\rho_{\gminik}(u)]=0$.
Let $T,T(u)\subset G$ denote the tori corresponding to
$C_{\gminig}(\gminia)$
and $C_{\gminig}(u)$,
respectively.
Let $T_c$ and $T(u)_c$ denote the maximal compact subgroups of
$T$ and $T(u)$, respectively.
Because $[u,\rho_{\gminik}(u)]=0$,
we obtain
$\rho_{\gminik}(C_{\gminig}(u))=C_{\gminig}(u)$.
It implies $T(u)_c\subset K_0^{\sigma}$.
Because $K\in \Herm(G,C_{\gminig}(\gminia))$,
we also have $T_c\subset K_0^{\sigma}$.
Both $T_c$ and $T(u)_c$ are maximal tori of $K_0^{\sigma}$,
there exists $k\in K_0^{\sigma}$
such that $kT(u)_ck^{-1}=T_c$.
Then, we obtain
$\Ad(k)C_{\gminig}(u)=C_{\gminig}(\gminia)$.
Because $\Ad(k)$ preserves the eigen decomposition of $\sigma$,
we obtain
$\Ad(k)C_{\gminig}(u)_1=\gminia$.
It implies
$\Ad(k)u\in\gminia$.
\hfill\qed

\vspace{.1in}
Let $H\subset G$
be any algebraically closed subgroup such that
$G^{\sigma}_0\subset H\subset G^{\sigma}_Z$.
\begin{prop}
\label{prop;25.6.15.10}
Let $g\in H$.
Let $u\in\gminig_1^{\rs}$.
Suppose that
$[u,\rho_{\gminik}(u)]=0$ 
and
$[\Ad(g)u,\rho_{\gminik}(\Ad(g)u)]=0$.
 Then,
\[
 g\in K_H\exp\Bigl(
 \sqrt{-1}\gminik^{\sigma}
 \cap
 C_{\gminig}(u)_{\real}
 \Bigr).
\] 
\end{prop}
\pf
We express $g=k\exp(v)$
for $k\in K_H$ and $v\in\sqrt{-1}\gminik^{\sigma}$.
We obtain
$\bigl[
\Ad(\exp(v))u,\,
\Ad(\exp(-v))
\rho_{\gminik}(u)
\bigr]
=0$.
It implies
$\bigl[
 \Ad(\exp(2v))u,
 \rho_{\gminik}(u)
 \bigr]=0$.
Because $u\in\gminig_1^{\rs}$ and $[u,\rho_{\gminik}(u)]=0$,
we have
$C_{\gminig}(u)=C_{\gminig}(\rho_{\gminik}(u))
=\rho_{\gminik}(C_{\gminig}(u))$.
We obtain
$\Ad(\exp(2v))C_{\gminig}(u)=C_{\gminig}(u)$,
i.e.,
$\exp(2v)\in N_G(C_{\gminig}(u))$.
There exists $k\in N_K(C_{\gminig}(u))$
and $w\in C_{\gminig}(u)_{\real}$
such that
$\exp(2v)=k\exp(w)$.
By the uniqueness of the Cartan decomposition,
we obtain
$k=1$ and $2v=w$.
Hence,
$v\in  \sqrt{-1}\gminik^{\sigma}
 \cap
 C_{\gminig}(u)_{\real}$.
\hfill\qed

\subsection{Split automorphisms}

We introduce the following condition.
\begin{condition}\mbox{{}}
\label{condition;25.5.13.2}
We say that $(\sigma,\omega)$ is split if
$\gminig^{\rs}_1\neq\emptyset$,
and if 
$C_{\gminig}(u)^{\sigma}:=C_{\gminig}(u)\cap\gminig_0=0$
for any $u\in\gminig^{\rs}_1$.
\hfill\qed
\end{condition}

\begin{lem}
Suppose that $\gminig^{\rs}_1\neq\emptyset$.
Then, $(\sigma,\omega)$ is split if and only if 
one of the following equivalent conditions hold.
\begin{itemize}
 \item There exists an element $u_0\in\gminig_1^{\rs}$
such that 
$C_{\gminig}(u_0)^{\sigma}=0$.
 \item For any Cartan subspace $\gminia\subset\gminig_1$,
       we obtain $C_{\gminig}(\gminia)^{\sigma}=0$.
 \item There exists a Cartan subspace $\gminia\subset\gminig_1$
       such that
       $C_{\gminig}(\gminia)^{\sigma}=0$.
\end{itemize}
\end{lem}
\pf
For any $u\in\gminig^{\rs}_1$,
we obtain the decomposition
$C_{\gminig}(u)
 =\bigoplus_{\ell=0}^{m-1}
 C_{\gminig}(u)_{\ell}$.
The dimensions are independent of $u\in \gminig_1^{\rs}$.
Hence, $(\sigma,\omega)$ is split
if and only if the first condition holds.
For any Cartan subspace $\gminia$ and any $u\in\gminia^{\rs}$,
we have $C_{\gminia}(u)=C_{\gminia}(\gminia)$.
Hence, 
both the second and third conditions are equivalent to
the condition that $(\sigma,\omega)$ is split.
\hfill\qed

\vspace{.1in}
For any $u\in\gminig_1$,
we set $C_H(u)=\{g\in H\,|\,\Ad(g)u=u\}$.
The following lemma is clear by the definition.
\begin{lem}
If $(\sigma,\omega)$ is split,
$C_{H}(u)$ is a finite group
for any $u\in\gminig_1^{\rs}$.
\hfill\qed
\end{lem}

\subsubsection{Finiteness of some groups}

Let $H\subset G^{\sigma}_Z$
be any algebraically closed subgroup such that $G^{\sigma}_0\subset H$.
Let $\gminia$ be any Cartan subspace.

\begin{lem}
$N_H(\gminia)$ and $C_H(\gminia)$ 
are finite groups.
\end{lem}
\pf
Because $C_H(\gminia)=C_H(u)$
for any $u\in\gminia^{\rs}$ by Lemma \ref{lem;25.6.17.2},
it is finite.
Because $W_H(\gminia)$ is finite,
$N_H(\gminia)$ is also finite.
\hfill\qed

\begin{cor}
\label{cor;25.6.14.1}
For any $a\in\gminig_1^{\rs}$,
$M_H(a,\gminia)$ is finite,
and
we have $|M_H(a,\gminia)|=|N_H(\gminia)|$.
\hfill\qed
\end{cor}

\subsubsection{Finite covering map and properness}

Let us consider the $H$-equivariant holomorphic map
\[
 \nu:H\times\gminia^{\rs}
 \to
 \gminig_1^{\rs},
 \quad
 \nu(g,u)
 =\Ad(g)u.
\]

\begin{prop}
\label{prop;25.5.13.31}
$\nu$ is a finite covering map.
The number of any fiber is
$|N_H(\gminia)|$.
\end{prop}
\pf
By Corollary \ref{cor;25.6.14.1},
the number of any fiber is $|N_H(\gminia)|$.
In particular, it is surjective.
By Sard's theorem, we obtain
$\dim (H\times\gminia^{\rs})
=\dim\gminig_1^{\rs}$.

For any $u\in\gminia^{\rs}$,
the morphism
$\ad(u):\gminig_0\to \gminig_1$
is injective.
It implies that
the derivative of $\nu$ at $(1,u)$ is injective
for any $u\in\gminia^{\rs}$,
where $1\in H$ denotes the unit element.
By comparing the dimensions,
we obtain that the derivative of $\nu$ at $(1,u)$ is an isomorphism
for any $u\in\gminia^{\rs}$.
Because $\nu$ is $H$-equivariant,
we obtain that $\nu$ is locally bi-holomorphic.
It implies that $\nu$ is a local homeomorphism.
Because the numbers of the fibers are constant,
$\nu$ is a covering map.
\hfill\qed

\begin{prop}
\label{prop;25.6.14.20}
Let $\nbiga$ be any compact subset in $\gminia^{\rs}$.
Then, 
$\nu(H\times\nbiga)$ is closed in $\gminig_1$.
As a result, 
the induced map
$H\times\nbiga\to \gminig_1$
is proper.
\end{prop}
\pf
The inclusion $\gminia^{\rs}\subset \gminig_1^{\rs}$
induces an isomorphism
$\gminia^{\rs}/W_H(\gminia)\simeq
\gminig_1^{\rs}/H$.
In particular,
there exists the natural algebraic morphism
$\gminig_1^{\rs}\to \gminia^{\rs}/W_H(\gminia)$.
By Theorem \ref{thm;25.5.13.3},
it extends to a morphism of complex algebraic varieties
$\varphi:\gminig_1\to \gminia/W_H(\gminia)$.
Moreover,
$\varphi^{-1}\Bigl(
 \gminia^{\rs}/W_H(\gminia)
 \Bigr)
=\gminig_1^{\rs}/H$.

Let $\pi_1:\gminia^{\rs}\to \gminia^{\rs}/W_H(\gminia)$
denote the projection.
We obtain a compact subset
$\pi_1(\nbiga)\subset \gminia^{\rs}/W_H(\gminia)$.
It is closed in
$\gminia/W_H(\gminia)$.
Hence, 
$\nu(H\times\nbiga)=
\varphi^{-1}(\pi_1(\nbiga))$
is closed in $\gminig_1$.

Because
$\nu$ is proper,
$H\times\nbiga\to \nu(H\times\nbiga)$
is proper.
Because $\nu(H\times\nbiga)$ is closed in $\gminig_1$,
the inclusion map
$\nu(H\times\nbiga)\to\gminig_1$ is proper.
Hence,
the induced map
$H\times\nbiga\to\gminig_1$ is proper.
\hfill\qed

\begin{cor}
\label{cor;25.9.20.1}
Let $\nbigb\subset\gminig_1^{\rs}$ be any compact subset.
Let $\nu_1:H\times\nbigb\to \gminig_1$
be defined by $\nu_1(g,u)=\Ad(g)(u)$.
Then, $\nu_1$ is proper.
\hfill\qed
\end{cor}

\subsubsection{Canonical maximal compact subgroups}
\label{subsection;25.6.18.11}

\begin{prop}
\label{prop;25.5.13.1}
For any abelian subspace $\gminia$,
$\Herm(G,\sigma)\cap\Herm(G,C_{\gminig}(\gminia))$
consists of a unique maximal compact subgroup
$K(\gminia)$.
In particular,
$\Herm(G,\sigma)\cap\Herm(G,C_{\gminig}(\gminia))\neq \emptyset$.
\end{prop}
\pf
We set
$\gminit=C_{\gminig}(\gminia)$.
It is a Cartan subalgebra of $\gminig$.
We have $\sigma(\gminit)=\gminit$ because
$\sigma(\gminia)=\gminia$.
Because $(\sigma,\omega)$ is split,
we have $\gminit^{\sigma}=\{0\}$.
Let $T\subset G$ be the maximal complex torus.
Let $T_c\subset T$ be the maximal compact subgroup of $T$.
We have
$\sigma(T)=T$ and $\sigma(T_c)=T_c$.

\begin{lem}
\label{lem;25.6.14.2}
The morphisms
$\id-\sigma:\gminit\to\gminit$ and
$\id-\sigma:\gminit_{\real}\to\gminit_{\real}$
are isomorphisms.
The homomorphisms
$\id-\sigma:T\to T$
and $\id-\sigma:T_c\to T_c$
are surjective. 
\end{lem}
\pf
It follows from $\gminit^{\sigma}=0$.
\hfill\qed

\vspace{.1in}
Let us prove Proposition \ref{prop;25.5.13.1}.
Let $K$ be a maximal compact subgroup of $G$
compatible with $\gminit$,
i.e.,
$T_c\subset K$.
Because $\sigma(T_c)=T_c$,
we have $T_c\subset \sigma(K)$.
Hence, $\sigma(K)$ is also compatible with $\gminit$.
By Lemma \ref{lem;25.6.14.11},
there exists $v\in \gminit_{\real}$
such that
$\sigma(K)=\exp(v)K\exp(-v)$.
By Lemma \ref{lem;25.6.14.2},
there exists $v_1\in \gminit_{\real}$
such that 
$v_1-\sigma(v_1)=v$.
Note that $v_1$ and $\sigma(v_1)\in\gminit_{\real}$
are commuting.
We set
$K_1=\exp(v_1)K\exp(-v_1)$.
We obtain $\sigma(K_1)=K_1$.
Because
$T_c=\exp(v_1)T_c\exp(-v_1)\subset K_1$,
$K_1$ is compatible with $\gminit$.

Let $K_2\in \Herm(G,\sigma)\cap\Herm(G,\gminit)$.
There exists $v\in \gminit_{\real}$
such that
$K_2=\exp(v)K_1\exp(-v)$
by Lemma \ref{lem;25.6.14.11}.
Because $\sigma(K_i)=K_i$,
we obtain $\sigma(v)=v$.
It implies $v=0$.
\hfill\qed

\vspace{.1in}

\begin{cor}
\label{cor;25.6.21.10}
For any $u\in \gminig_1^{\rs}$,
there exists a unique 
$K(u)=K(C_{\gminig}(u))\in\Herm(G,\sigma)$
such that the following condition holds.
\begin{itemize}
 \item $K(u)$ is compatible with
       the Cartan subalgebra $C_{\gminig}(u)$.
       It is equivalent to the condition
       $[u,\rho_{\gminik(u)}(u)]=0$
       according to
       Lemma {\rm\ref{lem;25.6.14.10}}.
\end{itemize} 
Clearly, 
we have $K(u')=K(u)$
for any $u'\in \gminig_1^{\rs}\cap C_{\gminig}(u)$.
 \hfill\qed
\end{cor}

\begin{cor}
For any $K\in \Herm(G,\sigma)$,
we obtain
the diffeomorphism
$\sqrt{-1}\gminik^{\sigma}
\simeq
\Herm(G,\sigma)$
given by
$v\mapsto \exp(v)K\exp(-v)$. 
\end{cor}
\pf
Because $\Herm(G,\sigma)\neq\emptyset$,
the claim follows from Corollary \ref{cor;25.6.17.10}.
\hfill\qed

\vspace{.1in}
For any Cartan subspace $\gminia\subset\gminig_1$,
let $K(\gminia)\in\Herm(G,\sigma)$
be determined by the condition
$K(\gminia)=K(u)$ for any $u\in\gminia^{\rs}$.
The following lemma is clear by the characterization of
$K(u)$.
\begin{lem}
\label{lem;25.9.27.1}
For any $g\in G_Z$,
we have $K(\Ad(g)u)=gK(u)g^{-1}$
and $K(\Ad(g)\gminia)=gK(\gminia)g^{-1}$.
\hfill\qed
\end{lem}

\begin{cor}
\label{cor;25.9.27.2}
The induced map
$K^{\can}:\gminig^{\rs}\to \Herm(G,\sigma)$
is $C^{\infty}$.
\end{cor}
\pf
It follows from Proposition \ref{prop;25.5.13.31}
and Lemma \ref{lem;25.9.27.1}.
\hfill\qed

\subsubsection{Normality}

Let $\gminia$ be a Cartan subspace.
We set $K=K(\gminia)$.
Let $H\subset G$ be a closed complex algebraic subgroup
such that $G_0^{\sigma}\subset H\subset G^{\sigma}_Z$.
We set $K_H=K\cap H$.

\begin{prop}
\label{prop;25.6.14.21}
Let $u\in\gminig_1^{\rs}$.
\begin{itemize}
 \item There exists $g\in G_0^{\sigma}$ such that
       $[\Ad(g)u,\rho_{\gminik}(\Ad(g)u)]=0$.
 \item If $[\Ad(g_i)u,\rho_{\gminik}(\Ad(g_i)u)]=0$
       for $g_i\in H$ $(i=1,2)$,
       then $g_1g_2^{-1}\in K_H$.
\end{itemize}
\end{prop}
\pf
We have the abelian subspace $C_{\gminig}(u)_1$
containing $u$.
There exists $g\in G^{\sigma}_0$
such that
$\Ad(g)C_{\gminig}(u)_1=\gminia$.
Then, we have
$[\Ad(g)u,\rho_{\gminik}(\Ad(g)u)]=0$.

Suppose that 
$[\Ad(g_i)u,\rho_{\gminik}(\Ad(g_i)u)]=0$
for $g_i\in H$.
We may assume $(g_1,g_2)=(g,1)$,
i.e.,
$[u,\rho_{\gminik}(u)]=0$
and
$[\Ad(g)u,\rho_{\gminik}(\Ad(g)u)]=0$
for some $g\in G_0^{\sigma}$.
By Proposition \ref{prop;25.6.15.1},
there exist $a_i\in K_0^{\sigma}$ $(i=1,2)$
such that
$\Ad(a_1)u\in\gminia$
and
$\Ad(a_2g)u\in\gminia$.
By Proposition \ref{prop;25.6.15.10}
and 
$ \sqrt{-1}\gminik^{\sigma}
 \cap
 C_{\gminig}(u)_{\real}
 \subset
 C_{\gminig}(u)^{\sigma}=0$,
$a_2ga_1^{-1}$ is contained in $K_H$.
Hence, we obtain $g\in K_ H$.
\hfill\qed

\begin{cor}
Let $u\in \gminia^{\rs}$
such that $[u,\rho_{\gminik}(u)]=0$.
Let $g\in H$.
Then, we have $[\Ad(g)u,\rho_{\gminik}(\Ad(g)u)]=0$
if and only if
$g\in K_H$.
\hfill\qed
\end{cor}

\subsection{Some estimates}

We continue to assume that
$(\sigma,\omega)$ is split.

\subsubsection{An estimate from the properness}

\begin{prop}
\label{prop;25.5.13.20}
Let $\nbigb$ be any compact subset in $\gminig_1^{\rs}$.
For any $C_1>0$, there exists $C_2>0$ such that
the following holds for any $u\in\nbigb$.
\begin{itemize}
\item If $K\in \Herm(G,\sigma)$ satisfies
      $|u|_{h_{\gminig,K}}\leq C_1$,
      then 
       $|v(K(u),K)|_{h_{\gminig,K(u)}}\leq C_2$ holds.
\end{itemize}
\end{prop}
\pf
Let $\gminia\subset\gminig_1$ be a Cartan subspace.
By Proposition \ref{prop;25.5.13.31},
it is enough to consider the case
$\nbigb\subset \gminia^{\rs}$.
By Proposition \ref{prop;25.5.13.1},
there exists unique $K(\gminia)\in\Herm(G,\sigma)$ such that
$[u,\rho_{\gminik(\gminia)}(u)]=0$ for any $u\in \gminia$.
For any $C_3>0$,
we consider
\[
 \nbigc_1=
 \Bigl\{
 (g,u)\in H\times \nbigb\,\Big|\,
 |\Ad(g)u|_{h_{\gminig,K(\gminia)}}
 \leq C_3
 \Bigr\}.
\]
Because
$H\times\nbigb\to \gminig_1$ is proper
by Corollary \ref{cor;25.9.20.1},
$\nbigc_1$ is compact.
Hence,
there exists $C_4>0$
such that
\[
\bigl|v(K(u),g^{-1}K(u)g)
\bigr|_{h_{\gminig,K(u)}}
=
\bigl|
v(K(\gminia),g^{-1}K(\gminia)g)
\bigr|_{h_{\gminig,K(\gminia)}}\leq C_4
\]
for any $(g,u)\in\nbigc_1$.
We also not that
$|\Ad(g)u|_{h_{\gminig,K(\gminia)}}
=|u|_{h_{\gminig},g^{-1}K(u)g}$
for $(g,u)\in\nbigc_1$.
Hence, we obtain the claim of Proposition \ref{prop;25.5.13.20}.
\hfill\qed

\subsubsection{Defect of the normality}

Let $\gminia$ be a Cartan subspace of $\gminig_1$.
Let $\nbiga\subset\gminia^{\rs}$ be a compact subset.
We consider the map
$F:G_0^{\sigma}\times\nbiga\to \real_{\geq 0}$
defined by
\[
 F(g,u)
 =h_{\gminig,K(\gminia)}
 \bigl(\Ad(g)\cdot u,\Ad(g)\cdot u\bigr)-h_{\gminig,K(\gminia)}(u,u).
\]
By Proposition \ref{prop;25.6.14.21}
and Proposition \ref{prop;25.5.13.20},
we obtain the following lemma.
\begin{lem}
\label{lem;25.6.20.10}
$F$ is proper.
Moreover,
$F^{-1}(0)=K(\gminia)_0^{\sigma}\times \nbiga$.
\hfill\qed
\end{lem}

Let 
$F_0:\sqrt{-1}\gminik(\gminia)^{\sigma}\times\nbiga\to \real_{\geq 0}$
be defined by
$F_0(v,u)=F(\exp(v), u)$.

\begin{prop}
For any $r_1>0$, there exists $C_1>0$
such that the following holds.
\begin{itemize}
 \item On $\bigl\{
       (v,u)\in\sqrt{-1}\gminik(\gminia)^{\sigma}\times\nbiga\,\big|\,
       F_0(v,u)\leq r_1
       \bigr\}$,
       we have
\[
       |v|_{h_{\gminig,K(\gminia)}}
       \leq
       C_1
       \Bigl|
       \bigl[u,\rho_{\gminik(\gminia)}(\Ad(\exp(2v))u)\bigr]
       \Bigr|_{h_{\gminig,K(\gminia)}}.
\]       
\end{itemize} 
\end{prop}
\pf
Let $v\in \sqrt{-1}\gminik(\gminia)^{\sigma}$ and $u\in\nbiga$.
We have
$\rho_{\gminik(\gminia)}\bigl(
\bigl[
 -\rho_{\gminik(\gminia)}(u),[u,v]
 \bigr]
 \bigr)
 =\bigl[
 u,[\rho_{\gminik(\gminia)}(u),v]
 \bigr]$
and
$\bigl[
 \rho_{\gminik(\gminia)}(u),u
\bigr]=0$.
By the Jacobi identity
$\bigl[
 u,[\rho_{\gminik(\gminia)}(u),v]
 \bigr]
 +\bigl[
 v,[u,\rho_{\gminik(\gminia)}(u)]
 \bigr]
 +\bigl[
 \rho_{\gminik(\gminia)}(u),[v,u]
 \bigr]=0$,
we obtain
$\rho_{\gminik(\gminia)}\bigl(
 \bigl[-\rho_{\gminik(\gminia)}(u),[u,v]\bigr]
 \bigr)
=-\bigl[-\rho_{\gminik(\gminia)}(u),[u,v]\bigr]$,
i.e.,
$\bigl[-\rho_{\gminik(\gminia)}(u),[u,v]\bigr]
\in\sqrt{-1}\gminik(\gminia)$.
Because $\sigma(K(\gminia))=K(\gminia)$,
we obtain
$\sigma\circ\rho_{\gminik(\gminia)}
=\rho_{\gminik(\gminia)}\circ\sigma$.
Hence, for $u\in\nbiga\subset\gminig_1$,
we have
$\sigma(\rho_{\gminik(\gminia)}(u))=
\rho_{\gminik(\gminia)}(\sigma (u))
=\rho_{\gminik(\gminia)}(\omega u)
=\omega^{-1} \rho_{\gminik(\gminia)}(u)$.
Therefore,
\[
\bigl[-\rho_{\gminik(\gminia)}(u),[u,v]\bigr]
\in\sqrt{-1}\gminik(\gminia)^{\sigma}.
\]

We obtain the linear maps
$\gamma_u:
\sqrt{-1}\gminik(\gminia)^{\sigma}
\to
\sqrt{-1}\gminik(\gminia)^{\sigma}$
$(u\in\nbiga)$
defined by
$\gamma_u(v)
=\bigl[
 -\rho_{\gminik(\gminia)}(u),[u,v]
\bigr]$.
Because
$h_{\gminig,K(\gminia)}\bigl(
 \gamma_u(v),
 v
 \bigr)
 =h_{\gminig,K(\gminia)}\bigl(
 \ad(u)v,\ad(u)v
 \bigr)$,
the linear maps $\gamma_u$  are invertible.
The norms of $\gamma_u$ and $\gamma_u^{-1}$
are uniformly dominated on $\nbiga$.

Therefore,
by Lemma \ref{lem;25.6.20.10},
there exists $\delta>0$
and $C_2>0$ such that
\[
       |v|_{h_{\gminig,K(\gminia)}}
       \leq
       C_2
       \Bigl|
       \bigl[u,\rho_{\gminik(\gminia)}(\Ad(\exp(2v))u)\bigr]
       \Bigr|_{h_{\gminig,K(\gminia)}}
\]       
on 
$\bigl\{
(v,u)\in\sqrt{-1}\gminik(\gminia)^{\sigma}\times\nbiga\,\big|\,
F_0(v,u)\leq \delta
\bigr\}$.
By Proposition \ref{prop;25.5.13.31}
and Proposition \ref{prop;25.6.14.20},
\begin{equation}
\label{eq;25.6.20.12}
 \bigl\{
(v,u)\in\sqrt{-1}\gminik(\gminia)^{\sigma}\times\nbiga\,\big|\,
\delta\leq F_0(v,u)\leq r_1
\bigr\}
\end{equation}
is compact.
By Lemma \ref{lem;25.6.20.10},
$\Bigl|
\bigl[u,\rho(\Ad(\exp(2v))u)\bigr]
\Bigr|_{h_{\gminig,K(\gminia)}}>0$
on the set (\ref{eq;25.6.20.12}).
Then, we obtain the claim of the proposition.
\hfill\qed

\begin{cor}
\label{cor;25.5.15.5}
Let $\nbigb$ be a compact subset of $\gminig_1^{\rs}$.
For any $r_1>0$,
there exist $C,C'>0$
such that the following holds
for any $u\in\nbigb$.
\begin{itemize}
 \item If $K_1\in\Herm(G,\sigma)$
       satisfies
       $|u|_{h_{\gminig,K_1}}^2-|u|_{h_{\gminig,K(u)}}^2\leq r_1$,
       then we obtain
       \[
       \bigl|
       \Ad(s(K(u),K_1))-\id
       \bigr|_{h_{\gminig,K(u)}}
       \leq
       C\bigl|
       [u,\rho_{\gminik_1}(u)]
       \bigr|_{h_{\gminig,K_1}},
       \]
\[
       |v(K(u),K_1)|_{h_{\gminig,K(u)}}
       \leq
       C'\bigl|
       [u,\rho_{\gminik_1}(u)]
       \bigr|_{h_{\gminig,K_1}}.
\]             
       \hfill\qed
\end{itemize}
\end{cor}

\subsection{Split real forms}

One of the main examples of
split automorphisms
coming from split real forms.

\subsubsection{Split real forms of semisimple complex Lie algebras}
\label{subsection;25.6.15.30}

Let $\gminig$ be a semisimple complex Lie algebra.
Let $\gminig_{\real}$ be a real form of $\gminig$,
i.e., $\gminig_{\real}$ is a real Lie subalgebra of $\gminig$
such that $\gminig_{\real}\otimes_{\real}\cnum=\gminig$.
\begin{df}
\label{df;25.6.15.40}
The real form $\gminig_{\real}$ is called a split real form
if the following conditions are satisfied.
\begin{itemize}
 \item There exists an abelian subalgebra
       $\gminia_{\real}\subset\gminig_{\real}$
       such that
       $\gminia=\gminia_{\real}\otimes\cnum$
       is a Cartan subalgebra of $\gminig$.
 \item The real part of $\gminia$ equals
       $\gminia_{\real}$.
       (See {\rm\S\ref{subsection;25.6.15.20}}
       for the real part of a Cartan subalgebra.)
       \hfill\qed
\end{itemize}
\end{df}
\begin{prop}[\mbox{\cite[Corollary 6.10, Chapter VI Problem 1]{Knapp-Book}}]
\label{prop;25.6.15.41}
Any semisimple complex Lie algebra $\gminig$
has a split real form $\gminig_{\real}$.
It is unique up to automorphisms of $\gminig$.
\hfill\qed
\end{prop}

\subsubsection{Split real forms and Cartan involutions}

Let $B_{\gminig_{\real}}$ be the Killing form of $\gminig_{\real}$.
A Lie algebra involution $\sigma_{\real}$ of $\gminig_{\real}$
is called a Cartan involution if
the symmetric bilinear form
$-B_{\gminig_{\real}}(u,\sigma_{\real}(v))$
is positive definite on $\gminig_{\real}$.
The eigen decomposition
$(\gminig_{\real},\sigma_{\real})
=(\gminig_{\real}^{\sigma_{\real}=1},\id)
\oplus
(\gminig_{\real}^{\sigma_{\real}=-1},-\id)$
is called the Cartan decomposition
associated with $\sigma_{\real}$.
Any semisimple Lie algebra has a Cartan involution
\cite[Corollary 6.18]{Knapp-Book}.

\begin{lem}
Suppose that $\gminig_{\real}$ is a split real form of $\gminig$.
Let $\sigma_{\real}$ be a Cartan involution of $\gminig_{\real}$.
Then,
 $\gminig^{\sigma_{\real}=-1}_{\real}$ contains 
an abelian subspace $\gminia_{\real}$
as in Definition {\rm\ref{df;25.6.15.40}}.
Namely, any split real form $\gminig_{\real}$ of $\gminig$
is a normal real form of $\gminig$
in the sense of
{\rm\cite[Page 426]{Helgason1978}}.
\end{lem}
\pf
According to \cite[Chapter IX, Theorem 5.10]{Helgason1978},
$\gminig$ has a normal real form
which is unique up to automorphisms.
Because a normal real form is a split real form in the sense of Definition
\ref{df;25.6.15.40},
any split real form is normal
by Proposition \ref{prop;25.6.15.41}.
\hfill\qed

\vspace{.1in}

Suppose that $\gminig_{\real}$ is a split real form.
Let $\sigma_{\real}$ be a Cartan involution of $\gminig_{\real}$.
The induced $\cnum$-linear involution on $\gminig$
is denoted by $\sigma$.
\begin{prop}
\label{prop;25.9.27.3}
$(\sigma,-1)$ is split.
\end{prop}
\pf
We obtain the decomposition
$(\gminig,\sigma)=
(\gminig_0,\id)
\oplus
(\gminig_1,-\id)$.
Let
$\gminia_{\real}\subset\gminig^{\sigma_{\real}=-1}_{\real}
=\gminig_1\cap \gminig_{\real}$
be an abelian subspace
as in Definition \ref{df;25.6.15.40}.
We obtain the Cartan subalgebra
$\gminia=\gminia_{\real}\otimes\cnum
\subset
\gminig_1$.
Because $\gminia^{\rs}\neq\emptyset$,
we obtain $\gminig_1^{\rs}\neq\emptyset$.

For any $u\in\gminig^{\rs}_1$,
we consider
$C_{\gminig}(u)\cap\gminig_j$ $(j=0,1)$.
We can observe that
$\dim C_{\gminig}(u)\cap\gminig_j$
are locally constant on $\gminig_1^{\rs}$.
If $u\in\gminia^{\rs}$,
we obtain $C_{\gminig}(u)\cap \gminig_0
=\gminia\cap\gminig_0=0$.
Because $\gminig_1^{\rs}$ is connected,
we obtain
$\dim C_{\gminig}(u)\cap\gminig_0=0$
for any $u\in\gminig_1^{\rs}$.
\hfill\qed

\subsubsection{Split involutions and split real forms}

Let us study the converse.
Let $\sigma$ be a Lie algebra involution of $\gminig$.
We have the eigen decomposition
\[
(\gminig,\sigma)
=(\gminig_0,\id)
\oplus
(\gminig_{1},-\id).
\]
Suppose that $(\sigma,-1)$ is split.
Let $u\in\gminig_1^{\rs}$.
We obtain the Cartan subspace $\gminia=C_{\gminig}(u)$.
By Proposition \ref{prop;25.5.13.1},
we also obtain the unique maximal compact subgroup
$K\in\Herm(G,\sigma)\cap\Herm(G,\gminia)$.
Because $\sigma(K)=K$,
we have
$\sigma\circ\rho_{\gminik}=\rho_{\gminik}\circ\sigma$.
We obtain the decomposition
$\gminik
=(\gminik\cap\gminig_0)
\oplus
(\gminik\cap\gminig_1)$.
We obtain the following real form
\[
 \gminig_{\sigma,\real}
 =(\gminik\cap\gminig_0)
 \oplus
 \sqrt{-1}
 \bigl(
 \gminik\cap\gminig_1
 \bigr).
\]

\begin{prop}
\label{prop;25.9.27.4}
$\gminig_{\sigma,\real}$
is a split real form of $\gminig$. 
\end{prop}
\pf
Because $(\sigma,-1)$ is split,
$C_{\gminig}(\gminia)=\gminia$ holds,
and hence $\gminia$ is a Cartan subalgebra.
Because $\rho_{\gminik}(\gminia)=\gminia$,
we obtain the decomposition
$\gminia=
 (\gminia\cap\gminik)\oplus
 \sqrt{-1}(\gminia\cap\gminik)$.
We set
$\gminia_{\real}
=\sqrt{-1}(\gminia\cap\gminik)$.
It is the real part of $\gminia$
in the sense of \S\ref{subsection;25.6.15.20}.
Because $\gminia\subset\gminig_1$,
we obtain
$\gminia_{\real}\subset
\sqrt{-1}(\gminik\cap\gminig_1)
\subset
\gminig_{\sigma,\real}$.
Hence, $\gminig_{\sigma,\real}$ is a split real form of $\gminig$.
\hfill\qed

\subsection{Cyclic case}
\label{subsection;25.6.18.20}

Let us explain another main example of
split automorphism
appeared in \cite{Kostant-TDS}.
In this subsection, we assume the following.
\begin{itemize}
 \item $\gminig$ is a simple complex Lie algebra.
 \item $G$ is the associated Lie group of adjoint type,
       i.e.,
       the connected closed algebraic subgroup
       of $\GL(\gminig)$
       whose Lie algebra is
       $\gminig\subset\gl(\gminig)$.
\end{itemize}
Let $\gminit$ be a Cartan subalgebra of $\gminig$.
Let $\gminig_{\phi}\subset\gminig$ denote the root space
corresponding to $\phi\in\Delta(\gminit,\gminig)$.
There exists the root decomposition 
$\gminig=
\gminit\oplus
 \bigoplus_{\phi\in\Delta(\gminit,\gminig)}\gminig_{\phi}$.

\subsubsection{Grading by the lengths}
We fix a Weyl chamber.
Let $\alpha_1,\ldots,\alpha_{\ell}\in\gminit^{\lor}$
denote the corresponding set of the positive simple roots.
Let $\epsilon_1,\ldots,\epsilon_{\ell}\in\gminit$
be the dual base,
i.e.,
$\alpha_i(\epsilon_i)=1$
and $\alpha_i(\epsilon_j)=0$ $(i\neq j)$.
We set $\ttx_0=\sum_{i=1}^{\ell} \epsilon_i$.
We set $\gminig^{(0)}=\gminit$.
For $j\in\seisuu\setminus\{0\}$,
we set
\[
 \gminig^{(j)}=\bigoplus_{\phi(\ttx_0)=j}
 \gminig_{\phi}.
\]
Because $\gminig$ is assumed to be simple,
there exists the highest root $\psi$.
We set $\tth=\psi(\ttx_0)\in\seisuu_{>0}$.
We have $\gminig^{(j)}=0$ if $|j|>\tth$.

\subsubsection{An automorphism of $\gminig$}
\label{subsection;25.5.4.2}
Let $G$ be the adjoint  group of $\gminig$.
Let $T\subset G$ be the Lie group corresponding to $\gminit$
which is a maximal torus in $G$.
We set
\[
 \ttw=\exp\bigl(
 2\pi\sqrt{-1}(\tth+1)^{-1}\ttx_0
 \bigr)\in T.
\]
We have $\ttw^{\tth+1}=1$.

We put $\omega=\exp(2\pi\sqrt{-1}(\tth+1)^{-1})$.
The automorphism $\Ad(\ttw)$ preserves
the decomposition $\bigoplus \gminig^{(j)}$,
and the eigenvalues on
$\gminig^{(j)}$ are $\omega^j$.

\begin{prop}
\label{prop;25.9.27.10}
$(\Ad(\ttw),\omega)$ is split.
\end{prop}
\pf
We consider the eigen decomposition (\ref{eq;25.5.13.11})
for the endomorphism $\Ad(\ttw)$ on $\gminig$.
We have
$\gminig_1
=\bigoplus_{i=1}^{\ell}\gminig_{\alpha_i}
\oplus\gminig_{-\psi}$.
According to the decomposition,
any element $u\in\gminig_1$
is described as
$u=\sum_{i=1}^{\ell} u_i+u_{-\psi}$
According \cite[\S6]{Kostant-TDS},
$u$ is regular semisimple
if and only if
$u$ is cyclic in the sense
$u_i\neq 0$ $(i=1,\ldots,\ell)$
and $u_{-\psi}\neq 0$.
Because $\gminig_0=\gminit$,
we can check
$C_{\gminig}(u)\cap\gminig_0=0$
for any $u\in\gminig_1^{\rs}$
by a direct computation.
\hfill\qed

\section{Harmonic $G$-bundles}

\subsection{Preliminary for linear algebraic groups and Lie algebras}

\subsubsection{Linear algebraic groups}

Let $\hyperk$ be a field of characteristic $0$.
Let $V$ be a finite dimensional $\hyperk$-vector space.
Let $\GL(V)$ denote the group of $\hyperk$-linear automorphisms of $V$.
Let $G$ be a Zariski closed algebraic subgroup of $\GL(V)$.
We recall some Tannakian arguments
by following Simpson in \cite{s5}.

For $(a,b)\in\seisuu_{\geq 0}^2$,
we set $T^{a,b}V=\Hom(V^{\otimes b},V^{\otimes a})$.
For $g\in \GL(V)$ and for any $a\in\seisuu_{\geq 1}$,
we obtain $g^{\otimes a}\in\GL(V^{\otimes \,a})$
satisfying $g^{\otimes a}(v_1\otimes\cdots\otimes v_a)
=g(v_1)\otimes\cdots \otimes g(v_a)$.
We set $g^{\otimes 0}=1$.
we obtain
$T^{a,b}(g)\in \GL(T^{a,b}V)$
satisfying
$T^{a,b}(g)(f)
=g^{\otimes a}\circ f\circ (g^{-1})^{\otimes b}$.

For any $\veca=(a_i)\in\seisuu_{\geq 0}^m$
and $\vecb=(b_i)\in\seisuu_{\geq 0}^m$,
we set
$T^{\veca,\vecb}(V)=\bigoplus_{i=1}^mT^{a_i,b_i}(V)$.
For any $g\in\GL(V)$,
we set
$T^{\veca,\vecb}(g)=\bigoplus_{i=1}^m T^{a_i,b_i}(g)$.
Let $\nbigs_G(\veca,\vecb)$ denote the set of
$G$-subrepresentations of $T^{\veca,\vecb}(V)$.

We recall some fundamental and useful facts
by following Deligne \cite{DMOS} and Simpson \cite{s5}.

\begin{prop}[Chevalley's theorem, \mbox{\cite[\S I. Proposition 3.1]{DMOS}}]
For any closed subgroup $H$ of $G$,
there exists a finite dimensional representation $W$ of $G$
and a one dimensional subspace $L\subset W$ such that
$H$ is characterized as the stabilizer of $L$ in $G$,
 i.e.,
$H=\bigl\{g\in G\,|\,gL\subset L\bigr\}$. 
\hfill\qed
\end{prop}

\begin{prop}[\mbox{\cite[\S I. Proposition 3.1]{DMOS}}]
\label{prop;25.9.21.1}
Any finite dimensional $G$-representation is
contained in $\nbigs_G(\veca,\vecb)$
for some $(\veca,\vecb)$.
\hfill\qed
\end{prop}

\begin{cor}
\label{cor;25.6.5.1}
There exist $(\veca,\vecb)\in (\seisuu_{\geq 0}^m)^2$
and $W\in\nbigs_G(\veca,\vecb)$
such that
$G=\bigl\{
 g\in \GL(V)\,\big|\,T^{\veca,\vecb}(g)W\subset W
 \bigr\}$.
\end{cor}
\pf
There exist a finite dimensional $\GL(V)$-representation $U$
and a one dimensional subspace $L\subset U$
such that $G=\{g\in \GL(V)\,|\,gL\subset L\}$.
There exists
$(\veca,\vecb)\in(\seisuu_{\geq 0}^m)^2$
such that
$U$ is isomorphic to
$W\in \nbigs_{\GL(V)}(\veca,\vecb)$.
Then, the claim follows from the previous propositions.
\hfill\qed

\vspace{.1in}

Let $H\subset G$ be any closed subgroup.
For any $(a,b)\in\seisuu_{\geq 0}^2$,
and for any $W\in\nbigs_G(a,b)$,
let $W^H$ denote the $H$-invariant part of $W$.
Let $H'\subset G$ denote the subgroup of
$g\in G$
such that $T^{a,b}(g)w=w$ holds
for any $(a,b)\in\seisuu_{\geq 0}^2$, $W\in\nbigs_G(a,b)$ and $w\in W^H$.
We always have $H\subset H'$.

\begin{prop}[\mbox{\cite[\S I. Proposition 3.1]{DMOS}}]
If $H$ is reductive,
we obtain $H=H'$.
\hfill\qed
\end{prop}

\begin{cor}
\label{cor;25.9.21.2}
If $H$ is reductive,
there exist $(\veca,\vecb)\in(\seisuu_{\geq 0}^m)^{2}$,
and $w\in (T^{\veca,\vecb}V)^H$
such that $H=\{g\in G\,|\,T^{\veca,\vecb}(g)w=w\}$.
\end{cor}
\pf
Because $H$ is finite dimensional,
there exist finite tuple $(a_i,b_i)$ $(i=1,\ldots,m)$,
$W_i\in \nbigs_G(a_i,b_i)$
and $w_i\in W_i^H$
such that
$H=\bigcap_{i=1}^m\bigl\{
 g\in G\,\big|\,T^{a_i,b_i}(g)w_i=w_i
 \bigr\}$. 
Then, the claim of the lemma follows.
\hfill\qed

\begin{cor}
If $G$ is reductive,
there exist $(\veca,\vecb)\in(\seisuu_{\geq 0}^m)^2$
and $w\in (T^{\veca,\vecb}V)^G$
such that 
\[
 G=\bigl\{
g\in \GL(V)\,\big|\,T^{\veca,\vecb}(g)w=w
 \bigr\}. 
\]
In particular,
we obtain
$G=\bigl\{
 g\in \GL(V)\,\big|\,T^{\veca,\vecb}(g)(w)=w,
 \,\,\,\forall w\in  (T^{\veca,\vecb}V)^G
\bigr\}$
for the tuple $(\veca,\vecb)$.
\hfill\qed
\end{cor}

We also recall a Tannakian duality,
which directly follows from the previous propositions.
Any finite dimensional $G$-representation
is a pair $\vecU=(U,\chi_{\vecU})$ of
a finite dimensional $\hyperk$-vector space $U$
and an algebraic homomorphism $\chi_{\vecU}:G\to\GL(U)$.
A morphism $F:\vecU_1\to\vecU_2$
of finite dimensional $G$-representations
is a $\hyperk$-linear map $F:U_1\to U_2$
such that $\chi_{U_2}(g)\circ F=F\circ\chi_{U_1}(g)$
for any $g\in G$.
Let $\Rep(G)$ denote
the category of finite dimensional $G$-representations.
For any $\vecU_1,\vecU_2\in\Rep(G)$,
the induced $G$-representations
$\vecU_1\oplus\vecU_2$,
$\vecU_1\otimes\vecU_2$
and $\Hom(\vecU_1,\vecU_2)$
are defined in the following natural ways.
\begin{itemize}
 \item For any $\vecU_1,\vecU_2\in\Rep(G)$,
       we have
       $\chi_{\vecU_1\oplus \vecU_2}(g)
       =\chi_{\vecU_1}(g)\oplus \chi_{\vecU_2}(g)$
       and 
       $\chi_{\vecU_1\otimes \vecU_2}(g)=
       \chi_{\vecU_1}(g)\otimes \chi_{\vecU_2}(g)$.
       We also have
       $\chi_{\Hom(\vecU_1,\vecU_2)}(a)=
       \chi_{\vecU_2}(g)\circ a\circ \chi_{\vecU_1}(g)^{-1}$
       for any $a\in\Hom(U_1,U_2)$.
\end{itemize}

\begin{cor}
Suppose that each $\vecU=(U,\chi_{\vecU})\in\Rep(G)$ is equipped with
$g_{\vecU}\in\GL(U)$ satisfying the following conditions.
\begin{itemize}
 \item For any monomorphism $F:\vecU_1\to \vecU_2$ in $\Rep(G)$,
       $F\circ g_{\vecU_1}=g_{\vecU_2}\circ F$ holds.
 \item For any $\vecU_1,\vecU_2\in\Rep(G)$,
       we have
       $g_{\vecU_1\oplus \vecU_2}=g_{\vecU_1}\oplus g_{\vecU_2}$
       and 
       $g_{\vecU_1\otimes \vecU_2}=g_{\vecU_1}\otimes g_{\vecU_2}$.
       We also have
       $g_{\Hom(\vecU_1,\vecU_2)}(a)=
       g_{\vecU_2}\circ a\circ g_{\vecU_1}^{-1}$
       for any $a\in\Hom(U_1,U_2)$.
\end{itemize}
Then, there uniquely exists $g\in G$
such that $g_{\vecU}=\chi_{\vecU}(g)$ $(\vecU\in \Rep(G))$.
\end{cor}
\pf
We regard $V$ as a faithful representation of $G$
by the inclusion $G\subset\GL(V)$.
We have $g_V\in\GL(V)$.
For any $(\veca,\vecb)$ and $U\in\nbigs_G(\veca,\vecb)$,
the inclusion
$\iota_U:U\to T^{\veca,\vecb}(V)$
is a morphism of $G$-representations.
By the assumption,
we have
$\iota_U\circ g_U=
g_{T^{\veca,\vecb}V}\circ\iota_U
=T^{\veca,\vecb}(g_V)\circ\iota_U$.
It implies $T^{\veca,\vecb}(g_{V})(U)\subset U$.
By Corollary \ref{cor;25.6.5.1},
there exists $g\in G$
which induces $g_{V}$ by the inclusion $G\subset \GL(V)$.

Let $\vecV$ denote the faithful $G$-representation
given by $V$ and $G\subset\GL(V)$.
By Proposition \ref{prop;25.9.21.1},
for any $\vecV_1=(V_1,\chi_{\vecV_1})\in\Rep(G)$,
there exist $(\veca,\vecb)\in (\seisuu_{\geq 0}^m)^2$
and a monomorphism of $G$-representations
$\iota:\vecV_1\to T^{\veca,\vecb}(\vecV)$.
Because
$\iota\circ g_{\vecV_1}
=\chi_{T^{\veca,\vecb}(\vecV)}(g)\circ\iota
=\iota\circ\chi_{\vecV_1}(g)$,
we obtain $g_{\vecV_1}=\chi_{\vecV_1}(g)$.
\hfill\qed

\subsubsection{Lie algebras}

We obtain the homomorphism of Lie algebras
$t^{\veca,\vecb}:
\End(V)\to
\End(T^{\veca,\vecb}V)$
as the derivative of the homomorphism of algebraic groups
$\GL(V)\to \GL(T^{\veca,\vecb}(V))$.
Let $\gminig$ denote the Lie algebra of $G$.
We obtain the following from Corollary \ref{cor;25.6.5.1}.
\begin{cor}
\label{cor;24.5.27.20}
There exist
$(\veca,\vecb)\in(\seisuu_{\geq 0}^m)^2$
and $W\in \nbigs_G(\veca,\vecb)$
such that 
 $\gminig=\bigl\{
 f\in \End(V)\,\big|\,
 t^{\veca,\vecb}(f)W\subset W
 \bigr\}$. 
\hfill\qed
\end{cor}
We obtain the following from Corollary \ref{cor;25.9.21.2}.
\begin{cor}
Suppose that $G$ is reductive.
There exist $(\veca,\vecb)\in(\seisuu_{\geq 0}^m)^2$
and $w\in (T^{\veca,\vecb}V)^G$
such that 
\[
 \gminig=\bigl\{
f\in \End(V)\,\big|\,t^{\veca,\vecb}(f)w=0
 \bigr\}. 
\]
In particular,
we obtain
$\gminig=\bigl\{
f\in \End(V)\,\big|\,t^{\veca,\vecb}(f) (T^{\veca,\vecb}V)^G=0
\bigr\}$.
\hfill\qed
\end{cor}

For $\vecU=(U,\chi_{\vecU})$,
we obtain the Lie algebra homomorphism
$\chi_{\vecU}:\gminig\to\End(U)$.
For any $\vecU_1,\vecU_2\in\Rep(G)$ and $f\in\gminig$,
we have
\[
       \chi_{\vecU_1\oplus \vecU_2}(f)=
       \chi_{\vecU_1}(f)\oplus \chi_{\vecU_2}(f),\quad
       \chi_{\vecU_1\otimes \vecU_2}(f)
       =\chi_{\vecU_1}(f)\otimes\id_{U_2}
       + \id_{U_1}\otimes \chi_{\vecU_2}(f).
\]
       We also have
       $\chi_{\Hom(\vecU_1,\vecU_2)}(f)(a)=
       \chi_{\vecU_2}(f)\circ a-a\circ \chi_{\vecU_1}(f)$
       for any $a\in\Hom(U_1,U_2)$.

\begin{cor}
\label{cor;24.5.30.1}
Suppose that each $\vecU\in\Rep(G)$ is equipped with
$f_{\vecU}\in\End(U)$ satisfying the following conditions.
\begin{itemize}
 \item For any monomorphism $F:\vecU_1\to \vecU_2$ in $\Rep(G)$,
       $F\circ f_{\vecU_1}=f_{\vecU_2}\circ F$ holds.
 \item For any $\vecU_1,\vecU_2\in\Rep(G)$,
       we have
       $f_{\vecU_1\oplus \vecU_2}=f_{\vecU_1}\oplus f_{\vecU_2}$
       and 
       $f_{\vecU_1\otimes \vecU_2}=f_{\vecU_1}\otimes\id_{U_2}
       + \id_{U_1}\otimes f_{\vecU_2}$.
       We also have
       $f_{\Hom(\vecU_1,\vecU_2)}(a)
       =f_{\vecU_2}\circ a-a\circ f_{\vecU_1}$
       for any $a\in\Hom(U_1,U_2)$.
\end{itemize}
Then, there uniquely exists $f\in\gminig$
such that
$\chi_{\vecU}(f)=f_{\vecU}$ $(\vecU\in\Rep(G))$.
\hfill\qed
\end{cor}

\subsubsection{Representations of reductive complex algebraic group}

Let us consider the case $\hyperk=\cnum$,
and $G$ is reductive.
Let $\Irr(G)$ denote the set of irreducible representations of $G$.
For any $\vecV\in\Rep(G)$,
we obtain the canonical decomposition of $G$-representations
\begin{equation}
\label{eq;25.5.5.20}
 T^{\veca,\vecb}\vecV
 =\bigoplus_{\vecW\in \Irr(G)}
 (T^{\veca,\vecb}\vecV)_{\vecW},
\end{equation}
where 
$(T^{\veca,\vecb}\vecV)_{\vecW}
\simeq
\vecW
\otimes \cnum^{m(\vecW,\veca,\vecb)}$.
Let $(\cnum,\II)$ denote the trivial $G$-representation,
i.e., $\II:G\to\{1\}\subset\cnum^{\ast}$.
The $G$-invariant part
$(T^{\veca,\vecb}\vecV)^G$ equals
$(T^{\veca,\vecb}\vecV)_{(\cnum,\II)}$.

\begin{cor}
\mbox{{}}\label{cor;25.6.16.2}
 \begin{itemize}
 \item 
$g\in \GL(V)$ is contained in $\chi_V(G)$
       if and only if
       for any $(\veca,\vecb)\in(\seisuu_{\geq 0}^{m})^2$,
       $T^{\veca,\vecb}(g)$ preserves the decomposition
       {\rm(\ref{eq;25.5.5.20})},
       and
       $T^{\veca,\vecb}(g)w=w$ holds
       for any $w\in (T^{\veca,\vecb}V)^{G}$.
 \item
      $f\in \End(V)$ is contained in $\chi_V(\gminig)$
      if and only if
      for any $(\veca,\vecb)\in(\seisuu_{\geq 0}^{m})^2$,
      $t^{\veca,\vecb}(f)$
      preserves the decomposition
      {\rm(\ref{eq;25.5.5.20})},
      and      
      $t^{\veca,\vecb}(f)w=0$ holds
      for any $w\in (T^{\veca,\vecb}V)^G$.
      \hfill\qed
 \end{itemize}
\end{cor}

Let $K\subset G$ be a maximal compact subgroup.
The Cartan involution is denoted by $\rho$.
For any $\vecU=(U,\chi_{\vecU})\in\Rep(G)$,
we obtain the homomorphism
$\chi_{\vecU|K}:K\to \GL(U)$
and the $K$-representation
$\vecU_{|K}:=(U,\chi_{\vecU|K})$.
Let $\vecU_i\in\Rep(G)$.
If $\vecU_{i|K}$ $(i=1,2)$ are isomorphic as $K$-representations,
$\vecU_i$ $(i=1,2)$ are isomorphic as $G$-representations
because $G$ is the complexification of $K$.

Let $\vecV\in\Rep(G)$ be a faithful representation.
Let $h_V$ be a $K$-invariant Hermitian metric of $V$.
Let $T^{\veca,\vecb}(h_V)$
denote the induced Hermitian metric
on $T^{\veca,\vecb}(V)$
for any $(\veca,\vecb)$.

\begin{lem}
The decomposition {\rm(\ref{eq;25.5.5.20})}
is orthogonal with respect to $T^{\veca,\vecb}h_V$.
\end{lem}
\pf
For any $\vecW\in\Irr(G)$,
the induced $K$-representations $W_{|K}$
is irreducible.
Moreover, if $\vecW\neq \vecW'$ in $\Irr(G)$,
the induced $K$-representations
$(W,\chi_{\vecW|K})$ and $(W',\chi_{\vecW'|K})$
are not isomorphic.
Then, we obtain the desired orthogonality.
\hfill\qed

\begin{lem}
Let $g=k\exp(v)\in G$ be the Cartan decomposition,
i.e., $k\in K$ and $v\in \sqrt{-1}\gminik$. 
If $g\in G$ satisfies
$\chi_V(g)^{\dagger}_{h_V}=\chi_V(g)$,
then $k^2=1$ and $\Ad(k)(v)=v$.
If moreover $\chi_V(g)$ is positive definite,
then $k=1$,
i.e., $g\in \exp(\sqrt{-1}\gminik)$. 
\end{lem}
\pf
Let $\gminiu(V,h_V)\subset\End(V)$ denote the real Lie subalgebra
of endomorphisms of $V$ which are anti-Hermitian with respect to $h_V$.
We obtain a homomorphism of the real Lie algebras
$\chi_V:\gminik\to \gminiu(V,h_V)$
whose complexification is
$\chi_V:\gminig\to\End(V)$.
We have
$\chi_V(\rho(u))
 =-\chi_V(u)^{\dagger}_{h_V}$
for any $u\in \gminig$.
We also obtain
$\chi_V(\rho(a))
=(\chi_V(a)^{\dagger}_{h_V})^{-1}$
for any $a\in G$.

Suppose $g=k\cdot \exp(v)$
satisfies $\chi_V(g)^{\dagger}_{h_V}=\chi_V(g)$,
where $k\in K$ and $v\in\sqrt{-1}\gminik$.
We obtain
\[
\chi_V(k\exp(v))=
\chi_V(k\exp(-v))^{-1}
=\chi_V(\exp(v)k^{-1})
=\chi_V(k^{-1}\exp(\Ad(k)v)).
\]
We obtain $k^2=1$ and $\Ad(k)v=v$
by the uniqueness of Cartan decomposition.
If moreover $\chi_V(g)$ is positive definite,
we obtain $\chi_V(k)=1$
by the uniqueness of the Cartan decomposition.
\hfill\qed

\begin{lem}
\label{lem;25.5.5.40}
For any $g\in G$,
we set $h_V^g(v_1,v_2)=h_V(\chi_V(g)v_1,\chi_V(g)v_2)$.
\begin{itemize}
 \item $h_V^g$ is a $g^{-1}Kg$-invariant Hermitian metric of $V$.
 \item $h_V^g$ depends only on
       $\rho(g^{-1})g\in\exp(\sqrt{-1}\gminik)$.
 \item The decomposition {\rm(\ref{eq;25.5.5.20})}
       is orthogonal with respect to $T^{\veca,\vecb}h_V^g$.
 \item
       Moreover,
       $T^{\veca,\vecb}h_V^g=T^{\veca,\vecb}h_V$
      on the $G$-invariant part
      $T^{\veca,\vecb}(V)^G$.
\end{itemize} 
\end{lem}
\pf
The first two claims follow from the definition of $h_V^g$.
The other claims follow from Corollary \ref{cor;25.6.16.2}.
\hfill\qed

\begin{lem}
\label{lem;25.5.5.41}
Let $h$ be any Hermitian metric of $V$
satisfying the following conditions.
\begin{itemize}
 \item The decomposition {\rm(\ref{eq;25.5.5.20})}
       is orthogonal with respect to $T^{\veca,\vecb}h$.
 \item $T^{\veca,\vecb}(h_V)=T^{\veca,\vecb}(h)$
       on the $G$-invariant part
       $T^{\veca,\vecb}(V)^G$.
\end{itemize}
Then, there exists a unique element
$g\in\exp(\sqrt{-1}\gminik)$ 
such that $h=h_V^g$.
\end{lem}
\pf
Let $s(h_V,h)$ be the automorphism of $V$
determined by
$h(v_1,v_2)=h_V(s(h_V,h)v_1,v_2)$.
We obtain that
$T^{\veca,\vecb}(s(h_V,h))$ preserves
the decomposition (\ref{eq;25.5.5.20}),
and the restriction of
$T^{\veca,\vecb}(s(h_V,h))$ to
$T^{\veca,\vecb}(V)^G$ equals the identity.
Hence, there exists $g_1\in G$
such that $s(h_V,h)=\chi_V(g_1)$.
Because $s(h_V,h)$ is self-adjoint
and positive definite with respect to $h_V$,
we obtain that $g_1=\exp(b)$
for some $b\in\sqrt{-1}\gminik$.
We obtain
$h=h_V^{\exp(b/2)}$.
\hfill\qed

\subsection{Preliminary for principal $G$-bundles}

\subsubsection{Principal $G$-bundles and the associated bundles}

Let $G$ be a Lie group.
Let $X$ be any differentiable manifold.
Let $\gbigp_G$ be a principal $G$-bundle on $X$
which is equipped with a free right $G$-action $\kappa$.
For any space $Y$ equipped with a $G$-action $\kappa_Y$,
let $\gbigp_G(Y)$ denote
the quotient space of
$\gbigp_G\times Y$
by the left $G$-action given by
$g(x,y)=(xg^{-1},gy)$,
where $xg^{-1}=\kappa(g^{-1},x)$
and $gy=\kappa_Y(g,y)$.
It is also denoted by $\gbigp_G(Y,\kappa_Y)$
when we emphasize the action $\kappa_Y$.
There exists the natural
$G$-equivariant isomorphism
$\gbigp_G\times_X\gbigp_G(Y)
\simeq
\gbigp_G\times Y$.

Let $\Gamma(\gbigp_G(Y))$ denote the space of
sections of the bundle $\gbigp_G(Y)$ on $X$.
Note that a section of the bundle $\gbigp_G(Y)$ on $X$
corresponds to
a map $f:\gbigp_G\to Y$ satisfying $f(xg^{-1})=gf(x)$.

\subsubsection{Extensions and reductions of principal bundles}

Let $f:G\to G_0$ be any morphism of Lie groups.
The space $G_0$ is equipped with the left $G$-action induced by $f$,
and the natural right $G_0$-action.
We obtain the induced principal $G_0$-bundle $\gbigp_G(G_0)$.

For a Lie subgroup $G_1\subset G$,
a $G_1$-reduction of $G$
is a principal $G_1$-bundle $\gbigp_{G_1}$
with an isomorphism
$\gbigp_{G_1}(G)\simeq \gbigp_G$.
We naturally regard $\gbigp_{G_1}$
as a subset of $\gbigp_G$. 
A $G_1$-reduction $\gbigp_{G_1}$ of $\gbigp_G$
naturally corresponds to a section of
$\gbigp_{G}(G/G_1)=(\gbigp_G)/G_1$ on $X$.
Indeed, $\gbigp_{G_1}$ equals
the fiber product of
the projection $\gbigp_G\to \gbigp_{G}(G/G_1)$
and the section $X\to\gbigp_{G_1}(G/G_1)$.

\begin{rem}
A section $h$ of $\gbigp_G(G/G_1)$
is also called a $G_1$-reduction of $\gbigp_G$.
The corresponding principal $G_1$-bundle
is denoted by $\gbigp^h_{G_1}$. 
\hfill\qed 
\end{rem}

\subsubsection{Gauge transformations}

We have the left $G$-action $\Ad$ on $G$
given by $\Ad(g)(x)=gxg^{-1}$.
We obtain the bundle of Lie groups $\gbigp_G(G,\Ad)$.
We set $\nbigg(\gbigp_G):=\Gamma(\gbigp_G(G,\Ad))$,
which is naturally a group,
called the group of gauge transformations of $\gbigp_G$.
For any space $Y$ equipped with a left $G$-action,
we obtain the natural left action
of $\nbigg(\gbigp_G)$ on $\Gamma(\gbigp_G(Y))$.
In particular,
any section of
$\gbigp_G(G,\Ad)$
induces an automorphism of $\gbigp_G(Y)$.

\subsubsection{Connections of principal $G$-bundles}

Let us recall the notion of connections on
a principal $G$-bundle $\gbigp_G$ on $X$.
We obtain the bundle of Lie algebras
$\gbigp_G(\gminig)$
by the adjoint representation of $G$ on $\gminig$.
Let $A^j(\gbigp_G(\gminig))$
denote the space of $\gbigp_G(\gminig)$-valued $j$-forms.

Let $\kappa:\gbigp_G\times G\to \gbigp_G$
denote the right $G$-action on $\gbigp_G$.
Let $\pi:\gbigp_G\to X$ denote the projection.
There is the natural $G$-equivariant isomorphism
$\pi^{\ast}\gbigp_G(\gminig)\simeq\gbigp_G\times\gminig$.
For any $\tau\in A^j(\gbigp_G(\gminig))$,
$\pi^{\ast}(\tau)$ is regarded
as a $\gminig$-valued $j$-form on $\gbigp_G$
satisfying
$\kappa_g^{\ast}\bigl(\pi^{\ast}(\tau)\bigr)
=\ad(g^{-1})\pi^{\ast}(\tau)$.

Let $T^{\ttv}\gbigp_G$ be the subbundle of the tangent bundle $T\gbigp_G$
in the fiber direction, obtained as $\Ker(d\pi)$
for the induced bundle map $d\pi:T\gbigp_G\to \pi^{-1}(TX)$.
For any $x\in \gbigp_G$,
we obtain the induced map
$\kappa_x:G\to \gbigp_G$
by $\kappa_x(g)=\kappa(x,g)$.
We obtain the isomorphism
$T_1\kappa_x:\gminig\simeq T^{\ttv}_x\gbigp_G$,
where $1\in G$ denote the unit element.

One of the equivalent definitions of a connection of $\gbigp_G$
is a $\gminig$-valued $1$-form $\omega$
on $\gbigp_G$ satisfying the following conditions.
\begin{itemize}
 \item $\omega_x(T_1\kappa_x(A))=A$
       for any $x\in\gbigp_G$ and any $A\in\gminig$.
 \item $\kappa_g^{\ast}(\omega)
       =\ad(g^{-1})\circ\omega$
       for any $g\in G$,
       where $\kappa_g:\gbigp_G\to\gbigp_G$
       is defined by
       $\kappa_g(x)=\kappa(x,g)$.
\end{itemize}
There exists a unique element $R(\omega)\in A^2(\gbigp_G(\gminig))$
such that
$\pi^{\ast}R(\omega)=d\omega+\frac{1}{2}[\omega,\omega]$,
which is called the curvature of the connection.

Let $\nbiga(\gbigp_G)$ denote the set of connections of $\gbigp_G$.
For any $\omega_1,\omega_2\in\nbiga(\gbigp_G)$,
$\omega_1-\omega_2$ is a $\gminig$-valued $1$-form on $\gbigp_G$
satisfying
$\kappa_g^{\ast}(\omega_1-\omega_2)
=\ad(g^{-1})\circ(\omega_1-\omega_2)$
and
$(\omega_1-\omega_2)_{|T^{\ttv}\gbigp}=0$.
Hence, there uniquely exists
$f\in A^1(\gbigp_G(\gminig))$
such that
$\omega_1-\omega_2=\pi^{\ast}(f)$.
In this way,
$\nbiga(\gbigp_G)$ is naturally
an affine space over $A^1(\gbigp_G(\gminig))$.

\subsection{Complex linear algebraic case}
\label{subsection;26.3.9.2}

Let $G$ be a complex linear algebraic group,
i.e.,
a closed algebraic subgroup of $\GL_n(\cnum)$ for some $n$.
The contents of \S\ref{subsection;26.3.9.2}--\ref{subsection;26.3.9.3}
are essentially contained in \cite{mochi4}.

\subsubsection{The associated vector bundles}

For any $\vecV\in\Rep(G)$,
we obtain a $C^{\infty}$ complex vector bundle
$\gbigp_G(\vecV)$ on $X$.
Let $A^k(\gbigp_G(\vecV))$ denote the space of
$\gbigp_G(\vecV)$-valued $k$-forms.
For any $\vecV_1,\vecV_2\in\Rep(G)$,
we have
$\gbigp_G(\vecV_1\oplus \vecV_2)
 = \gbigp_G(\vecV_1)\oplus\gbigp_G(\vecV_2)$,
$\gbigp_G(\vecV_1\otimes \vecV_2)
 = \gbigp_G(\vecV_1)\otimes\gbigp_G(\vecV_2)$
and
$\gbigp_G(\Hom(\vecV_1,\vecV_2))
= \Hom(\gbigp_G(\vecV_1),\gbigp_G(\vecV_2))$.
For the trivial representation $(\cnum,\II)$,
$\gbigp_G(\cnum,\II)$ equals the product line bundle
$X\times\cnum$.

\subsubsection{$k$-forms of the adjoint bundles in the linear algebraic case}

For any $\vecV\in\Rep(G)$,
we obtain the bundle map
$\chi_{\vecV}:\gbigp_G(\gminig)\to\End(\gbigp_G(\vecV))$.
For any $\tau\in A^k(\gbigp_G(\gminig))$,
we obtain
$\chi_{\vecV}(\tau)\in A^k(\End(\gbigp_G(\vecV)))$.
If $\vecV$ is a faithful representation,
we may regard $A^k(\gbigp_G(\gminig))$
as a subspace of
$A^k(\End(\gbigp_G(\vecV)))$.

\begin{lem}
\label{lem;25.6.5.10}
Suppose that $\tau_{\vecV}\in A^k(\End(\gbigp_G(\vecV)))$
is attached to each $\vecV\in\Rep(G)$
such that
the following conditions are satisfied.
\begin{itemize}
 \item For any monomorphism $F:\vecV_1\to \vecV_2$ in $\Rep(G)$,
       we obtain
       $F\circ\tau_{\vecV_1}=\tau_{\vecV_2}\circ F$.
 \item For any $V_1,V_2\in\Rep(G)$,
       we obtain
       $\tau_{\vecV_1\oplus \vecV_2}=
       \tau_{\vecV_1}\oplus\tau_{\vecV_2}$
       and
       $\tau_{\vecV_1\otimes \vecV_2}=
       \tau_{\vecV_1}\otimes \id_{V_2}+\id_{V_1}\otimes\tau_{\vecV_2}$.
       Moreover,
       we have
       $\tau_{\Hom(\vecV_1,\vecV_2)}(a)
       =\tau_{\vecV_2}\circ a-a\circ\tau_{\vecV_1}$.
\end{itemize}
Then, there uniquely exists
$\tau\in A^k(\gbigp_G(\gminig))$
such that $\tau_{\vecV}=\chi_{\vecV}(\tau)$.
\end{lem}
\pf
It follows from Corollary \ref{cor;24.5.30.1}.
\hfill\qed

\vspace{.1in}

Let $\vecV\in\Rep(G)$.
We have
\[
 T^{\veca,\vecb}(\gbigp_G(\vecV)):=
 \bigoplus_{i=1}^m
 \Hom\bigl(\gbigp_G(\vecV)^{\otimes b_i},\gbigp_G(\vecV)^{\otimes a_i}\bigr)
\simeq
 \gbigp_G(T^{\veca,\vecb}(\vecV)).
\]
Let $\nbigs_G(\veca,\vecb)$ denote
the set of $G$-subrepresentations of $T^{\veca,\vecb}(\vecV)$.
For any $U\in \nbigs_G(\veca,\vecb)$,
we have the corresponding $\vecU\in\Rep(G)$,
and we obtain the subbundle
$\gbigp_G(\vecU)\subset
T^{\veca,\vecb}(\gbigp_G(\vecV))$.

\begin{lem}
\label{lem;24.5.16.11}
If $\vecV$ is a faithful representation,
there exist $(\veca,\vecb)\in(\seisuu_{\geq 0}^m)^2$
and $U\in\nbigs_G(\veca,\vecb)$
such that
the inclusion
$\gbigp_G(\gminig)\subset \End(\gbigp_G(\vecV))$
induces 
\begin{equation}
\label{eq;24.5.16.12}
 A^k(\gbigp_G(\gminig))\simeq
 \Bigl\{
 f\in A^k\bigl(\End(\gbigp_G(\vecV))\bigr)
 \,\Big|\,
 t^{\veca,\vecb}(f)\bigl(\gbigp_G(\vecU)\bigr)
 \subset
 \gbigp_G(\vecU)\otimes\bigwedge^kT^{\ast}X
 \Bigr\}.
 \end{equation}
 \end{lem}
\pf
It follows from Corollary \ref{cor;24.5.27.20}. 
\hfill\qed

\subsubsection{Connections of principal $G$-bundles and the associated
vector bundles}

A connection $\omega$ of $\gbigp_G$
induces a connection $\nabla^{\omega}_{\gbigp_G(\vecV)}$ of
the vector bundle $\gbigp_G(\vecV)$
for any $\vecV=(V,\chi_V)\in\Rep(G)$.
Indeed, a section of $\gbigp_G(\vecV)$
is regarded as a map $f:\gbigp_G\to V$
satisfying $\kappa_g^{\ast}f=\chi_{\vecV}(g^{-1}) f$.
We obtain the $V$-valued $1$-form $df+\omega f$ on $\gbigp_G$
satisfying
$\kappa_g^{\ast}(df+\omega f)=\chi_{\vecV}(g^{-1})(df+\omega f)$.
Because $(df+\omega f)_{|T^{\ttv}\gbigp_G}=0$,
there exists $\gbigp_G(\vecV)$-valued $1$-form
$\nabla^{\omega}_{\gbigp_G(\vecV)}(f)$
such that
$df+\omega f=\pi^{\ast}(\nabla^{\omega}_{\gbigp_G(\vecV)}(f))$.
The connections $\nabla^{\omega}_{\gbigp_G(\vecV)}$
satisfy the following conditions.
\begin{itemize}
 \item For any $F:\vecV_1\to \vecV_2$ in $\Rep(G)$,
       we obtain
       $F\circ\nabla^{\omega}_{\gbigp_G(\vecV_1)}
       =\nabla^{\omega}_{\gbigp_G(\vecV_2)}\circ F$.
 \item For any $\vecV_1,\vecV_2\in\Rep(G)$,
       the connections
       $\nabla^{\omega}_{\gbigp_G(\vecV_1\oplus \vecV_2)}$,
       $\nabla^{\omega}_{\gbigp_G(\vecV_1\otimes \vecV_2)}$
       and
       $\nabla^{\omega}_{\gbigp_G(\Hom(\vecV_1,\vecV_2))}$
       are the connections
       on the bundles
       $\gbigp_G(\vecV_1\oplus \vecV_2)$,
       $\gbigp_G(\vecV_1\otimes \vecV_2)$
       and
       $\gbigp_G(\Hom(\vecV_1,\vecV_2))$
       naturally induced by
       $\nabla^{\omega}_{\gbigp(\vecV_1)}$
       and $\nabla^{\omega}_{\gbigp(\vecV_2)}$.
 \item For the trivial representation $(\cnum,\II)$,
       the connection
       $\nabla^{\omega}_{\gbigp_G(\cnum,\II)}$
       is the exterior derivative.  
\end{itemize}

\begin{prop}
\label{prop;24.5.30.2}
Suppose that the vector bundles $\gbigp_G(\vecV)$
$(\vecV\in\Rep(G))$
are equipped with connections $\nabla_{\gbigp_G(\vecV)}$
satisfying the following conditions.
\begin{itemize}
 \item For any monomorphism
       $F:\vecV_1\to \vecV_2$ in $\Rep(G)$,
       we obtain
       $F\circ\nabla_{\gbigp_G(\vecV_1)}
       =\nabla_{\gbigp_G(\vecV_2)}\circ F$.
 \item For any $\vecV_1,\vecV_2\in\Rep(G)$,
       the connections
       $\nabla_{\gbigp_G(\vecV_1\oplus \vecV_2)}$,
       $\nabla_{\gbigp_G(\vecV_1\otimes \vecV_2)}$
       and
       $\nabla_{\gbigp_G(\Hom(\vecV_1,\vecV_2))}$
       are the connections
       on the bundles
       $\gbigp_G(\vecV_1\oplus \vecV_2)$,
      $\gbigp_G(\vecV_1\otimes \vecV_2)$
      and
       $\gbigp_G(\Hom(\vecV_1,\vecV_2))$
       naturally induced by
       $\nabla_{\gbigp(\vecV_1)}$
       and $\nabla_{\gbigp(\vecV_2)}$.
\end{itemize}
Then, there exists a unique connection $\omega$ of $\gbigp_G$
such that 
$\nabla^{\omega}_{\gbigp_G(\vecV)}=\nabla_{\gbigp_G(\vecV)}$
 $(\vecV\in\Rep(G))$.
\end{prop}
\pf
Let $\omega_0$ be any connection of $\gbigp_G$.
By Lemma \ref{lem;25.6.5.10},
there uniquely exists
$\tau\in A^1(\gbigp_G(\gminig))$
such that
$\chi_{\vecV}(\tau)=
\nabla_{\gbigp_G(\vecV)}-\nabla^{\omega_0}_{\gbigp_G(\vecV)}$
for any $\vecV\in\Rep(G)$.
By setting $\omega=\omega_0+\pi^{\ast}(\tau)$,
we obtain
$\nabla^{\omega}_{\gbigp_G(\vecV)}
=\nabla_{\gbigp_G(\vecV)}$ $(\vecV\in\Rep(G))$.
The uniqueness is clear.
\hfill\qed

\begin{notation}
\label{notation;25.9.23.1}
Under the assumption that $G$ is complex linear algebraic,
we may regard a connection of $\gbigp_G$
as a tuple of connections
$\nabla=(\nabla_{\gbigp_G(\vecV)}\,|\,\vecV\in\Rep(G))$
satisfying the conditions in Proposition {\rm\ref{prop;24.5.30.2}}.
\hfill\qed
\end{notation}

We also have another characterization of
a connection of $\gbigp_G$.
We set $E:=\gbigp_G(\vecV)$
for a faithful representation $\vecV\in\Rep(G)$.
A connection $\nabla_E$ on $E$
induces a connection
$t^{\veca,\vecb}(\nabla_E)$
on $T^{\veca,\vecb}(E)$.
Let $\nbiga_G(E)$ be the set of
connections $\nabla_E$ of $E$
such that
for any $(\veca,\vecb)\in(\seisuu_{\geq 0}^m)^2$
and $U\in\nbigs_G(\veca,\vecb)$
the induced connection $t^{\veca,\vecb}(\nabla_E)$
preserves $\gbigp_G(\vecU)$.
It is easy to see that
if $\nabla_E$ is induced by a connection of $\gbigp_G$,
it is contained in $\nbiga_G(E)$.
Thus, we obtain a map
$\nbiga(\gbigp_G)\to \nbiga_G(E)$.

\begin{lem}
The map $\nbiga(\gbigp_G)\to \nbiga_G(E)$
is a bijection. 
\end{lem}
\pf
By choosing a connection of $\gbigp_G$,
we obtain an isomorphism of affine spaces
$\nbiga(\gbigp_G)\simeq A^1(\gbigp_G(\gminig))$.
We also obtain an isomorphism of
$\nbiga_G(E)$
and the right hand side of (\ref{eq;24.5.16.12})
with $k=1$.
Then, the claim follows from Lemma \ref{lem;24.5.16.11}.
\hfill\qed

\subsubsection{$K$-reductions in the complex reductive case}
\label{subsection;25.6.9.23}

Suppose that $G$ is a reductive complex algebraic group
with a Cartan involution $\rho$.
We set $K=G^{\rho}$,
which is a maximal compact subgroup.
Let $\gminik$ denote the Lie algebra of $K$.

Let $h$ be a $K$-reduction of $\gbigp_G$.
Because $\gbigp_G(\gminig)=\gbigp^h_K(\gminig)$,
we obtain the involution $\rho_h$ on $\gbigp_G(\gminig)$.
We obtain the decomposition
\begin{equation}
\label{eq;24.5.30.3}
 \gbigp_G(\gminig)
 =\gbigp^h_K(\gminik)
 \oplus
 \gbigp^h_K(\sqrt{-1}\gminik),
\end{equation}
and $\gbigp^h_K(\gminik)$ is the $\rho_h$-invariant part,
and $\gbigp^h_K(\sqrt{-1}\gminik)$ is the $\rho_h$-anti-invariant part.
We obtain
the anti-$\cnum$-linear endomorphism $\rho_h$
on $A^k(\gbigp_G(\gminig))$
induced by $\rho_h$
and the complex conjugation for $k$-forms.

\begin{lem}
\label{lem;24.5.30.4}
For any connection $\omega$ of $\gbigp_G$,
there exists a unique connection $\omega_h$ of $\gbigp^h_K$
and a section $\Phi_h$ of
$A^1(\gbigp^h_K(\sqrt{-1}\gminik))$ 
such that $\omega=\omega_h+\Phi_h$.
\end{lem}
\pf
Let $\omega'_K$ be any connection of $\gbigp^h_K$.
It induces a connection of $\gbigp_G=\gbigp^h_K(G)$,
which is also denoted by $\omega_K'$.
We obtain
$\Psi\in A^1(\gbigp_G(\gminig))$
such that
$\pi^{\ast}(\Psi)=\omega-\omega_K'$.
We have the decomposition
$\Psi=\Psi_{\gminik}+\Psi_{\sqrt{-1}\gminik}$
according to the decomposition (\ref{eq;24.5.30.3}).
We obtain the desired connection
$\omega_h=\omega_K'+\pi^{\ast}\Psi_{\gminik}$
of $\gbigp^h_K$,
and the desired $1$-form
$\Phi_h=\Psi_{\sqrt{-1}\gminik}$.
The uniqueness is clear.
\hfill\qed

\subsubsection{$K$-reductions and the induced metrics
of the associated vector bundles}

Let $\gbigp_G$ be a principal $G$-bundle.
Let $h$ be a $K$-reduction of $\gbigp_G$.
Let $\vecV=(V,\chi_V)\in\Rep(G)$.
Let $h_{\vecV}$ be a Hermitian metric of $V$
which is $K$-invariant under the representation
$\chi_{\vecV|K}$.
We have the compatibility of
the Cartan involutions,
i.e.,
$\chi_{\vecV}\circ\rho(v)
=-\chi_{\vecV}(v)^{\dagger}_{h_{\vecV}}$.
We obtain the Hermitian metric
$h_{\vecVtilde}$ of the vector bundle
$\vecVtilde=\gbigp_G(\vecV)=\gbigp^h_K(\vecV)$
induced by $h$ and $h_{\vecV}$.
For any section $f$ of $\End(\vecVtilde)$,
let $f^{\dagger}_{h_{\vecVtilde}}$
denote the section of $\End(\vecVtilde)$
obtained as the adjoint of $f$ with respect to $h_{\vecVtilde}$.
For the natural morphism
$\chi_{\vecV}:\gbigp_G(\gminig)\to \End(\vecVtilde)$,
we obtain
$\chi_{\vecV}(\rho_h(s))
=-(\chi_{\vecV}(s))^{\dagger}_{h_{\vecVtilde}}$
for any section $s$ of $\gbigp_G(\gminig)$.

Let $\omega$ be a connection of $\gbigp_G$.
It induces a connection $\nabla^{\omega}_{\vecVtilde}$ of $\vecVtilde$.
There exists the unique decomposition
$\nabla^{\omega}_{\vecVtilde}
=\nabla_{h_{\vecVtilde}}+\Phi_{h_{\vecVtilde}}$,
where
$\nabla_{h_{\vecVtilde}}$ is a unitary connection of $\vecVtilde$
with respect to $h_{\vecVtilde}$,
and $\Phi_{h_{\vecVtilde}}\in A^1(\End(\vecVtilde))$
satisfying
$(\Phi_{h_{\vecVtilde}})^{\dagger}_{h_{\vecVtilde}}
=\Phi_{h_{\vecVtilde}}$.
The following lemma is clear by the uniqueness.
\begin{lem}
$\nabla_{h_{\vecVtilde}}$
and $\Phi_{h_{\vecVtilde}}$
are induced by $\omega_h$ and $\Phi_h$
in Lemma {\rm\ref{lem;24.5.30.4}}.
\hfill\qed
\end{lem}

The following lemma is obvious.
\begin{lem}
\label{lem;24.5.30.5}
The connection $\omega$
is induced by a connection of $\gbigp^h_K$
if and only if
there exists a faithful representation $\vecV\in\Rep(G)$
such that 
$\nabla^{\omega}_{\vecVtilde}$ preserves $h_{\vecVtilde}$.
\hfill\qed
\end{lem}

\subsubsection{A characterization of $K$-reductions}

Let $\gbigp_G$ be a principal $G$-bundle on
a $C^{\infty}$-manifold $X$.
Let $\vecV$ be any faithful $G$-representation.
We fix a $K$-invariant Hermitian metric $h_V$ of $V$.
We obtain the vector bundle $E=\gbigp_G(\vecV)$.
Corresponding to (\ref{eq;25.5.5.21}),
we obtain the decomposition
\begin{equation}
 \label{eq;25.5.5.21}
  T^{\veca,\vecb}(E)
  =\bigoplus_{\vecW\in\Irr(G)}
T^{\veca,\vecb}(E)_{\vecW}.
\end{equation}
Note that $T^{\veca,\vecb}(E)_{(\cnum,\II)}$
is the product bundle
$X\times T^{\veca,\vecb}(V)^G$.
It is equipped with the constant metric
obtained as the restriction of
$T^{\veca,\vecb}(h_V)$
to $T^{\veca,\vecb}(V)^G$
denoted by
$T^{\veca,\vecb}(h_V)_{(\cnum,\II)}$.

\begin{prop}
\label{prop;25.5.5.22}
A Hermitian metric $h_E$ of $E$ is
induced by $h_V$ and
a $K$-reduction of $\gbigp_G$
if and only if the following conditions are satisfied for any
$(\veca,\vecb)$.
\begin{description}
 \item[(i)] The decomposition
       {\rm(\ref{eq;25.5.5.21})}
       is orthogonal with respect to
       $T^{\veca,\vecb}h_E$.
 \item[(ii)] The restriction of $T^{\veca,\vecb}h_E$
       to $T^{\veca,\vecb}(E)_{(\cnum,\II)}$
       equals
       $T^{\veca,\vecb}(h)_{(\cnum,\II)}$.
\end{description}
\end{prop}
\pf
If $h_E$ is induced by $h_V$
and a $K$-reduction of $\gbigp_G$,
then (i) and (iii) are satisfied by Lemma \ref{lem;25.5.5.40}.
The converse follows from Lemma \ref{lem;25.5.5.41}.
\hfill\qed

\subsubsection{Description of the induced metrics}

Let $\gbigp_G$ be a principal $G$-bundle on
a $C^{\infty}$-manifold $X$.
Let $h$ be a $K$-reduction of $\gbigp_G$.
Let $\vecV$ be any $G$-representation.
We fix a $K$-invariant Hermitian metric $h_V$ of $V$.
We obtain the Hermitian metric $h_E$ 
of $E=\gbigp_G(\vecV)=\gbigp^h_K(V)$
induced by $h$ and $h_V$.

Let $a$ be a section of $\gbigp_G$.
We obtain the map
$X\times V\to \gbigp_G\times V$
by $(x,v)\longmapsto (a(x),v)$.
It induces an isomorphism of vector bundles
$\Psi_a:X\times V\simeq \gbigp_G(\vecV)$.
Let $b$ be a section of $\gbigp^h_K$.
We obtain the map $\alpha:X\to G$
determined by $a=b\cdot\alpha$.
We can check the following lemma directly from the construction.
\begin{lem}
\label{lem;25.9.25.20}
For sections $u_1,u_2$ of $X\times V$,
we have 
$\Psi_a^{\ast} (h_E)(u_1,u_2)
=h_V(\alpha^{-1}u_1,\alpha^{-1}u_2)$.
\hfill\qed 
\end{lem}

\subsection{Holomorphic principal $G$-bundles}

\subsubsection{Holomorphic structures of principal $G$-bundles}

Let $G$ be a complex linear algebraic group.
Suppose that $X$ is a complex manifold.
Recall that a holomorphic principal $G$-bundle is a principal $G$-bundle
in the category of complex manifolds.
If $\gbigp_G$ is a holomorphic principal $G$-bundle on $X$,
$\gbigp_G(\vecV)$ is naturally a holomorphic vector bundle on $X$
for any $\vecV\in\Rep(G)$.
Moreover, the following holds.
\begin{itemize}
 \item For any $\vecV_1\to \vecV_2$ in $\Rep(G)$,
       the induced map
       $\gbigp_G(\vecV_1)\to\gbigp_G(\vecV_2)$ is holomorphic.
 \item For any $\vecV_1,\vecV_2\in\Rep(G)$,
       the natural isomorphisms
       $\gbigp_G(\vecV_1\oplus \vecV_2)\simeq
       \gbigp_G(\vecV_1)\oplus\gbigp_G(\vecV_2)$,
       $\gbigp_G(\vecV_1\otimes \vecV_2)
       \simeq
       \gbigp_G(\vecV_1)\otimes\gbigp_G(\vecV_2)$,
       and
       $\gbigp_G(\Hom(\vecV_1,\vecV_2))
       \simeq
       \Hom(\gbigp_G(\vecV_1),\gbigp_G(\vecV_2))$
       are holomorphic.
 \item $\gbigp_G(\cnum,\II)=\nbigo_X$.
\end{itemize}
The following proposition is a variant of Nori's representability theorem
\cite{Nori}.
\begin{prop}
\label{prop;24.5.27.30}
Let $\gbigp_G$ be a principal $G$-bundle on $X$ in the category of
$C^{\infty}$-manifolds.
Suppose that
the vector bundles $\gbigp_G(\vecV)$ $(\vecV\in\Rep(G))$
are equipped with holomorphic structures
satisfying the following conditions.
\begin{itemize}
 \item For any monomorphism $\vecV_1\to \vecV_2$ in $\Rep(G)$,
       the induced map
       $\gbigp_G(\vecV_1)\to\gbigp_G(\vecV_2)$ is holomorphic.
 \item For any $\vecV_1,\vecV_2\in\Rep(G)$,
       the natural isomorphisms
       $\gbigp_G(\vecV_1\oplus \vecV_2)\simeq
       \gbigp_G(\vecV_1)\oplus\gbigp_G(\vecV_2)$,
       $\gbigp_G(\vecV_1\otimes \vecV_2)
       \simeq
       \gbigp_G(\vecV_1)\otimes\gbigp_G(\vecV_2)$,
       and
       $\gbigp_G(\Hom(\vecV_1,\vecV_2))
       \simeq
       \Hom(\gbigp_G(\vecV_1),\gbigp_G(\vecV_2))$
       are holomorphic.
\item $\gbigp_G(\cnum,\II)=\nbigo_X$.
\end{itemize}
Then, $\gbigp_G$ has a unique complex structure
such that $\gbigp_G$ is a a holomorphic principal $G$-bundle
and that it induces the holomorphic structures of
$\gbigp_G(\vecV)$ $(\vecV\in\Rep(G))$.
\end{prop}
\pf
Let $\vecV$ be a faithful representation of $G$.
We set $E=\gbigp_G(\vecV)$
which is a holomorphic vector bundle.
Set $n=\dim V$,
and we choose an isomorphism $\GL(V)\simeq\GL_n(\cnum)$.
Let $\gbigp_E$ denote the frame bundle of $E$,
which is a holomorphic principal $\GL_n(\cnum)$-bundle
by the holomorphic structure of $E$.
There exists an isomorphism 
$\gbigp_G(\GL(V))\simeq\gbigp_E$
of principal $\GL_n(\cnum)$-bundles.
By this isomorphism,
we may regard $\gbigp_G$ as a closed submanifold of $\gbigp_E$.
Let us prove that
$\gbigp_G$ is a complex submanifold of $\gbigp_E$.

We set $W:=\End(V)\oplus \End(V^{\lor})$.
Let $\chi_W:G\to\GL(W)$
be the homomorphism defined by
$\chi_W(g)(f,f^{\lor})
=(\chi_{\vecV}(g)\circ f,\chi_{\vecV}(g^{-1})^{\lor}\circ f^{\lor})$.
We obtain the representation $\vecW=(W,\chi_W)$.

Let $\iota:G\to\GL(V)$ denote the inclusion.
Let $\varphi:\GL(V)\to \End(V)\times \End(V^{\lor})=W$
be the closed embedding defined by
$\varphi(A)=(A,(A^{\lor})^{-1})$.
We obtain the induced closed embedding of algebraic varieties
$\varphi_G=\varphi\circ\iota:G\to W$.
We consider the holomorphic
vector bundle $\Etilde=\gbigp_G(\vecW)$.
By $G\subset \GL(V)\subset W$,
we obtain
$\gbigp_G\subset\gbigp_E\subset \Etilde$.
Because $\gbigp_E$ is a complex submanifold of $\Etilde$,
it is enough to prove that
$\gbigp_G$ is a complex submanifold of $\Etilde$.

Let $\cnum[W]$ denote the ring of polynomial functions on $W$.
Let $\cnum[G]$ denote the ring of algebraic functions on $G$.
There exists the surjection
$\varphi_G^{\ast}:\cnum[W]\to\cnum[G]$.
Let $I_G\subset \cnum[W]$ denote the kernel ideal.
Let $\cnum[W]_{\leq k}\subset\cnum[W]$
denote the subspace of polynomial functions on $W$
of degree at most $k$.
We set $(I_G)_{\leq k}:=I_G\cap\cnum[W]_{\leq k}$.
If $k$ is sufficiently large,
$(I_G)_{\leq k}$ generates the ideal $I_G$.
For any generator $f_1,\ldots,f_N$ of
$(I_G)_{\leq k}$ over $\cnum$, we obtain
\[
 G=\bigcap_{i=1}^N  f_i^{-1}(0)\subset W.
\]

By the $G$-action on $W$,
we may regard $\cnum[W]_{\leq k}$ and
$(I_G)_{\leq k}$ are $G$-representations.
We obtain the following holomorphic vector bundles on $X$:
\[
 \gbigp_G\bigl(
 (I_G)_{\leq k}
 \bigr)
 \subset
 \gbigp_G\bigl(
 \cnum[W]_{\leq k}
 \bigr).
\]
We may regard $C^{\infty}$-sections of
$\gbigp_G\bigl(
 \cnum[W]_{\leq k}
 \bigr)$
as $C^{\infty}$-functions on $\Etilde$
whose restriction to any fibers are polynomials of degree at most $k$.
We may regard $C^{\infty}$-sections of
$\gbigp_G\bigl(
 (I_G)_{\leq k}
 \bigr)$
as $C^{\infty}$-functions on $\Etilde$
whose restriction to any fibers are polynomials of degree at most $k$
and vanishing on $\gbigp_G$.

We have the holomorphic structures of
$\gbigp_G\bigl(
 (I_G)_{\leq k}
 \bigr)$
and
$\gbigp_G\bigl(
 \cnum[W]_{\leq k}
 \bigr)$.
We may regard
holomorphic sections of
$\gbigp_G\bigl(\cnum[W]_{\leq k}\bigr)$
as holomorphic sections on $\Etilde$
whose restriction to any fibers are polynomials of degree
at most $k$.
We may regard holomorphic sections of
$\gbigp_G\bigl(
 (I_G)_{\leq k}
 \bigr)$
as holomorphic functions on $\Etilde$
whose restriction to any fibers are polynomials of degree at most $k$
and vanishing on $\gbigp_G$.

Let $P$ be any point of $X$.
There exist a neighbourhood $X_P$
and a holomorphic frame
$s_1,\ldots,s_{\ell}$ of
$\gbigp_G\bigl(
 (I_G)_{\leq k}
 \bigr)$.
We regard $s_i$ as holomorphic functions on
$\gbigp_{\Etilde|X_P}$.
Then,
\[
(\gbigp_G)_{|X_P}
=\bigcap_{i=1}^{\ell}s_i^{-1}(0)
=\Etilde_{|X_P}.
\]
Hence, $\gbigp_G$ is a closed complex analytic subset of $\gbigp_E$,
and hence a complex submanifold of $\Etilde$.
Because $\gbigp_G$ is a complex submanifold of $\gbigp_E$,
the complex structure of $\gbigp_G$
induces the holomorphic structure of $E$.

Let $\vecV_i$ $(i=1,2)$ be faithful $G$-representations.
Note that the isomorphism
$\gbigp_G(\vecV_1)\oplus\gbigp_G(\vecV_2)
\simeq
\gbigp_G(\vecV_1\oplus \vecV_2)$
is holomorphic.
We set $E_i=\gbigp_G(\vecV_i)$
and $E_0=\gbigp_G(\vecV_1\oplus \vecV_2)$.
Let $(\gbigp_G)_i\subset\gbigp_{E_i}$
denote the complex submanifolds
obtained as the image of $\gbigp_G$.
The image of the holomorphic embedding
$\gbigp_{E_1}\times_X\gbigp_{E_2}
\subset
 \gbigp_{E_0}$
contains $(\gbigp_G)_0$.
The holomorphic projections
$\gbigp_{E_1}\times_X\gbigp_{E_2}\to\gbigp_{E_i}$
maps $(\gbigp_G)_0$ onto $(\gbigp_G)_i$.
Hence,
$(\gbigp_G)_i$ are naturally isomorphic
as complex manifolds.
Hence, 
we obtain that
the complex structure of $\gbigp_G$
is independent of the choice of
a faithful representation $\vecV$.

Let $\vecU\in\Rep(G)$.
By taking a faithful representation $\vecV\in\Rep(G)$,
we obtain a faithful representation $\vecU\oplus \vecV$.
Because
$\gbigp_G(\vecU)\subset\gbigp_G(\vecU\oplus \vecV)$
is holomorphic,
it is easy to see that holomorphic structure of
$\gbigp_G(\vecU)$ is induced by the complex structure of
$\gbigp_G$.
\hfill\qed

\vspace{.1in}
Let $A^{p,q}(\gbigp_G(\gminig))$
denote the space of $\gbigp_G(\gminig)$-valued $(p,q)$-forms.

\begin{cor}
\label{cor;25.6.10.1}
Let $\omega$ be a connection of $\gbigp_G$.
Suppose that the curvature $R(\omega)$ is contained 
in $A^{1,1}(\gbigp_G(\gminig))$.
\begin{itemize}
 \item For any $\vecV\in\Rep(G)$,
       let $(\nabla^{\omega}_{\vecVtilde})^{0,1}$
       denote the $(0,1)$-part of
       the induced connection $\nabla^{\omega}_{\vecVtilde}$ of
       the vector bundle
       $\vecVtilde=\gbigp_G(\vecV)$.
       Then,
       $(\nabla^{\omega}_{\vecVtilde})^{0,1}
       \circ(\nabla^{\omega}_{\vecVtilde})^{0,1}=0$ holds.
       As a result,
       $\vecVtilde$ is equipped with the holomorphic structure
       defined by 
       $(\nabla^{\omega}_{\vecVtilde})^{0,1}$.
 \item $\gbigp_G$ has a unique complex structure
       such that
       $\gbigp_G$ is a holomorphic principal $G$-bundle
       and that
       it induces the holomorphic structures of
       $\vecVtilde$ $(\vecV\in\Rep(G))$.
       \hfill\qed
\end{itemize}
\end{cor}

\subsubsection{$K$-reductions and holomorphic structures}
\label{subsection;25.6.9.22}

Suppose moreover that $G$ is the complexification of
its compact maximal subgroup $K$.
If $G$ is connected and reductive,
this condition is satisfied.
Let $\rho$ denote the Cartan involution.
Let $h$ be a $K$-reduction of $\gbigp_G$.

Let $\vecV\in\Rep(G)$ be a representation.
Let $h_V$ be a $K$-invariant Hermitian metric of $V$.
It induces a Hermitian metric
$\htilde_V$ of $\gbigp_G(\vecV)$.
\begin{lem}
\label{lem;25.6.9.11}
The Chern connection of $(\gbigp_G(\vecV),\htilde_V)$ 
is independent of the choice of
a $K$-invariant Hermitian metric $h_V$.
\end{lem}
\pf
Let $h'_V$ be another $K$-invariant Hermitian metric of $V$.
There exists the $K$-invariant decomposition
\begin{equation}
\label{eq;25.6.9.10}
V=\bigoplus_{j=1}^mV_j
\end{equation}
such that
(i) the decomposition (\ref{eq;25.6.9.10}) is orthogonal with respect to 
both $h_V$ and $h_V'$,
(ii) $h'_{V|V_j}=c_jh_{V|V_j}$ for some $c_j>0$.
Note that the decomposition (\ref{eq;25.6.9.10}) is also $G$-invariant,
i.e.,
we obtain the decomposition
$\vecV=\bigoplus \vecV_j$ in $\Rep(G)$.
Let $\htilde_V$ and $\htilde'_V$
denote the induced Hermitian metrics of $\gbigp_G(\vecV)$.
The decomposition
\[
\gbigp_G(\vecV)=\bigoplus_{j=1}^m\gbigp_K(\vecV_j)
\]
is holomorphic and orthogonal with respect to
both $\htilde_V$ and $\htilde_{V'}$.
Moreover,
$\htilde'_{V|\gbigp_K(V_j)}
=c_j\htilde_{V|\gbigp_K(V_j)}$.
Then, the Chern connections
of the induced Hermitian metrics are the same.
\hfill\qed

\begin{prop}
\label{prop;25.6.9.20}
There exists
a unique connection $\omega_h$ of $\gbigp^h_K$
which induces 
the Chern connections of
$(\gbigp_G(V),\htilde_V)$
for any $V\in\Rep(G)$
and any $K$-invariant Hermitian metric $h_V$ of $V$. 
The connection $\omega_h$ is called
the Chern connection of
$\gbigp_G$ with the $K$-reduction $\gbigp^h_K$.
\end{prop}
\pf
By Proposition \ref{prop;24.5.30.2}
and Lemma \ref{lem;25.6.9.11},
there exists a connection $\omega_h$ of $\gbigp_G$
which induces the Chern connections of
$(\gbigp_G(\vecV),\htilde_V)$ $(\vecV\in\Rep(G))$.
By Lemma \ref{lem;24.5.30.5},
it is induced by a connection of $\gbigp^h_K$.
\hfill\qed

\subsection{Harmonic $G$-bundles}
\label{subsection;26.3.9.3}

Let $X$ be a complex manifold.
Let $G$ be a reductive complex linear group
with a compact real form $K$ and the Cartan involution $\sigma$.
We shall use the convention in Notation \ref{notation;25.9.23.1}
to describe a connection of a principal $G$-bundle.

\subsubsection{Higgs case}

Let $\gbigp_G$ be a holomorphic principal $G$-bundle.
For any section $\tau_i$ $(i=1,2)$
of $\gbigp_G(\gminig)\otimes\Omega^{p(i),q(i)}$,
we obtain
the section $\tau_1\wedge\tau_2$
of $\gbigp_G(\gminig)\otimes\Omega^{p(1)+p(2),q(1)+q(2)}$
induced by the Lie bracket on
$\gbigp_G(\gminig)$
and the exterior product of $\Omega^{\bullet,\bullet}$.

\begin{df}
A holomorphic section $\theta$
of $\gbigp_G(\gminig)\otimes\Omega^{1,0}_X$ 
is called a Higgs field 
$\gbigp_G$
if $\theta\wedge\theta=0$.
Such a pair $(\gbigp_G,\theta)$ is called a $G$-Higgs bundle.
\hfill\qed
\end{df}

Let $h$ be a $K$-reduction of $\gbigp_G$.
We obtain the Chern connection $\nabla_h$ of $\gbigp_G$
in Proposition \ref{prop;25.6.9.20}.
We also obtain
$\theta^{\dagger}_h
=-\rho_h(\theta)\in A^{0,1}(\gbigp_G(\gminig))$.
We obtain another connection
$\DD^1_h=\nabla_h+\theta+\theta^{\dagger}_h$ of $\gbigp_G$.
\begin{df}
$h$ called is a pluri-harmonic reduction of
a $G$-Higgs bundle $(\gbigp_G,\theta)$
if $\DD^1_h$ is a flat connection.  
In that case,
$(\gbigp_G,\theta,h)$ is called a harmonic $G$-bundle.
\hfill\qed
\end{df}

The following lemma is clear.
\begin{lem}
The following conditions are equivalent.
\begin{itemize}
 \item $(\gbigp_G,\theta,h)$ is a harmonic $G$-bundle.
 \item $(\gbigp_G(\vecV),\chi_{\vecV}(\theta),\htilde_V)$
       is a harmonic bundle
       for any $\vecV\in\Rep(G)$ with a $K$-invariant Hermitian
       metric $h_V$ of $V$.
 \item $(\gbigp_G(\vecV),\chi_{\vecV}(\theta),\htilde_V)$
       is a harmonic bundle
       for a faithful $G$-representation $\vecV$
       with a $K$-invariant Hermitian metric $h_V$ of $V$.
\end{itemize}
(See {\rm\S\ref{subsection;25.6.9.22}} for
the metric $\htilde_V$.)
\hfill\qed
\end{lem}

\subsubsection{Flat case}

Let $\gbigp_G$ be a principal $G$-bundle
with a flat connection $\nabla$.
Let $h$ be a $K$-reduction of $\gbigp_G$.
The connection $\nabla$ is decomposed as
$\nabla=\nabla_h+\Phi_h$,
where $\nabla_h$ is induced by the connection of $\gbigp^h_K$,
and $\Phi_h\in A^{1}(\gbigp^h_K(\sqrt{-1}\gminik))$.
(See \S\ref{subsection;25.6.9.23}.)
There exists the decomposition
$\Phi_h=\theta_h+\theta^{\dagger}_h$
where
$\theta_h\in A^{1,0}(\gbigp_G(\gminig))$
and
$\theta^{\dagger}_h\in A^{0,1}(\gbigp_G(\gminig))$.
We obtain the connection
$\nabla_{h,\gbigp_G(\gminig)}$
of the vector bundle $\gbigp_G(\gminig)$ induced by $\nabla_h$,
which is decomposed as
$\nabla_{h,\gbigp_G(\gminig)}
=\nabla_{h,\gbigp_G(\gminig)}^{1,0}
+\nabla_{h,\gbigp_G(\gminig)}^{0,1}$
into the $(1,0)$-part and the $(0,1)$-part.

\begin{df}
$h$ is called a pluri-harmonic reduction
of $(\gbigp_G,\nabla)$
 if $\nabla_{h,\gbigp_G(\gminig)}^{0,1}\theta_h=0$.
 In the case,
$(\gbigp_G,\nabla,h)$  is called a harmonic $G$-bundle.
 \hfill\qed
\end{df}

Let $h$ be a pluri-harmonic reduction of $(\gbigp_G,\nabla)$.
Let $\vecV\in\Rep(G)$,
and we set $\vecVtilde=\gbigp_G(\vecV)$.
Let $\nabla_{\vecVtilde}$ be the induced flat connection.
Let $h_V$ be a $K$-invariant Hermitian metric of $V$,
which induces a Hermitian metric $\htilde_V$ of $\vecVtilde$.
We obtain the unique decomposition
$\nabla_{\vecVtilde}=\nabla_{\htilde_{V}}+\Phi_{\htilde_{V}}$,
where $\nabla_{\htilde_{V}}$
denotes a unitary connection of $(\vecVtilde,\htilde_V)$,
and $\Phi_{\htilde_V}\in A^1(\End(\vecVtilde))$
is self-adjoint
with respect to $\htilde_{V}$.
We also obtain the decomposition
$\nabla_{\htilde_V}=\del_{\vecVtilde}+\delbar_{\vecVtilde}$
and
$\Phi_{\htilde_V}=\theta_{\vecVtilde}+\theta_{\vecVtilde}^{\dagger}$
into the $(1,0)$-part and the $(0,1)$-part.

\begin{lem}
\label{lem;24.5.27.40}
$(\vecVtilde,\delbar_{\vecVtilde},\theta_{\vecVtilde})$
is a Higgs bundle.
\end{lem}
\pf
By the condition
$\nabla^{0,1}_{h,\gbigp_G(\gminig)}\theta_h=0$,
we obtain
$\delbar_{\vecVtilde}(\theta_{\vecVtilde})=0$.
Then, we obtain
$\delbar_{\vecVtilde}\circ\delbar_{\vecVtilde}=0$
and $\theta_{\vecVtilde}\wedge\theta_{\vecVtilde}=0$
as explained in \cite{Mochizuki-KH-Higgs}.
\hfill\qed

\begin{lem}
The tuple $(\gbigp_G,\theta_h)$ is a $G$-Higgs bundle.
\end{lem}
\pf
By Proposition \ref{prop;24.5.27.30},
we obtain the complex structure of $\gbigp_G$.
We obtain $\theta_h\wedge\theta_h=0$ from Lemma \ref{lem;24.5.27.40}
by taking a faithful representation $\vecV\in\Rep(G)$.
Because $\theta_{\vecVtilde}$ 
are holomorphic sections of
$\End(\vecVtilde)\otimes\Omega^{1,0}$
for a faithful $G$-representation $\vecV$,
we obtain that $\theta$ is holomorphic section of
$\gbigp_G(\gminig)\otimes\Omega^{1,0}$.
\hfill\qed

\begin{rem}
Let $(\gbigp_G,\theta_G)$
be a $G$-Higgs bundle
with  
a pluri-harmonic reduction $h$.
We obtain the flat connection $\DD^1_h$ of $\gbigp_G$,
and $h$ is a pluri-harmonic reduction of
$(\gbigp_G,\DD^1_h)$.
Conversely,
let $\gbigp_G$ be a principal $G$-bundle
with a flat connection $\nabla$
and  a pluri-harmonic reduction $h$.
We have  the complex structure of
$\gbigp_G$ determined by $\nabla_h^{0,1}$,
$\theta_h$ is a Higgs field
of the holomorphic principal $G$-bundle $\gbigp_G$,
and $h$ is a pluri-harmonic reduction of
$(\gbigp_G,\theta_h)$.
The constructions are mutually inverse.
In this sense,
two definitions of harmonic $G$-bundles are equivalent.
\hfill\qed
\end{rem}

\subsection{Hermitian metrics of principal bundles in the semisimple case}

\subsubsection{Hermitian metrics in the semisimple case}

Suppose that $G$ is semisimple complex algebraic group.
Let $\gbigp_G$ be a principal $G$-bundle.
We set $\Herm(\gbigp_G)=\gbigp_G(\Herm(G))$.
\begin{df}
\label{df;25.6.17.30}
A section $h$ of $\Herm(\gbigp_G)$ is called
a Hermitian metric of  $\gbigp_G$.
\hfill\qed
\end{df}

Let $K$ be a maximal compact subgroup of $G$.
The left $G$-action on $\Herm(G)$
induces a bijection $G/K\to\Herm(G)$
by $g\mapsto gKg^{-1}$.
It induces
\[
\gbigp_G(G/K)\simeq \Herm(\gbigp_G).
\]
Hence, Hermitian metrics of $\gbigp_G$
are equivalent to $K$-reductions of $\gbigp_G$.

For any Hermitian metric $h$ of $\gbigp_G$,
the corresponding principal $K$-bundle
is denoted by $\gbigp^h_K$.
The Hermitian metric $h_{\gminig,K}$ on $\gminig$
induces a Hermitian metric on
the vector bundle $\gbigp_G(\gminig)=\gbigp^h_{K}(\gminig)$,
which is also denoted by $h$.
By the naturally induced action of $K$ on $\exp(\sqrt{-1}\gminik)$,
we obtain the bundle
$\gbigp^h_{K}(\exp(\sqrt{-1}\gminik))$ over $X$.
There exists the natural isomorphism
\[
 \gbigp^h_{K}(\exp(\sqrt{-1}\gminik))
 \simeq
 \gbigp^h_{K}(\Herm(G))
 \simeq
 \gbigp_{G}(\Herm(G))
=\Herm(\gbigp_G).
\]
We also have
\[
 \gbigp^h_{K}(\exp(\sqrt{-1}\gminik))
 \subset
 \gbigp^h_{K}(G,\Ad)
 =\gbigp_G(G,\Ad).
\]

Let $h_1,h_2$ be Hermitian metrics of $\gbigp_G$.
By Lemma \ref{lem;25.5.7.1},
we obtain
\[
s(h_1,h_2)\in
\gbigp^{h_1}_{K}(\exp(\sqrt{-1}\gminik))
\subset
\gbigp_G(G,\Ad)
\]
as the difference of $h_1$ and $h_2$.
We define
\[
 v(h_1,h_2)\in\Gamma(\gbigp^{h_1}_{K}(\sqrt{-1}\gminik))
 \subset
 \Gamma(\gbigp_G(\gminig))
\]
by the condition
$\exp(v(h_1,h_2))=s(h_1,h_2)$.
\begin{lem}
We have $s(h_2,h_1)=s(h_1,h_2)^{-1}$
in $\Gamma(\gbigp_G(G,\Ad))$,
and $v(h_2,h_1)=-v(h_1,h_2)$
in $\Gamma(\gbigp_G(\gminig))$.
\hfill\qed
\end{lem}

\subsubsection{Compatibility with an automorphism of finite order}

Let $\sigma$ be a complex algebraic automorphism of
the complex algebraic group $G$
of finite order.
Let $H$ be a complex algebraic closed subgroup of $G$
such that
$G_0^{\sigma}\subset H\subset G_Z^{\sigma}$.
Let $\gbigp_{H}\subset \gbigp_G$ be an $H$-reduction.
We set
$\Herm(\gbigp_G,\sigma)=\gbigp_{H}(\Herm(G,\sigma))$.
We have the natural inclusion
$\Herm(\gbigp_G,\sigma)\subset \Herm(\gbigp_G)$.

\begin{df}
A Hermitian metric of $\gbigp_G$ is called compatible with $\sigma$
if it is a section of $\Herm(\gbigp_G,\sigma)$.
\hfill\qed 
\end{df}

Let $K\in\Herm(G,\sigma)$.
We set $K_{H}=K\cap H$.
\begin{lem}
Suppose $\Herm(G,\sigma)\neq\emptyset$,
then Hermitian metrics of $\gbigp_G$ compatible with $\sigma$
bijectively correspond to
$K_{H}$-reductions of $\gbigp_{H}$.
\end{lem}
\pf
Because $H/K_H\simeq \Herm(G,\sigma)$
by Corollary \ref{cor;25.6.17.10},
we obtain
\[
 \gbigp_H(H/K_{H})
 \simeq
 \gbigp_H(\Herm(G,\sigma))
 \simeq
 \Herm(\gbigp_G,\sigma).
\]
Then, we obtain the claim of the lemma.
\hfill\qed

\vspace{.1in}

Suppose moreover $H\subset G^{\sigma}$.
Then, the automorphism $\id\times\sigma$
of $\gbigp_H\times \Herm(G)$
induces an automorphism
of $\Herm(\gbigp_{G})=\gbigp_{H}(\Herm(G))$.
The induced automorphism is also denoted by $\sigma$.
Because the $\sigma$-invariant part of
$\Herm(\gbigp_G)$
equals
$\gbigp_H(\Herm(G,\sigma))$,
the following lemma is clear by definition.
\begin{lem}
\label{lem;25.6.16.3}
A Hermitian metric $h$ of $\gbigp_G$
is compatible with $\sigma$
if and only if it is $\sigma$-invariant.
\hfill\qed
\end{lem}

\subsection{Some estimates in the split case}

\subsubsection{Canonical Hermitian metrics associated with
regular semisimple sections}
\label{subsection;25.6.18.30}

Let $G$ be a semisimple complex Lie group.
Let $\sigma:G\to G$ be a complex algebraic homomorphism
of finite order $m$.
Let $\omega$ be a primitive $m$-th root of $1$.
We obtain the decomposition (\ref{eq;25.5.13.11}).
We assume that $(\sigma,\omega)$ is split.

Let $H$ be a complex algebraic closed subgroup of $G$
such that
$G_0^{\sigma}\subset H\subset G^{\sigma}_Z$.
Let $\gbigp_H$ be a principal $H$-bundle on $X$.
We set $\gbigp_G=\gbigp_H(G)$.
The adjoint action of $H$ on $\gminig_1$
induces
$\gbigp_{H}(\gminig_1)\subset \gbigp_G(\gminig)$.
It contains an open subset
$\gbigp_{H}(\gminig_1^{\rs})$.

Let $u$ be any section of
$\gbigp_{H}(\gminig_1^{\rs})$.
We obtain the following lemma from
Corollary \ref{cor;25.6.21.10}.

\begin{lem}
\mbox{{}}\label{lem;25.6.20.2}
\begin{itemize}
\item There exists
      a unique Hermitian metric $h^{u}$
      of $\gbigp_G$ compatible with $\sigma$
      such that
      $\bigl[u,\rho_{h^{u}}(u)\bigr]=0$.
 \item $h(u,u)\geq h^{u}(u,u)$
       for any Hermitian metric $h$ of $\gbigp_G$
       compatible with $\sigma$.
       The equality holds
       if and only if $h=h^u$.
       \hfill\qed
\end{itemize}
\end{lem}

\subsubsection{Estimates}

Let $u$ be any section of
$\gbigp_{H}(\gminig_1^{\rs})$.
Let $A\subset X$ be any compact subset.
We obtain the following proposition
from Proposition \ref{prop;25.5.13.20}.
\begin{prop}
\label{prop;25.5.13.21}
For any $C_1>0$,
there exists a compact subset
$\nbigk\subset\Herm(\gbigp_G,\sigma)_{|A}$ 
depending on $(u,A,C_1)$
such that the following holds.
\begin{itemize}
 \item Let $h$ be any Hermitian metric of $\gbigp_{G|A}$
       compatible with $\sigma$
       such that
       $h(u_{|A},u_{|A})\leq C_1$.
       Then,
       $v(h^u,h)(A)\subset\nbigk$.
\end{itemize}
As a result, for any $C_1>0$, there exists
$C_2>0$ depending on $(u,A,C_1)$
such that the following holds.
\begin{itemize}
 \item Let $h$ be any Hermitian metric of $\gbigp_{G|A}$
       compatible with $\sigma$
       such that
       $h(u_{|A},u_{|A})\leq C_1$.
       Then, $|v(h^{u}_{|A},h)|_{h^{u}}\leq C_2$ holds.
       \hfill\qed
\end{itemize}
\end{prop}

We obtain the following proposition from
Corollary \ref{cor;25.5.15.5}.
\begin{prop}
\label{prop;25.6.18.50}
For any $r_1>0$,
there exists $C>0$ 
such that the following holds.
\begin{itemize}
 \item Let $h$ be 
       any Hermitian metric of $\gbigp_{G|A}$ compatible with $\sigma$
       such that 
       $h(u_{|A},u_{|A})-h^u(u,u)_{|A}\leq r_1$.
       Then,
       $\bigl|
       v(h^u_{|A},h)
       \bigr|_{h^u_{|A}}
       \leq
       C
       \bigl|
       [u_{|A},\rho_{h}(u_{|A})]
       \bigr|_{h}$
       holds.
       \hfill\qed
\end{itemize}
\end{prop}

\section{Basic results for harmonic $G$-bundles on Riemann surfaces}

\subsection{Local estimate of the norm of Higgs fields}
\label{subsection;25.9.24.10}

Let $G$ be a reductive complex algebraic group
with the Lie algebra $\gminig$.
For $R>0$,
we set $U(R)=\bigl\{z\in\cnum\,\big|\,|z|<R\bigr\}$.
We consider the principal $G$-bundle $\gbigp_G=U(R)\times G$ on $U(R)$.
We have
$\gbigp_G(\gminig)=U(R)\times\gminig$.
Let $\theta=f\,dz$ be
a holomorphic section of $\gbigp_G(\gminig)\otimes\Omega^1$.

Let $\gminit$ be a Cartan subalgebra of $\gminig$
with the Weyl group $W(\gminit)$.
Let $h_{\gminit}$ be a $W(\gminit)$-invariant Hermitian metric
of $\gminit$.
Let $|u|_{h_{\gminit}}$ denote the norm of $u\in\gminit$
with respect to $h_{\gminit}$.
It induces a map
$|\cdot|_{h_{\gminit}}:\gminit/W(\gminit)\to\real_{\geq 0}$.

Let $\Psi_{f}:U(R)\to \gminit/W(\gminit)$ denote the morphism
induced by $f$.

We directly obtain the following from Simpson's main estimate
\cite[Proposition 2.1]{Decouple}.
 \begin{prop}
\label{prop;25.5.13.50}
Suppose that 
$|\Psi_f|_{h_{\gminit}}$ is bounded on $U(R)$.
Then, for any $G$-representation $\vecV$
with a $K$-invariant Hermitian metric $h_V$,
and any $0<R_1<R$,
there exists a positive constant
$C_i$ $(i=1,2)$ depending on
the tuple
 $(h_{\gminit},R,R_1,V,\chi_V,h_V)$ such that
for any harmonic reduction $h$ of $(\gbigp_G,\theta)$,
the inequality
\[
       |\chi_V(f)|_{(h,h_V)}\leq
       C_1\sup_{U(R)}|\Psi_f|_{h_{\gminit}}+C_2
\]
  holds on $U(R_1)$.
  Here, $|\cdot |_{(h,h_V)}$
  denotes the norm induced by
  $h$ and $h_V$.
\hfill\qed
 \end{prop}

\subsection{Uniqueness of harmonic metrics under some condition}

Let $X$ be any Riemann surface.
Let $Y\subset X$ be a connected and relatively compact open subset
with smooth non-empty boundary $\del Y$.
Let $D\subset Y$ be a finite subset.
For any $P\in D$,
let $(X_P,z_P)$ be a holomorphic coordinate neighbourhood
around $P$
such that
(i) $X_P$ is relatively compact in $Y$,
(ii) $\Xbar_P\cap \Xbar_{P'}=\emptyset$ $(P\neq P')$.
By the coordinate $z_P$,
we assume $X_P\simeq\{|z|<1/2\}$.

\subsubsection{The case of Higgs bundles}

Let $(E,\theta)$ be a Higgs bundle on $X\setminus D$
in the sense that $E$ is a holomorphic vector bundle,
and $\theta$ is a holomorphic section of $\End(E)\otimes\Omega^1$.

\begin{prop}
\label{prop;25.6.18.1}
Let $h_1,h_2$ be harmonic metrics of
$(E,\theta)_{|Y\setminus D}$
such that the following holds.
\begin{itemize}
 \item $h_{1|\del Y}=h_{2|\del Y}$.
 \item $\log \Tr(s(h_1,h_2))=O\bigl(\log(-\log|z_P|)\bigr)$
       for any $P\in D$.
\end{itemize}
Then, we obtain $h_1=h_2$.
\end{prop}
\pf
First, let us consider the case $\rank E=1$.
Then, $\log s(h_1,h_2)$ is a harmonic function on $Y$
such that $\log s(h_1,h_2)_{|\del Y}=0$.
We have
$\log s(h_1,h_2)=O\bigl(\log(-\log|z_P|)\bigr)$
on $X_P\setminus\{P\}$ for any $P\in D$.
For any $\epsilon>0$,
$\log s(h_1,h_2)
+\epsilon\log|z_P|$
is bounded from above on $X_P\setminus\{P\}$,
and hence 
$\log s(h_1,h_2)
+\epsilon\log|z_P|$ is subharmonic on $X_P$.
By the maximum principle,
there exists $C_{10}>0$
depending only on
$\log s(h_1,h_2)_{|\del X_P}$
such that
$\log s(h_1,h_2)
+\epsilon\log|z_P|\leq C_{10}$
for any $\epsilon>0$.
As a result,
$\log s(h_1,h_2)$ is bounded from above.
Moreover, the maximum principle holds on $Y$.
Because $\log s(h_1,h_2)_{|\del Y}=0$,
we obtain that $s(h_1,h_2)=1$.

Let us consider the general rank case.
We have $\det s(h_1,h_2)=1$
by the consideration in the rank one case.
In particular, we have
$\Tr s(h_1,h_2)\geq \rank E$.
The function $\Tr(s(h_1,h_2))$ on $Y\setminus D$
is subharmonic.
By a similar argument,
we obtain that $\Tr s(h_1,h_2)$ is subharmonic on $Y$.
In particular, the maximum principle holds on $Y$.
Because $\Tr s(h_1,h_2)=\rank E$ on $\del Y$,
we obtain $\Tr s(h_1,h_2)=\rank E$ on $Y$,
which implies $s(h_1,h_2)=\id$.
\hfill\qed

\subsubsection{The case of $G$-Higgs bundles}

Let $G$ be a complex reductive linear algebraic group
with a maximal compact subgroup $K$.
Let $(\gbigp_G,\theta)$ be a holomorphic $G$-Higgs bundle on
$X\setminus D$.

\begin{prop}
\label{prop;25.5.6.1}
Let $h_1,h_2$ be harmonic $K$-reductions of
$(\gbigp_G,\theta)_{|Y\setminus D}$
such that the following holds.
\begin{itemize}
 \item $h_{1|\del Y}=h_{2|\del Y}$.
 \item For a faithful $G$-representation
       $\vecV$,
       we obtain
       $\log \Tr\bigl(\chi_V (s(h_1,h_2))\bigr)
       =O\bigl(\log(-\log|z_P|)\bigr)$
       for any $P\in D$.
\end{itemize}
 Then, we obtain $h_1=h_2$.
\end{prop}
\pf
It is enough to consider the Higgs bundles
$(\gbigp_G(V),\chi_V(\theta))$ with
the induced harmonic metrics.
\hfill\qed 

\subsection{Reduction}

\subsubsection{Centralizers of semisimple elements}

Let $G$ be a connected reductive complex linear algebraic group.
Let $\gminig$ denote the Lie algebra.
For any $u\in\gminig$,
we set
$C_G(u)=\bigl\{
a\in G\,\big|\,
\Ad(a)u=u
\bigr\}$.

\begin{lem}
\label{lem;25.6.18.2}
If $u$ is semisimple,
then $C_G(u)$ is a connected reductive complex linear algebraic group.
\end{lem}
\pf
Let $A\subset G$ denote the Zariski closure of
the subset $\exp\bigl(\cnum\cdot u\bigr)$.
We have $C_G(u)=C_G(A)$.
Let $\gminit$ be a Cartan subalgebra of $\gminig$
which contains $u$.
Let $T$ be the corresponding maximal complex torus of $G$.
Then, $A\subset T$.
Hence, $A$ is a complex torus.
According to \cite[Theorem 11.31]{Nishiyama-Ohta},
$C_G(A)$ is connected.
According to \cite[Theorem 10.20]{Nishiyama-Ohta},
$C_G(A)$ is reductive.
\hfill\qed

\vspace{.1in}
Note that $u$ is not necessarily regular
in Lemma \ref{lem;25.6.18.2}.

\subsubsection{Maximal compact subgroups}

Let $K\subset G$ be a maximal compact subgroup.
Let $\gminik$ denote the Lie algebra of $K$.
Let $u\in\sqrt{-1}\gminik$.
Because $u$ is semisimple,
$C_G(u)\subset G$ is connected and reductive
by Lemma \ref{lem;25.6.18.2}.

\begin{lem}
\label{lem;25.9.24.2}
$K\cap C_G(u)$ is a compact real form of 
the complex reductive group $C_G(u)$,
and we have the Cartan decomposition
\[
 (K\cap C_G(u))\times \sqrt{-1}(\gminik\cap C_{\gminig}(u))
\simeq C_G(u).
\]
In particular,
$K\cap C_G(u)$ is a maximal compact subgroup
of $C_G(u)$.
\end{lem}
\pf
We set $\gminih=C_{\gminig}(u)$,
which is the Lie algebra of $H=C_G(u)$.
Because $\rho_{\gminik}(u)=-u$,
we obtain
$\rho_{\gminik}(\gminih)=\gminih$,
and
\[
 \gminih=(\gminih\cap\gminik)\oplus\sqrt{-1}(\gminih\cap\gminik).
\]
Because $\rho_{\gminik}(\gminih)=\gminih$,
we obtain
$\rho_K(H)=H$.

Let $g\in H$.
We have the unique expression $g=a\exp(v)$
by some $a\in K$ and $v\in\sqrt{-1}\gminik$.
Because $\rho_K(g)\in H$,
we obtain that
$\exp(2v)\in H$.
We obtain
$\Ad\exp(2v)(\sqrt{-1}u)=\sqrt{-1}u$,
which implies
\[
 \exp(2v)\exp(\sqrt{-1}tu)\exp(-2v)=\exp(\sqrt{-1}tu)
\]
for any $t\in\cnum$.
We obtain
$\exp\Bigl(
 \Ad\bigl(
 \exp(-\sqrt{-1}tu)\bigr)2v
 \Bigr)
=\exp(2v)$.
It implies
$\Ad\bigl(\exp(-\sqrt{-1}tu)\bigr)2v 
 =2v$.
We obtain that
$[u,2v]=0$,
i.e.,
$v\in \sqrt{-1}(C_{\gminig}(u)\cap \gminik)$.
We also obtain that $a\in K\cap C_G(u)$.
\hfill\qed

\subsubsection{Reductions associated with holomorphic and real semisimple sections}

Let $\gminit$ be a Cartan subalgebra of $\gminig$.
We have the real structure
$\gminit_{\real}\subset\gminit$.
We obtain the real subanalytic subset
$\gminit_{\real}/W(\gminit)
\subset
\gminit/W(\gminit)\simeq\cnum^{\dim\gminit}$.
Let $\gbigp_G$ be a holomorphic principal $G$-bundle on $X$.
Let $u$ be a holomorphic section of $\gbigp_G(\gminig)$.
It induces a map
$\Phi_u:X\to \gminit/W(\gminit)$.

\begin{lem}
\label{lem;25.9.24.3}
Assume that the image of $\Phi_u$ is contained in
$\gminit_{\real}/W(\gminit)$.
Then, $\Phi_u$ is constant.
\end{lem}
\pf
Let $Q$ be any polynomial function on $\gminit$.
By setting
\[
 \Psi(Q)(x):=\frac{1}{|W(\gminit)|}
 \sum_{a\in W(\gminit)}Q(a\cdot x),
\]
we obtain the $W(\gminit)$-invariant polynomial function
on $\gminit$.
It induces a $\cnum$-valued polynomial function
$\Psitilde(Q)$ on $\gminit/W(\gminit)$.
We obtain the holomorphic function $\Psitilde(Q)\circ\Phi_u$ on $X$.

We have the decomposition $Q=Q_1+\sqrt{-1}Q_2$,
where $Q_i$ are induced by
$\real$-valued polynomial functions on $\gminit_{\real}$.
We obtain the decomposition
$\Psitilde(Q)\circ\Phi_u
=\Psitilde(Q_1)\circ\Phi_u
+\sqrt{-1}\Psitilde(Q_2)\circ\Phi_u$.
Because the image of $\Phi_u$ is contained in $\gminit_{\real}/W(\gminit)$,
$\Psitilde(Q_j)\circ\Phi_u$
are $\real$-valued.
We obtain that $\Psitilde(Q_j)\circ\Phi_u$ are constant on $X$.
Hence, $\Psitilde(Q)\circ\Phi_u$ is constant
for any polynomial function $Q\in \gminit$.
It implies that $\Phi_u$ is constant.
\hfill\qed

\vspace{.1in}

We choose $a_0\in \gminit_{\real}$
such that
the induced point $[a_0]\in\gminit/W(\gminit)$
equals the image of $\Phi_u$.
We have the orbit
$G\cdot a_0\subset\gminig$.
It induces a subspace
$\gbigp_G(G\cdot a_0)
 \subset
 \gbigp_G(\gminig)$.

If $u_{|P}$ are semisimple for any $P\in X$,
then $u$ is a section of
$\gbigp_G(G\cdot a_0)
 =\gbigp_G(G/C_G(a_0))$.
It induces a $C_G(a_0)$-reduction
$\gbigp^u_{C_G(a_0)}\subset \gbigp_G$,
which is a holomorphic principal $C_G(a_0)$-bundle.
Note that
$C_G(a_0)$ is a connected reductive complex linear algebraic
subgroup of $G$.
The following lemma is clear.
\begin{lem}
\label{lem;25.9.24.1}
Let $L$ be any holomorphic line bundle on $X$.
If $f$ is a section of $\gbigp_G(\gminig)\otimes L$
such that $[u,f]=0$,
then $f$ is induced by a section of
d $\gbigp^u_{C_G(a_0)}(C_{\gminig}(a_0))\otimes L$.
\hfill\qed
\end{lem}

\subsubsection{Reductions}

Let $G$ be a connected reductive complex linear algebraic group.
Let $(\gbigp_G,\theta)$ be a $G$-Higgs bundle on $X$.

\begin{lem}
\label{lem;25.6.18.3}
Let $h_i$ $(i=1,2)$ be harmonic metrics of
$(\gbigp_G,\theta)$
such that the following holds
for a faithful $G$-representation $\vecV$:
\begin{itemize}
 \item $\chi_V(s(h_1,h_2))$ is holomorphic,
       and 
       $[\chi_V(\theta),\chi_V (s(h_1,h_2))]=0$.
 \end{itemize}
Then, the following holds.
\begin{itemize}
 \item $[v(h_1,h_2)]:X\to \gminit/W(\gminit)$
 is constant,
  and contained in $\gminit_{\real}/W(\gminit)$.
 \item Let $a_0\in\gminit_{\real}$ be a representative
 of $[v(h_1,h_2)]$.
       Then, we obtain a reduction
       $(\gbigp^{v(h_1,h_2)}_{C_G(a_0)},\theta)$
       of $(\gbigp_G,\theta)$.
       Moreover,
       $h_i$ are induced by
       harmonic metrics of
       the $C_G(a_0)$-Higgs bundle
       $(\gbigp^{v(h_1,h_2)}_{C_G(a_0)},\theta)$.
\end{itemize}
\end{lem}
\pf
By Lemma \ref{lem;25.9.24.3},
we obtain the first claim.
By Lemma \ref{lem;25.6.18.2}, 
and Lemma \ref{lem;25.9.24.1},
we obtain the $C_G(a_0)$-Higgs bundle
$(\gbigp^{v(h_1,h_2)}_{C_G(a_0)},\theta)$
which is a reduction of $(\gbigp_G,\theta)$.

Let us prove that
$(\gbigp^{v(h_1,h_2)}_{C_G(a_0)})^{h_i}_{C_G(a_0)\cap K}
:=\gbigp^{h_i}_K\cap
 \gbigp^{v(h_1,h_2)}_{C_G(a_0)}$
are $K\cap C_G(a_0)$-reductions.
Note that $v(h_1,h_2)$ is a section of
$\gbigp_{K}^{h_i}(\sqrt{-1}\gminik)$.
Let $P\in X$.
There exist $x_i\in\gbigp^{h_i}_K$
such that 
$v(h_1,h_2)(P)$ corresponds to
\[
(x_i,a_0)\in
\gbigp_K^{h_i}\times \gminit_{\real}
\subset
\gbigp_K^{h_i}\times \sqrt{-1}\gminik.
\]
Then,
we can observe that
$x_i\in (\gbigp^{v(h_1,h_2)}_{C_G(a_0)})^{h_i}_{K\cap C_G(a_0)}$,
and that
the fiber of
$(\gbigp^{v(h_1,h_2)}_{C_G(a_0)})^{h_i}_{K\cap C_G(a_0)}$
over $P$ is
the $K\cap C_G(a_0)$-orbit containing $x_i$.
Because
$\gbigp^{h_i}_{K}=
 (\gbigp^{v(h_1,h_2)}_{C_G(a_0)})^{h_i}_{K\cap C_G(a_0)}(K)$,
it is easy to check that
the $K\cap C_G(a_0)$-reduction
of $(\gbigp^{v(h_1,h_2)}_{C_G(a_0)},\theta)$
is harmonic.
\hfill\qed

\subsubsection{The case of compact Riemann surfaces with punctures}

Let $X$ be a compact Riemann surface.
Let $D\subset X$ be a finite subset.
Let $(\gbigp_G,\theta)$ be a $G$-Higgs bundle on $X\setminus D$.

\begin{prop}
\label{prop;25.6.18.31}
Let $h_i$ $(i=1,2)$ be harmonic metrics of
$(\gbigp_G,\theta)$
such that
\[
  \log \Tr\bigl(\chi_V(s(h_1,h_2))\bigr)
       =O\bigl(\log(-\log|z_P|)\bigr)
\]
for a faithful $G$-representation $\vecV$
around any $P\in D$.
Then, the following holds.
\begin{itemize}
 \item
       $h_1$ and $h_2$ are mutually bounded.
 \item $v(h_1,h_2)$ induces a constant map
       $[v(h_1,h_2)]:
       X\to \gminit/W(\gminit)$,
       whose image is
       contained in
       $\gminit_{\real}/W(\gminit)$.
 \item Let $a_0\in\gminit_{\real}$ be a representative
       of $[v(h_1,h_2)]$.
       Then,
       we obtain the reduction
       $(\gbigp^{v(h_1,h_2)}_{C_G(a_0)},\theta)$
       of $(\gbigp_G,\theta)$.
       Moreover,
       $h_i$ are induced by
       harmonic metrics of
       $(\gbigp^{v(h_1,h_2)}_{C_G(a_0)},\theta)$.
 \item $s(h_1,h_2)$ is a holomorphic section of
       $\gbigp_T(T,\Ad)$ such that
       $\Ad(s(h_1,h_2))\theta=\theta$.
\end{itemize}
\end{prop}
\pf
By using the arguments in the proof of Proposition \ref{prop;25.6.18.1},
we obtain the first claim.
Let $h_V$ be a $K$-invariant Hermitian metric of $V$.
We consider the Higgs bundle
$(E,\theta_E)=(\gbigp_G(\vecV),\chi_{\vecV}(\theta))$
with the induced harmonic metrics $\htilde_i=(h_i,h_V)$.
It is standard to obtain
$[\theta_E,s(\htilde_1,\htilde_2)]=0$
and
$\delbar_Es(\htilde_1,\htilde_2)$.
(See the proof of \cite[Proposition 2.5]{mochi4}.)
Then, we obtain the second and third claims
from Lemma \ref{lem;25.6.18.3}.
\hfill\qed

\subsection{Dirichlet problems}

\subsubsection{Smooth case}

Let $X$ be any Riemann surface.
Let $G$ be a reductive complex algebraic group
with a maximal compact subgroup $K$.
Let $(\gbigp_G,\theta)$ be a holomorphic principal $G$-Higgs bundle on $X$.
Let $Y$ be a relatively compact connected open subset of $X$
with smooth non-empty boundary $\del Y$.
Let $h_{\del Y}$ be a $K$-reduction of $\gbigp_{G|\del Y}$.
The following is a variant of Donaldson's theorem
\cite{Donaldson-boundary-value}.

\begin{prop}
\label{prop;25.9.24.12}
There exists a unique harmonic $K$-reduction $h$ of $(\gbigp_G,\theta)_{|Y}$
such that $h_{|\del Y}=h_{\del Y}$.
\end{prop}
\pf
We have already studied the uniqueness 
in Proposition \ref{prop;25.5.6.1}.
Let us study the existence.

Let $\vecV$ be a faithful $G$-representation
with a $K$-invariant Hermitian metric $h_V$.
We obtain the Higgs bundle
$(E,\theta_E)=(\gbigp_G(\vecV),\chi_{\vecV}(\theta))$.
We have
$T^{\veca,\vecb}E=\gbigp_{G}(T^{\veca,\vecb}\vecV)$
which is equipped with the induced Higgs field
$t^{\veca,\vecb}\theta_E$.
The induced Hermitian metric
$T^{\veca,\vecb}h_{E}$ is a harmonic metric of
$(T^{\veca,\vecb}E,t^{\veca,\vecb}\theta_E)$.
We have the decomposition (\ref{eq;25.5.5.21}) of
$T^{\veca,\vecb}(E)$
from the decomposition (\ref{eq;25.5.5.20}) of $T^{\veca,\vecb}V$.
It induces a decomposition of the Higgs bundle
\begin{equation}
\label{eq;25.5.4.20}
 (T^{\veca,\vecb}(E),t^{\veca,\vecb}\theta_E)
 =\bigoplus_{\vecW\in\Irr(G)}
 \bigl(
 T^{\veca,\vecb}(E)_{\vecW},
 t^{\veca,\vecb}(\theta_{E})_{\vecW}
 \bigr).
\end{equation}

Let $h_{E,\del Y}$ denote the Hermitian metric
of $E_{|\del Y}$ induced by $h_{\del Y}$ and $h_V$.
Because $h_{E,\del Y}$ comes from
a Hermitian metric of $\gbigp_{G|\del Y}$,
the restriction of the decomposition (\ref{eq;25.5.4.20})
to $\del Y$
is orthogonal with respect to $T^{\veca,\vecb}(h_{E,\del Y})$.
Moreover, the restriction of
$T^{\veca,\vecb}(h_{(E,\del Y)})$
to 
$T^{\veca,\vecb}(E)_{(\cnum,\II)|\del Y}$
is the constant Hermitian metric
$T^{\veca,\vecb}(h_V)_{(\cnum,\II)|\del Y}$.

According to the theorem of
Donaldson \cite{Donaldson-boundary-value}
(see also \cite{Li-Mochizuki2, Li-Mochizuki3}),
there exists a unique harmonic metric $h_E$ of $(E,\theta_E)_{|Y}$
such that $h_{E|\del Y}=h_{E,\del Y}$.
We have
$T^{\veca,\vecb}(h_E)_{|\del Y}=T^{\veca,\vecb}(h_{E,\del Y})$.
By the uniqueness of solutions of Dirichlet problem,
the decomposition (\ref{eq;25.5.4.20}) is orthogonal with respect to
$T^{\veca,\vecb}(h_E)$.
Moreover, the restriction of
$T^{\veca,\vecb}(h_E)$ to 
$T^{\veca,\vecb}(E)_{(\cnum,\II)}$
is the constant Hermitian metric
$T^{\veca,\vecb}(h)_{(\cnum,\II)}$.
By Proposition \ref{prop;25.5.5.22},
there exists a $K$-reduction $h$ of $\gbigp_G$
such that $h_V$ and $h$ induces $h_E$.
It is easy to check that
$h$ is a harmonic $K$-reduction of
$(\gbigp_G,\theta)_{|Y}$
such that $h_{|\del Y}=h_{\del Y}$.
\hfill\qed

\subsubsection{Removable singularity}

Let $D\subset X$ be a discrete subset.
Let $h_0$ be any $K$-reduction of $\gbigp_G$ on $X$.
\begin{prop}
\label{prop;25.5.10.2}
Let $h$ be a harmonic $K$-reduction of
$(\gbigp_G,\theta)_{|X\setminus D}$
such that the following holds
for a faithful representation $\vecV$.
\begin{itemize}
 \item $\log\Tr \chi_V(s(h_0,h))=O\bigl(
       \log(-\log|z_P|)
       \bigr)$ around any $P\in D$.
       We remark that this condition is independent of
       the choice of $h_0$.
 \end{itemize}
Then, $h$ induces a harmonic $K$-reduction of 
$(\gbigp_G,\theta)$ on $X$.
In other words,
there exists a harmonic $K$-reduction $\htilde$
of $(\gbigp_G,\theta)$ on $X$
such that $\htilde_{|X\setminus D}=h$.
\end{prop}
\pf
There exists a harmonic $K$-reduction $h_P$
of $(\gbigp_G,\theta)_{|X_P}$
such that
$h_{P|\del X_P}=h_{|\del X_P}$.
We obtain $h_{P|X_P\setminus\{P\}}=h_{|X_P\setminus\{P\}}$
by Proposition \ref{prop;25.5.6.1}.
It implies the claim of the proposition
by Proposition \ref{prop;25.5.6.1}.
\hfill\qed

\subsubsection{Compatibility with an automorphism of finite order}

Let $\sigma$ be a complex algebraic automorphism of $G$
of order $m$.
We use the notation in \S\ref{subsection;25.6.18.10}.
We consider the case where
$(\sigma,\omega)$ is not necessarily split.

Let $H$ be closed complex algebraic subgroup of $G$
such that $G^{\sigma}_0\subset H\subset G^{\sigma}$.
Let $\gbigp_{H}$
be a holomorphic principal $H$-bundle.
Let $\theta$ be a holomorphic section of
$\gbigp_{H}(\gminig_1)\otimes\Omega^1_X$.

We obtain the principal $G$-bundle
$\gbigp_{G}=\gbigp_{H}(G)$
which is equipped with the automorphism $\sigma$.
We have the induced holomorphic section
$\theta$ of $\gbigp_G(\gminig)\otimes\Omega^1_X$.
It satisfies $\sigma(\theta)=\omega\theta$.

\begin{thm}
\label{thm;25.5.13.40}
For any $K$-reduction $h_{\del Y}$
of $\gbigp_{G|\del Y}$
compatible with $\sigma$,
there exists a unique harmonic $K$-reduction $h$
of $(\gbigp_G,\theta)$ compatible with $\sigma$
such that 
$h_{|\del Y}=h_{\del Y}$.
\end{thm}
\pf
There exists a harmonic $K$-reduction $h$
of $(\gbigp_G,\theta)_{|Y}$
such that $h_{|\del Y}=h_{\del Y}$.
By the automorphism $\sigma$,
we obtain the induced harmonic $K$-reduction
$\sigma(h)$ of $(\gbigp_G,\omega\theta)$,
and $\sigma(h)$ is also a harmonic $K$-reduction of 
$(\gbigp_G,\theta)$.
Because $h_{\del Y}$ is $\sigma$-invariant,
$\sigma(h)=h$ holds on $\del Y$.
By the uniqueness,
we obtain $\sigma(h)=h$ on $Y$.
By Lemma \ref{lem;25.6.16.3},
$h$ is a harmonic $K$-reduction of $(\gbigp_G,\theta)$
compatible with $\sigma$.
\hfill\qed

\section{Local estimates and global existence in the split case}

\subsection{Local estimates}

\subsubsection{Regular semisimple Higgs fields}

Let $G$ be a semisimple complex Lie group
with a complex algebraic automorphism $\sigma$
of order $m>0$.
Let $\gminig$ denote the associated Lie algebra.
Let $\omega$ be a primitive $m$-th root of $1$.
We assume that $(\sigma,\omega)$ is split.

Let $H$ be a complex algebraic closed subgroup of $G$
such that
$G_0^{\sigma}\subset H\subset G^{\sigma}_Z$.
Let $\gbigp_{H}$ be
a holomorphic principal $H$-bundle
on a Riemann surface $X$.
Let $\theta$ be a holomorphic section of
$\gbigp_{H}(\gminig_1)\otimes\Omega^1$.
\begin{df}
\label{df;25.6.20.1}
We say that $\theta$ is regular semisimple
if for any local holomorphic coordinate $(U,z)$,
$\theta_{|U}$ is expressed as $f_U\,dz$,
where $f_U$ is a holomorphic section of
$\gbigp_{H}(\gminig_1^{\rs})_{|U}$.
We say that $\theta$ is generically regular semisimple
if there exists a discrete subset $D\subset X$
such that $\theta_{|X\setminus D}$ is regular semisimple.
\hfill\qed
\end{df}
The following lemma is clear.
\begin{lem}
$\theta$ is generically regular semisimple
if and only if
there exists $P\in X$ such that
$\theta$ is regular semisimple at $P$.
\hfill\qed
\end{lem}

\subsubsection{Canonical decoupled harmonic metrics}

Suppose that $\theta$ is regular semisimple.
For a holomorphic local coordinate $(U,z)$,
we have the unique Hermitian metric $h^{\can}_U$
of $\gbigp_{G|U}$
compatible with $\sigma$
such that
$[f_U,\rho_{h^{\can}_U}(f_U)]=0$
as in Lemma \ref{lem;25.6.20.2}.
By gluing the Hermitian metrics
$h^{\can}_U$ of $\gbigp_{G|U}$ 
compatible with $\sigma$ for varied $U$,
we obtain
a global Hermitian metric
$h^{\can}$ of $\gbigp_G$.

\begin{lem}
\label{lem;25.5.13.30}
$h^{\can}$ is a unique decoupled harmonic metric of
$(\gbigp_G,\theta_G)$ compatible with $\sigma$.
 \end{lem}
\pf
We have $[\theta,\rho_{h^{\can}}(\theta)]=0$.
Let $\gminia\subset\gminig_1$ be a Cartan subspace.
Let $K(\gminia)$
be the maximal compact subgroup
as in \S\ref{subsection;25.6.18.11}.
By Proposition \ref{prop;25.5.13.31},
after a holomorphic gauge transform,
there exists a holomorphic isomorphism
$\gbigp_{H|U}\simeq
U\times H$
such that
$f_U$ is a section of
$\gminia\otimes\nbigo_U$.
In this case,
$h^{\can}_U$ is constantly $K(\gminia)$.
Hence, the associated Chern connection is flat.
\hfill\qed

\subsubsection{Comparison with canonical metrics in the regular semisimple case}

For $R>0$,
we set $U(R)=\bigl\{z\in\cnum\,\big|\,|z|<R\bigr\}$.
We consider the holomorphic principal
$H$-bundle $\gbigp_H=U(R)\times H$.
Let $\theta=f\,dz$ be a holomorphic section of
$\gbigp_{H}(\gminig_1)\otimes\Omega^1$.
Let $\gbigp_G=\gbigp_H(G)$.
Let $\Harm(\gbigp_G,\theta,\sigma)$
denote the set of harmonic metrics of $(\gbigp_G,\theta)$
compatible with $\sigma$.
Any $K$-reduction $h$ of $\gbigp_G$
and the Killing form of $\gminig$ induce
a Hermitian metric of $\gbigp_G(\gminig)$,
which is also denoted by $h$.

Suppose that $\theta$ is regular semisimple.
We have
$h^{\can}\in\Harm(\gbigp_G,\theta,\sigma)$.
Let $\nbigk$ be a compact subset of $\gminit^{\rs}/W(\gminit)$.
Let $\Psi_f:U(R)\to\gminit/W(\gminit)$ be the holomorphic map
induced by $f$ as in \S\ref{subsection;25.9.24.10}.
We suppose that the image of
$\Psi_f$ is contained in $\nbigk$.
There exists the natural $\cnum^{\ast}$-action on $\gminit^{\rs}/W(\gminit)$
induced by the $\cnum^{\ast}$-multiplication on $\gminit$.
Note that $h^{\can}$ is
the unique decoupled harmonic metric of
$(\gbigp_G,a\theta)$ compatible with $\sigma$
for any $a\in\cnum^{\ast}$.
We choose a $W(\gminit)$-invariant Hermitian metric $h_{\gminit}$
of $\gminit$.

\begin{thm}
\label{thm;25.6.20.20}
There exist $C_i>0$ $(i=1,2)$
depending only on
$(R,R_1,\nbigk,h_{\gminit})$
such that 
\begin{equation}
\label{eq;25.6.20.10}
\bigl|
 v(h^{\can},h)
 \bigr|_{h^{\can}}
 \leq
 C_1\exp(-C_2|a|)
\end{equation}
for $a\in\cnum^{\ast}$ with $|a|\geq 1$
and any $h\in \Harm(\gbigp_G,a\theta,\sigma)$  on $U(R_1)$.
\end{thm}
\pf
By Proposition \ref{prop;25.5.13.50},
there exit $C_i>0$ $(i=3,4)$ depending only on
$(R,R_1,\nbigk,h_{\gminit})$
such that the following holds on $U(R_1)$:
\[
 |af|_{h}
 \leq
 C_3|a|
+C_4.
\]
There exits $C_5>0$ depending only on
$(R,R_1,\nbigk,h_{\gminit})$
such that the following holds on $U(R_1)$:
\[
 |f|_{h}
 \leq
 C_5.
\]
By Proposition \ref{prop;25.5.13.21},
there exists $C_6>0$ depending only on
$(R,R_1,\nbigk,h_{\gminit})$
such that
\[
 |v(h^{\can},h)|_{h^{\can}}\leq C_6.
\]

Because $\ad(f)$ is semisimple,
by \cite[Corollary 2.6]{Decouple},
there exist $C_i>0$ $(i=7,8)$
depending on $(R,R_1,\nbigk,h_{\gminit})$
such that 
\[
 \bigl|
 [f,f^{\dagger}_h]
 \bigr|_h
 \leq
 C_7\exp(-C_8|a|).
\]
We obtain the desired estimate (\ref{eq;25.6.20.10})
from Proposition \ref{prop;25.6.18.50}.
\hfill\qed

\subsubsection{An estimate in the generically regular semisimple case}

Let $X$ be a Riemann surface.
Let $\gbigp_H$ be a holomorphic principal $H$-bundle on $X$.
Let $\theta$ be a holomorphic section of
$\gbigp_H(\gminig_1)\otimes \Omega_X^1$.
Suppose that $\theta$ is generically regular semisimple.
Let $h_0$ be any Hermitian metric of $\gbigp_G$
compatible with $\sigma$.

\begin{prop}
\label{prop;25.6.20.30}
For any relatively compact open subset $Y\subset X$,
there exists $C>0$ such that 
\begin{equation}
\label{eq;25.9.24.11}
 \bigl|
 v(h_0,h)
 \bigr|_{h_0}
 \leq
 C
\end{equation}
holds on $Y$ for any $h\in\Harm(\gbigp_G,\theta,\sigma)$.
\end{prop}
\pf
Let $D\subset X$ be a discrete subset
such that $\theta_{|X\setminus D}$ is regular semisimple.
We may assume that
$\del Y$ is smooth and $D\cap\del Y=\emptyset$.

We consider the induced Higgs bundle
$(\gbigp_G(\gminig),\ad\theta)$.
The Hermitian metric $h_0$ of $\gbigp_G$
induces a Hermitian metric $\htilde_0$
of the vector bundle $\gbigp_G(\gminig)$.
Any $h\in\Herm(\gbigp_G,\theta,\sigma)$
induces a harmonic metric $\htilde$ of
$(\gbigp_G(\gminig),\ad\theta)$.
By \cite[Lemma 3.1]{s1},
we obtain
\[
 \sqrt{-1}\Lambda_{g_X}\delbar\del
 \log\Tr(\Ad s(h_0,h))
 \leq
 |\ad F(h_0)|_{\htilde_0,g_X}
\]
on $X$,
where $g_X$ denotes a conformal metric of $X$.

For $P\in Y\cap D$,
let $Y_P$ be a relatively compact neighbourhood of $P$ in $Y$.
Let $Y_P'$ be a relatively compact neighbourhood of $P$ in $Y_P$.
By Theorem \ref{thm;25.6.20.20},
there exists $C_1>0$,
which is independent of $h$,
such that
$\Tr(s(\htilde_0,\htilde))\leq C_1$
on $Y\setminus \bigcup_{P\in D}Y_P'$.
There exist functions $\alpha_P$ on $Y_P$
such that
$\sqrt{-1}\Lambda\delbar\del \alpha_P=
|\ad F(h_0)|_{\htilde_0}$.
Then, by using the maximum principle
for $\log\Tr (\Ad s(h_0,h))-\alpha_P$,
we can prove that 
there exists $C_2>0$,
which is independent of $h$,
such that
$\Tr(s(\htilde_0,\htilde))\leq C_2$
on $Y$.
(See \cite[\S2.3.4]{Li-Mochizuki3}.)
We obtain the estimate (\ref{eq;25.9.24.11}).
\hfill\qed

\subsection{Existence of harmonic metrics
in the generically regular semisimple case}

Let $H$ be a complex algebraic closed subgroup of $G$
such that $G_0^{\sigma}\subset H\subset G^{\sigma}$.
Let $\gbigp_{H}$ be
a holomorphic principal $H$-bundle
on a non-compact Riemann surface $X$.
Let $\theta$ be a section of
$\gbigp_{H}(\gminig_1)\otimes\Omega^1$.
Assume the following conditions.
\begin{itemize}
 \item $(\sigma,\omega)$ is split.
 \item $\theta$ is generically regular semisimple,
       i.e.,
       there exists a discrete subset $D\subset X$
       such that
       $\theta_{|X\setminus D}$ is regular semisimple.
       It is equivalent to the existence of a point $P\in X$
       such that $\theta$ is regular semisimple at $P$.
\end{itemize}

The following theorem is a generalization of 
\cite[Theorem 2.34]{Li-Mochizuki3}.
\begin{thm}
\label{thm;25.5.13.32}
There exists a harmonic metric of $(\gbigp_G,\theta)$
compatible with $\sigma$.
\end{thm}
\pf
First, we consider the non-compact case.
We take a smooth exhaustion family $Y_i\subset X$
$(i=1,2,\ldots)$,
i.e.,
$Y_i$ are connected relatively compact open subsets of $X$
with smooth boundary $\del Y_i$
such that $Y_{i+1}$ are relatively compact in $Y_i$
and that $\bigcup Y_i=X$.
We also assume that $D\cap \del Y_i=\emptyset$.
Let $h_i$ be harmonic metrics of
$(\gbigp_G,\theta)_{|Y_i}$
such that $h_{i}$ are compatible with $\sigma$
as in Theorem \ref{thm;25.5.13.40}.
Let $\htilde_i$ denote the harmonic metric of
the harmonic metric of
$(\gbigp_G(\gminig),\ad\theta)$ on $Y_i$
induced by $h_i$ and the Killing form of $\gminig$.
Due to the estimates in Proposition \ref{prop;25.6.20.30},
there exists a convergent subsequence $\htilde_{j'}$
by \cite[Proposition 2.6]{Li-Mochizuki2}.
The limit is induced by
a harmonic metric $h_{\infty}$ of
$(\gbigp_G,\theta)$ compatible with $\sigma$.

Let us consider the case where $X$ is compact.
Let $P\in X$.
We set $U=X\setminus\{P\}$.
Let $X_P$ be a neighbourhood of $P$.
Let $X_P'$ be a neighbourhood of the closure of $X_P$.
We may assume that
$\theta$ is regular semisimple on $X_P'\setminus\{P\}$.
Let $h_0$ be a harmonic metric of 
$(\gbigp_G,\theta)_{|X'_P}$.
Note that the existence follows from
Proposition \ref{prop;25.9.24.12}.

Let $X_i\subset X$ $(i=1,2,\ldots)$
be an exhaustive family with smooth boundary
such that $\del X_i\subset X_P\setminus\{P\}$.
Let $h_i$ be a harmonic metric of
$(\gbigp_G,\theta)_{|X_i}$
such that
$h_{i|\del X_i}=h_{0|\del X_i}$.

We obtain the functions
$\Tr\Bigl(
\Ad\bigl(s(h_{0|X_i},h_i)\bigr)\Bigr)$
on $X_i\cap X_P'$.
Note that $\det \Ad\bigl(s_{h_{0|X_i},h_i}\bigr)=1$.
By Theorem \ref{thm;25.6.20.20},
we may assume that there exists $C>0$
such that
$\Tr\Bigl(
\Ad\bigl(
s(h_{0|X_i},h_i)\bigr)
\Bigr)\leq C$
for any $i$ around $\del X_P$.
We also have
$\Tr\Bigl(
\Ad\bigl(
s(h_{0|X_i},h_i)\bigr)
\Bigr)=\dim\gminig$ on $\del X_i$.
Because
$\Tr\Bigl(
\Ad\bigl(
s(h_{0|X_i},h_i)\bigr)\Bigr)
$ are subharmonic,
we obtain that
\[
\Tr\Bigl(\Ad\bigl(
s(h_{0|X_i},h_i)\bigr)\Bigr)
\leq \dim\gminig+C
\]
on $X_i\cap X_P$ for any $i$.
We obtain that
$h_{\infty}$ and $h_0$ are mutually bounded.
By Proposition \ref{prop;25.5.10.2},
we obtain that $h_{\infty}$ induces
a harmonic metric $\htilde_{\infty}$ of $(\gbigp_G,\theta)$.
Because $h_{\infty}$ is compatible with $\sigma$,
$\htilde_{\infty}$ is also compatible with $\sigma$.
\hfill\qed

\section{Classification in the generically cyclic case}

\subsection{Cyclic Higgs fields and generically cyclic Higgs fields}

Suppose that $\gminig$ is a simple complex Lie algebra.
Let $G$ be the associated Lie group of adjoint type,
i.e., $G\subset\GL(\gminig)$.
Let $\gminit$ be a Cartan subalgebra of $\gminig$.
Let $T\subset G$ be the corresponding maximal complex torus.
We use the notation in \S\ref{subsection;25.6.18.20}.
Let $\psi_i\in\seisuu_{>0}$
be determined by $\psi=\sum \psi_i\alpha_i$
for the highest root $\psi$.
We define the one dimensional $\cnum$-vector space
\[
 \Upsilon(\gminig,\gminit):=
 \bigotimes_{i=1}^{\ell}
 \gminig^{\otimes\psi_i}_{\alpha_i}
 \otimes\gminig_{-\psi}.
\]
Note that the natural $T$-action on $\Upsilon(\gminig,\gminit)$
is trivial.

Let $X$ be a Riemann surface.
Let $\gbigp_T$ be a holomorphic principal $T$-bundle on $X$.
We set $\gbigp_G=\gbigp_T(G)$.
We obtain the decompositions
\begin{equation}
\label{eq;25.6.18.22}
\gbigp_G(\gminig)
=\bigoplus_{\phi\in\Delta}
\gbigp_T(\gminig_{\phi}),
\quad\quad
\gbigp_T(\gminig_1)
 =\bigoplus_{i=1}^{\ell}
 \gbigp_T(\gminig_{\alpha_i})
 \otimes
 \gbigp_T(\gminig_{-\psi}).
\end{equation}
For any line bundle $L$,
there exists the natural isomorphism
\[
 \bigotimes_{i=1}^{\ell}
 \bigl(
 \gbigp_T(\gminig_{\alpha_i})
 \otimes_{\nbigo_X} L
 \bigr)^{\otimes\psi_i}
 \otimes_{\nbigo_X}
 \bigl(
 \gbigp_T(\gminig_{-\psi})
 \otimes_{\nbigo_X} L
 \bigr)
 \simeq
 \Upsilon(\gminig,\gminit)
 \otimes_{\cnum}
 L^{\otimes(\tth+1)}.
\]

Let $\theta$ be a holomorphic section of
$\gbigp_T(\gminig_1)\otimes K_X$.
According to (\ref{eq;25.6.18.22}),
we obtain the decomposition
$\theta=\sum_{i=1}^{\ell} \theta_{\alpha_i}
+\theta_{-\psi}$.
We set
\[
 \gminio(\theta)
 :=\prod_{i=1}^{\ell}
 \theta_{\alpha_i}^{\otimes\psi_i}
 \cdot \theta_{-\psi}
 \in H^0\bigl(X,\Upsilon(\gminig,\gminit)\otimes K_X^{\otimes(\tth+1)}\bigr).
\]

\begin{df}
We say that $\theta$ is cyclic
if $\gminio(\theta)$ is nowhere vanishing.
We say that $\theta$ is generically cyclic
if $\gminio(\theta)$ is not constantly $0$.
\hfill\qed
\end{df}

\subsubsection{Canonical decoupled harmonic metrics}

As a special case of the construction in \S\ref{subsection;25.6.18.30},
we obtain the following lemma.
\begin{lem}
If $\theta$ is cyclic,
we obtain the canonical decoupled harmonic metric
$h^{\can}$ of
$(\gbigp_G,\theta)$
compatible with $\Ad(\ttw)$.
By the compatibility with $\Ad(\ttw)$,
$h^{\can}$ is induced by 
a $T_c$-reduction $\gbigp^{h^{\can}}_{T_c}\subset\gbigp_T$.
\hfill\qed
\end{lem}

\subsubsection{General existence result
in the generically regular semisimple case}

As a special case of Theorem \ref{thm;25.5.13.40},
we obtain the following.
\begin{prop}
\label{prop;25.6.19.1}
If $\theta$ is generically cyclic,
there exists a harmonic metric $h$ of
$(\gbigp_G,\theta)$
compatible with the reduction
$\gbigp_T\subset\gbigp_G$.
It is induced by a reduction
$\gbigp^h_{T_c}\subset \gbigp^h_T$.
\hfill\qed
\end{prop}

\subsubsection{Preliminary for the classification}

Let $X$ be a compact Riemann surface.
Let $D\subset X$ be a finite subset.
Let $\gbigp_T$ be a principal $T$-bundle on $X\setminus D$.
Let $\theta$ be a section of
$\gbigp_T(\gminig_1)\otimes\Omega^1$.
Assume that $\theta$ is generically cyclic on $X\setminus D$,
i.e.,
there exists a discrete subset $Z\subset X\setminus D$
such that
$\theta_{|X\setminus (D\cup Z)}$ is cyclic.
In this case, we have the following refinement of
Proposition \ref{prop;25.6.18.31}.
\begin{prop}
\label{prop;25.5.9.3}
Let $h_i$ $(i=1,2)$ be harmonic metrics of
$(\gbigp_G,\theta)$
compatible with the reduction $\gbigp_T\subset\gbigp_G$
satisfying the following conditions.
\begin{itemize}
 \item $\log\Tr s(h_1,h_2)=O\bigl(\log(-\log|z_P|)\bigr)$
       around any $P\in D$.
\end{itemize}
Then, we obtain $h_1=h_2$.
\end{prop}
\pf
According to Proposition \ref{prop;25.6.18.31},
we obtain that $s(h_1,h_2)$
is a holomorphic section of
$\gbigp_T(T,\Ad)$ such that
$\Ad(s(h_1,h_2))\theta=\theta$.
Note that $s(h_1,h_2)$
is a section of
$\gbigp^{h_1}_{T_c}(T_{\real})$.
We obtain that $s(h_1,h_2)=\id$
by the following lemma.
\hfill\qed

\begin{lem}
Let $s\in T_{\real}$.
If $\Ad(s)u=u$ for a cyclic element $u\in\gminig_1$,
we obtain $s=1$.
\end{lem}
\pf
Let $s=\exp(v)$ for $v\in\gminit_{\real}$.
We have the decomposition
$u=\sum_{i=1}^{\ell} u_{\alpha_i}+u_{-\psi}$,
where $u_{\phi}\in\gminig_{\phi}$
are non-zero.
We obtain
\[
\Ad(\exp(v))\Bigl(
\sum_{i=1}^{\ell} u_{\alpha_i}
+u_{-\psi}
\Bigr)
=\sum_{i=1}
e^{\alpha_i(v)}u_{\alpha_i}
+e^{-\psi(v)}u_{-\psi}
=\sum_{i=1}^{\ell}u_{\alpha_i}
+u_{-\psi}.
\]
It implies
$\alpha_i(v)=0$
and $-\psi(v)=0$.
We obtain that $v=0$.
\hfill\qed

\subsection{Model harmonic $G$-Higgs bundles}
\subsubsection{Normalization of compact real forms compatible with $\gminit$}
\label{subsection;25.5.3.10}

\begin{prop}
\label{prop;25.5.3.3}
For any $K\in\Herm(G,\gminit)$,
there exist root vectors $e_{\phi}$ in $\gminig_{\phi}$ $(\phi\in\Delta)$
satisfying the following conditions.
\begin{description}
 \item[(a)] $[e_{\phi},e_{-\phi}]=\phi$,
       or equivalently
       $B(e_{\phi},e_{-\phi})=1$.
 \item[(b)] 
	    Let $N_{\phi,\phi'}$
	    $(\phi+\phi'\in\Delta)$
	    be the numbers determined by
	    $[e_{\phi},e_{\phi'}]=N_{\phi,\phi'}e_{\phi+\phi'}$.
	    Then,
	    $N_{\phi,\phi'}$ are real,
	    and $N_{\phi,\phi'}=-N_{-\phi,-\phi'}$ holds.
 \item[(c)] Let $\gminik$ denote the Lie algebra of $K$.
	    Then,
\begin{equation}
\label{eq;25.5.3.2}
 \gminik
 =\bigoplus_{\phi\in\Delta}
 \Bigl(
 \real (e_{\phi}-e_{-\phi})
 \oplus
 \real \sqrt{-1}(e_{\phi}+e_{-\phi})
 \Bigr)
 \oplus
 \sqrt{-1}\gminit_{\real}.
\end{equation}
\end{description}
\end{prop}
\pf
According to \cite[Theorem 6.6]{Knapp-Book},
there exist root vectors
$e^{(0)}_{\phi}\in\gminig_{\phi}$ $(\phi\in\Delta)$
satisfying the conditions (a) and (b).
Then, as explained in \cite[\S6.1]{Knapp-Book},
\[
\gminik^{(0)}
 =\bigoplus_{\phi\in\Delta}
 \Bigl(
 \real (e^{(0)}_{\phi}-e^{(0)}_{-\phi})
 \oplus
 \real \sqrt{-1}(e^{(0)}_{\phi}+e^{(0)}_{-\phi})
 \Bigr)
 \oplus
 \sqrt{-1}\gminit_{\real}.
\]
is a compact real form of $\gminig$.
Let $K^{(0)}$ denote the corresponding maximal compact subgroup.
Clearly, $K^{(0)}\in \Herm(G,\gminit)$.

Let $K\in \Herm(G,\sigma)$.
Let $\gminik$ denote the Lie algebra of $K$.
By Corollary \ref{cor;25.5.3.1},
there exists $v\in \gminit_{\real}$
such that $\gminik=\Ad(\exp(v))(\gminik^{(0)})$.
We set
$e_{\phi}=\exp(v(\phi))e^{(0)}_{\phi}$.
The tuple $e_{\phi}$ $(\phi\in\Delta)$
satisfies the conditions (a) and (b).
Moreover, (\ref{eq;25.5.3.2}) holds.
\hfill\qed

\begin{rem}
If $\phi+\phi'\neq 0$ and $\phi+\phi'\not\in\Delta$,
we have $[\gminig_{\phi},\gminig_{\phi'}]=0$,
and hence
$[e_{\phi},e_{\phi'}]=0$.
\hfill\qed 
 \end{rem}

The following lemma is remarked in \cite[\S6]{Kostant-TDS}.
\begin{lem}
For $K\in\Herm(G,\gminit)$,
let $e_{\phi}$ $(\phi\in\Delta)$ be a tuple of root vectors
as in Proposition {\rm\ref{prop;25.5.3.3}}.
Then, $\rho_{\gminik}(e_{\phi})=-e_{-\phi}$ holds.
\end{lem}
\pf
We have
$\rho_{\gminik}(e_{\phi}-e_{-\phi})=e_{\phi}-e_{-\phi}$
and 
$\rho_{\gminik}\bigl(
 \sqrt{-1}(e_{\phi}+e_{-\phi})
 \bigr)
=\sqrt{-1}(e_{\phi}+e_{-\phi})$.
Then, we obtain the claim of the lemma.
\hfill\qed

\subsubsection{Some computations}

Let $e_{\phi}$ $(\phi\in\Delta)$
be a tuple of root vectors as in
Proposition {\rm\ref{prop;25.5.3.3}}
for $K\in\Herm(G,\gminit)$.
As remarked in \cite[\S6]{Kostant-TDS},
for $\phi_1,\phi_2\in \{\alpha_1,\ldots,\alpha_{\ell}\}\sqcup\{-\psi\}$,
we have
\[
[e_{\phi_1},e_{-\phi_2}]
=\left\{
\begin{array}{ll}
 \phi_1& (\phi_1=\phi_2) \\
 0 & (\phi_1\neq\phi_2).
\end{array}
\right.
\]
As a result,
we obtain the following lemma,
which is contained in \cite[\S6]{Kostant-TDS}.
\begin{lem}
For any
$u=\sum c_i \psi_i^{1/2}e_{\alpha_i}+c_{-\psi}e_{-\psi}
\in\gminig_1$,
we have
\begin{equation}
\label{eq;25.4.27.51}
 \bigl[u,\rho(u)\bigr]=
-\sum_{i=1}^{\ell}
 |c_i|^2\psi_i\alpha_i
 +|c_{-\psi}|^2\psi
=-\sum_{i=1}^{\ell}
 \bigl(|c_i|^2-|c_{-\psi}|^2\bigr)
 \psi_i\alpha_i.
\end{equation} 
As a result,
we have $[u,\rho(u)]=0$
if and only if
$|c_1|=\cdots=|c_{\ell}|=|c_{-\psi}|$. 
\hfill\qed 
\end{lem}

\begin{cor}
Let $g=\exp\bigl(\sum_{i=1}^{\ell} \beta_i\epsilon_i\bigr)\in T$.
We set
$u=\Ad(g)\Bigl(\sum_{i=1}^{\ell}\psi_i^{1/2}e_{\alpha_i}+e_{-\psi}\Bigr)$.
Then, $\bigl[u,\rho(u)\bigr]=0$ holds
if and only if
$\beta_i$ are purely imaginary.
\end{cor}
\pf
We have
$u=\sum_{i=1}^{\ell}e^{\beta_i}\psi_{i}^{1/2}e_{\alpha_i}
 +e^{-\sum \psi_i\beta_i}e_{-\psi}$.
By the previous lemma,
$[u,\rho(u)]=0$ holds
if and only if
\[
 \Re(\beta_j)
=-\sum_{i=1}^{\ell}\psi_i\Re (\beta_i)
\]
for any $j=1,\ldots,\ell$.
It implies that $\Re(\beta_j)=0$ for $j=1,\ldots,\ell$.
\hfill\qed

\subsubsection{Some $\minisl_2$-triples}

We identify $\gminit$ and $\gminit^{\lor}$
by using the restriction of the Killing form.
In particular, we regard $\Delta$ as a subset of $\gminit_{\real}$.
We set
$\Pi^Q=\{\alpha_1,\ldots,\alpha_{\ell},-\psi\}$.
Let us consider a non-empty subset
$S\subsetneq \Pi^Q$.
Let $\gminit_{\real,S}\subset\gminit_{\real}$
denote the $\real$-subspace generated by $\phi\in S$.
Because $S\neq\Pi^Q$,
the tuple $\phi\in S$ is a base of $\gminit_{\real,S}$.

Let $K\in\Herm(G,\gminit)$.
Let $\rho$ denote the Cartan involution associated with $K$.
Let $e_{\phi}$ $(\phi\in\Delta)$ be
a tuple of root vectors as in Proposition \ref{prop;25.5.3.3}.

\begin{prop}\mbox{{}}
\label{prop;25.5.8.1}
\begin{itemize}
 \item There exists a unique element
       $u_S\in \gminit_{\real,S}$
       such that
       $\phi(u_S)=1$ for any $\phi\in S$.
 \item There exists an element
       $v_S=\sum_{\phi\in S} v_{\phi}
       \in \bigoplus_{\phi\in S}\gminig_{\phi}$
       such that
       $[v_S,-\rho(v_S)]=2u_S$.
 \item  An element $v_S'=\sum_{\phi\in S}v'_{\phi}\in
	\bigoplus_{\phi\in S}\gminig_{\phi}$
	satisfies
	$[v'_S,-\rho(v'_S)]=2u_S$
	if and only if
	there exist $\gamma_{\phi}\in \cnum^{\ast}$
	$(\phi\in S)$ with $|\gamma_{\phi}|=1$
	such that $v'_{\phi}=\gamma_{\phi}v_{\phi}$
	for any $\phi\in S$.
 \item The tuple $u_S$, $v_S$, $-\rho(v_S)$
       determines a Lie subalgebra of $\gminig$
       isomorphic to
       $\minisl_2$.
       An isomorphism is given by
\[
       u_S\leftrightarrow
       \frac{1}{2}
       \left(
       \begin{array}{cc}
	1 & 0 \\ 0 & -1
       \end{array}
       \right),
       \quad
       v_S\leftrightarrow
       \left(
       \begin{array}{cc}
	0 & 1 \\ 0 & 0
       \end{array}
       \right),
       \quad
       -\rho(v_S)\leftrightarrow
       \left(
       \begin{array}{cc}
	0 & 0 \\ 1 & 0
       \end{array}
       \right).
\]
 \end{itemize}
\end{prop}
\pf
Let $\langle\cdot,\cdot\rangle$
denote the restriction of the Killing form
to $\gminit_{\real}$.
Because $\langle\cdot,\cdot\rangle$ is positive definite,
its restriction to $\gminit_{\real,S}$ is also positive definite.
Hence, there exists a unique element
$u_S\in \gminit_{\real,S}$
such that $\langle u_S,\phi\rangle=1$ $(\phi\in S)$.
It implies $\phi(u_S)=1$
by the identification of $\gminit$ and $\gminit^{\lor}$.

Let $v_S=\sum_{\phi\in S} \beta_{\phi}e_{\phi}$.
Because
$-\rho(v_S)
 =\sum_{\phi\in S}\overline{\beta_{\phi}}e_{-\phi}$,
we obtain
\[
\bigl[
v_S,-\rho(v_S)
\bigr]
=\sum_{\phi\in S}
 |\beta_{\phi}|^2[e_{\phi},e_{-\phi}]
 =\sum_{\phi\in S}|\beta_{\phi}|^2\phi
 \in\gminit_{\real,S}.
\]
Because the tuple $(\phi\in S)$ is a base of
$\gminit_{\real,S}$,
the second and third claims are clear.
Because $[u_S,v_S]=v_S$ and
$[u_S,-\rho(v_S)]=-(-\rho(v_S))$,
the fourth claim is also clear by the commutation relations.
\hfill\qed

\subsubsection{Model harmonic $G$-Higgs bundles}

Let $K\in\Herm(G,\gminit)$.
\begin{itemize}
 \item Let $S\subsetneq \Pi^Q$
       be a non-empty subset.
 \item Let $\vecm=(m(\phi)\,|\,\phi\in S)\in\seisuu^{S}$.
 \item Let $\beta\in \gminit_{\real}$
       be an element such that
       $\beta(\phi)=m(\phi)$ $(\phi\in S)$.
\end{itemize}

We set $U:=\{|z|<1\}$ and $U^{\ast}:=U\setminus\{0\}$.
We consider the holomorphic principal $T$-bundle
$\gbigp_T=U^{\ast}\times T$ on $U^{\ast}$.
We set
$\gbigp_G=\gbigp_T(G)=U^{\ast}\times G$.
Let $u_S\in\gminit_{\real,S}$
and
$v_S=\sum_{\phi\in S}v_{\phi}
\in\bigoplus_{\phi\in S}\gminig_{\phi}$
be as in Proposition \ref{prop;25.5.8.1}.

Let $t>0$.
We consider the holomorphic map
$\vtilde_{S,\vecm}=\sum_{\phi\in S}
 z^{m(\phi)}v_{\phi}:U^{\ast}\to \gminig$.
We obtain the following Higgs field of $\gbigp_G$:
\[
 \theta_{S,\vecm,t}=
 t\vtilde_{S,\vecm}\frac{dz}{z}.
\]
We consider the following sections
of $U^{\ast}\times T_{\real}$ on $U^{\ast}$:
\[
 g_{\beta,S,t}=
 \exp\Bigl(
 -\beta\log|z|
 -u_S\log(-\log|z|^{2t}).
 \Bigr)
\]
We obtain the Hermitian metric
\[
 h_{\beta,S,t}=
 g_{\beta,S,t}^{-1}Kg_{\beta,S,t}
 \in
 \Gamma\bigl(
 \Herm(\gbigp_G,\Ad\ttw)
 \bigr).
\]
\begin{prop}
\label{prop;25.9.25.30}
The $G$-Higgs bundle $(\gbigp_G,\theta_{S,\vecm,t})$
with the Hermitian metric $h_{\beta,S,t}$
is a harmonic $G$-bundle.
Moreover, we have the following estimate 
with respect to the induced harmonic metric
$\htilde_{\beta,S,t}$ on $\gbigp_G(\gminig)$:
\begin{equation}
\label{eq;25.9.25.31}
 \log|e_{\phi}|_{\htilde_{\beta,S,t}}
+\log|z|^{\beta(\phi)}
=O\bigl(\log(-\log|z|)\bigr)
\quad(\phi\in\Delta).
\end{equation}
\end{prop}
\pf
Let $h_{\gminig}$ denote
the Hermitian metric of $\gminig$
induced by the Killing form $B$
and the Cartan involution associated with $K$,
i.e.,
$h_{\gminig}(x_1,x_2)=B(x_1,-\rho(x_2))$.
By using the adjoint $G$-action on $\gminig$,
we obtain the Higgs bundle
$(\gbigp_T(\gminig),\ad\theta_{S,\vecm,t})$
with the induced Hermitian metric $\htilde_{\beta,S,t}$.
The trivialization $\gbigp_G=U^{\ast}\times G$
induces
the isomorphism
$\Psi:U^{\ast}\times\gminig\simeq\gbigp_G(\gminig)$.
It is equipped with the Hermitian metric
$\Psi^{\ast}(\htilde_{\beta,S,t})$.
The corresponding Higgs field is
$t\ad(\vtilde_{S,\vecm})dz/z$.

Let us check that
the tuple
$(U^{\ast}\times\gminig,t\ad(\vtilde_{S,\vecm})dz/z,
\Psi^{\ast}(\htilde_{\beta,S,t}))$
is a harmonic bundle.
By Lemma \ref{lem;25.9.25.20},
the induced Hermitian metric $\Psi^{\ast}(\htilde_{\beta,S,t})$
of $U^{\ast}\times\gminig$
is described as follows:
\[
 \Psi^{\ast}(\htilde_{\beta,S,t})(x_1,x_2)
 =h_{\gminig}(\Ad(g_{\beta,S,t})x_1,\Ad(g_{\beta,S,t})x_2)
 =h_{\gminig}\bigl(
 \Ad(g_{\beta,S,t}^2)x_1,x_2
 \bigr).
\]
We obtain the estimates (\ref{eq;25.9.25.31}).

The curvature
$R(\Psi^{\ast}\htilde_{\beta,S,t})$
of the Chern connection
of $(U^{\ast}\times\gminig,\Psi^{\ast}(\htilde_{\beta,S,t}))$
equals
\[
 \delbar\bigl(
 (\Ad g_{\beta,S,t}^2)^{-1}\del(\Ad g_{\beta,S,t}^2)
 \bigr)
=-2\ad(u_S)\delbar\del\log(-\log|z|^{2t})
=-2\ad(u_S)
\frac{1}{(\log|z|^{2})^2}\frac{dz\,d\zbar}{|z|^2}.
\]
Because
$\bigl(\ad \vtilde_{S,\vecm}\bigr)^{\dagger}_{\Psi^{\ast}\htilde_{\beta,S,t}}
=\ad\bigl(-\rho(\Ad(g_{\beta,S,t}^{2})\vtilde_{S,\vecm})\bigr)$,
we have
\[
  \bigl[
 \ad(t\vtilde_{S,\vecm}),
 \ad(t\vtilde_{S,\vecm})^{\dagger}_{\Psi^{\ast}\htilde_{\beta,S,t}}
 \bigr]
 =
t^2\ad\Bigl(
 \bigl[
 \vtilde_{S,\vecm},
 -\rho\bigl(\Ad(g_{\beta,S,t}^{2})\vtilde_{S,\vecm}\bigr)
 \bigr]
 \Bigr).
\]
Because
\[
 \Ad(g_{\beta,S,t}^{2})\vtilde_{S,\vecm}=
 \sum_{\phi\in S}
 (-\log|z|^{2})^{-2}t^{-2}
 |z|^{-2\beta(\phi)}
 z^{m(\phi)}v_{\phi},
\]
we obtain
\begin{multline}
 t^2\bigl[
 \vtilde_{S,\vecm},
 -\rho( \Ad(g_{\beta,S,t}^{2})\vtilde_{S,\vecm})
 \bigr]
 = (\log|z|^2)^{-2}
 \left[
  \sum_{\phi\in S}
 z^{m(\phi)}
 v_{\phi},\,\,
 -\sum_{\phi\in S}
 \zbar^{m(\phi)}|z|^{-2\beta(\phi)}
 \rho(v_{\phi})
 \right]
 \\
=(\log|z|^2)^{-2}
 \sum_{\phi\in S}
 |z|^{2m(\phi)-2\beta(\phi)}
 \bigl[
 v_{\phi},
 -\rho(v_{\phi})
 \bigr]
=(\log|z|^2)^{-2}
 \sum_{\phi\in S}
 \bigl[
 v_{\phi},
 -\rho(v_{\phi})
 \bigr]
\\
 =(\log|z|^2) ^{-2}
 [v_S,-\rho(v_S)]
=(\log|z|^2)^{-2}2u_S.
\end{multline}
We obtain
\[
  \bigl[
 \ad(tv_S),\ad(tv_S)^{\dagger}_{\Psi^{\ast}\htilde_{\beta,S,t}}
 \bigr]
=(\log|z|^2) ^{-2}\ad(2u_S).
\]
Hence, we obtain
\[
 R(\Psi^{\ast}\htilde_{\beta,S,t})
 +\Bigl[
 \ad(tv_S)dz/z,
 \ad(tv_S)^{\dagger}_{\Psi^{\ast}\htilde_{\beta,S,t}}(d\zbar/\zbar)
 \Bigr]
 =\\
 -2\ad(u_S)\frac{dz\,d\zbar}{(\log|z|^2)^2|z|^2}
+2\ad(u_S)\frac{dz\,d\zbar}{(\log|z|^2)^{2}|z|^2}
=0.
\]
It means that
$(U^{\ast}\times\gminig,t\ad(\vtilde_{S,\vecm})dz/z,
\Psi^{\ast}(\htilde_{\beta,S,t}))$
is a harmonic bundle.
Hence, we obtain that
$(\gbigp_G,\theta_{S,\vecm,t},h_{\beta,S,t})$
is a harmonic $G$-Higgs bundle.
\hfill\qed

\subsection{Asymptotic behaviour in the wild case on punctured disc}

\subsubsection{Principal torus bundles on a punctured disc}

Let $U$ be a neighbourhood of $0$ in $\cnum$.
Let $U^{\ast}=U\setminus\{0\}$.
Let $\gbigp_T$ be a holomorphic principal $T$-bundle on $U^{\ast}$.
\begin{lem}
There exists an isomorphism $\gbigp_T\simeq U^{\ast}\times T$
as a holomorphic principal $T$-bundle.
\end{lem}
\pf
Because $T$ is connected,
$\gbigp_T$ is topologically isomorphic to
$U^{\ast}\times T$.
By Oka-Grauert principle,
$\gbigp_T$ is holomorphically isomorphic to
$U^{\ast}\times T$.
\hfill\qed

\subsubsection{The associated wild Higgs bundles}

Let $\theta$ be a section of 
$\gbigp_T(\gminig_1)\otimes\Omega^1_{U^{\ast}}$.
Let $(\gbigp_G,\theta)$
denote the induced $G$-Higgs bundle on $U^{\ast}$.
We assume the following.
\begin{condition}
$\gminio(\theta)$ is section of 
$\Upsilon(\gminig,\gminit)\otimes\nbigo_U(\ast 0)$,
and not constantly $0$.
\hfill\qed
\end{condition}

By shrinking $U$,
we may assume that $\gminio(\theta)$ is nowhere vanishing
on $U^{\ast}$.
Let $\ord(\gminio(\theta))$ is the zero order of
$\gminio(\theta)$ at $0$,
i.e.,
$\gminio(\theta)
=z^{\ord(\gminio(\theta))}\cdot \gamma$
for some nowhere vanishing holomorphic function $\gamma$.
We set
\[
 c(\theta)=\frac{-\ord(\gminio(\theta))}{\tth+1}.
\]

Let $f$ be the section of
$\gbigp_T(\gminig_1)$
determined by
$\theta=f\,dz/z$.
We obtain the characteristic polynomial
\[
 \det(t\id-\ad(f))
 =t^{\dim\gminig}+\sum_{j=0}^{\dim\gminig-1}
 a_j(z)t^j
 \in \nbigo_{U}(\ast 0)[t].
\]
The following is proved in
\cite[Lemma 6.2]{Kostant-TDS}.
\begin{lem}
\label{lem;25.5.3.10}
For any $u,u'\in\gminig_1^{\rs}$,
there exist $g\in T$
and $\lambda\in\cnum^{\ast}$
such that
$\Ad(g)(u)=\lambda u'$.
The number $\lambda$ satisfies
$\gminio(u)=\lambda^{\tth+1}\gminio(u')$.
\hfill\qed 
\end{lem}

We have
$\gminio(\theta)=\gminio(f)(dz/z)^{\tth+1}$.
Note that
$\gminio(f)z^{(c(\theta)-1)(\tth+1)}
(e_{-\psi}\otimes
\bigotimes_{i=1}^{\ell} e_{\alpha_i}^{\otimes\psi_i})^{-1}$
is holomorphic and nowhere vanishing on $U$.
Let $\gamma$ be
a nowhere vanishing holomorphic function on $U$
such that
$\gamma^{\tth+1}
(e_{-\psi}\otimes
\bigotimes_{i=1}^{\ell} e_{\alpha_i}^{\otimes\psi_i})
=
\gminio(f)z^{(c(\theta)-1)(\tth+1)}$.
The following lemma is easy to see.
\begin{lem}
\mbox{{}}\label{lem;25.9.25.10}
\begin{itemize}
 \item 
The roots of $\det(t\id-\ad(f))$ are
multi-valued meromorphic functions
       of the form $\beta z^{-c(\theta)+1}\gamma$
       for some complex numbers $\beta$.
 \item
      $a_j(z)$ $(j=0,\ldots,\dim\gminig-1)$
      are holomorphic at $z=0$
      if and only if
      $c(\theta)\leq 1$.
      Moreover,
       $a_j(0)=0$ $(j=0,\ldots,\dim\gminig-1)$
      if and only if
      $c(\theta)<1$.
       \hfill\qed
\end{itemize}
\end{lem}

\subsubsection{The case $c(\theta)\geq 1$}

Let $h$ be a harmonic metric
of $(\gbigp_G,\theta_G)$ compatible with $\Ad(\ttw)$.

\begin{prop}
If $c(\theta)>1$,
there exits $\epsilon>0$
such that  
\begin{equation}
\label{eq;25.5.10.3}
  \bigl|v(h^{\can},h)\bigr|_{h^{\can}}
  =O\bigl(
 \exp(-\epsilon|z|^{-c(\theta)+1})
  \bigr).
\end{equation}
If $c(\theta)=1$,
there exists $\epsilon>0$
such that
\begin{equation}
\label{eq;25.5.10.4}
  |v(h^{\can},h)|_{h^{\can}}
=O\bigl(|z|^{\epsilon}
  \bigr).
\end{equation}
In particular,
$h$ and $h^{\can}$ are mutually bounded
if $c(\theta)\geq 1$.
\end{prop}
\pf
We may assume that
$U=\{|z|<1\}$.
Let $\Utilde=\{\zeta\in\cnum\,|\,\Image(\zeta)>0\}$.
Let $\varphi:\Utilde\to U$
given by $\varphi(\zeta)=\exp(\sqrt{-1}\zeta)$.
We set
$\varphi^{\ast}(z)^{-c(\theta)+1}
=\exp\bigl(\sqrt{-1}(-c(\theta)+1)\zeta\bigr)$.

We set
$v_0=\sum_{i=1}^{\ell} e_{\alpha_i}+e_{-\psi}$.
We consider the semisimple endomorphism 
$\ad(v_0)$ on $\gminig$.
We obtain the eigen decomposition
\[
 \gminig
 =\bigoplus_{\beta\in\Sp(v_0)}
 \EE_{\beta}(\gminig,v_0).
\]
Here $\Sp(v_0)\subset\cnum^{\ast}$.
We set $m(\beta)=\dim \EE_{\beta}(\gminig,v_0)$.
We have the decomposition
\[
 \det(t\id-\ad(f))
 =t^{\dim\gminit}
 \cdot
 \prod_{\beta\in \Sp(v_0)}
 (t-\varphi^{\ast}(z)^{-c(\theta)+1}\varphi^{\ast}(\gamma)\beta)^{m(\beta)}.
\]
We set
$\betatilde=\varphi^{\ast}(z)^{-c(\theta)+1}\varphi^{\ast}(\gamma) \beta$.
By Lemma \ref{lem;25.5.3.10},
we obtain the eigen decomposition
\[
 \varphi^{\ast}\bigl(\gbigp_G(\gminig),\ad(f)\bigr)
 =(E_0,0)\oplus
 \bigoplus_{\beta\in\Sp(v_0)}
 (E_{\beta},\betatilde\id_{E_{\beta}}).
\]

Let us consider the case $c(\theta)>1$.
For $\Image(\zeta_0)>1$,
let $U(\zeta_0)=\{|\zeta-\zeta_0|<1\}$.
By \cite[Corollary 2.6]{Decouple},
we obtain that
\[
\Bigl|
\bigl[
\ad(\varphi^{\ast}f),
\ad(\varphi^{\ast}f)^{\dagger}_{\varphi^{\ast}(h)}
\bigr]
\Bigr|_{\varphi^{\ast}(h)}
=O\Bigl(
\exp\bigl(
 -\epsilon
 e^{(c(\theta)-1)\Image\zeta_0}
\bigr)
\Bigr).
\]
Hence, we obtain that
\[
 \bigl|
[\ad(f),\ad(f)^{\dagger}_{h}]
\bigr|_h=
O\bigl(
 \exp(-\epsilon |z|^{-c(\theta)+1})
\bigr).
\]
We obtain (\ref{eq;25.5.10.3})
by Proposition \ref{prop;25.6.18.50}.

Let us consider the case $c(\theta)=1$.
Let $\zeta_0\in\Utilde$
with $\Image(\zeta_0)>10$.
We consider
\[
 U'(\zeta_0)=\bigl\{
 |\zeta-\zeta_0|<\Image(\zeta_0)/2
 \bigr\}.
\]
Let $B=\{|w|<1\}$.
Let $\rho_{\zeta_0}:B\to U'(\zeta_0)$
be an isomorphism
$\rho_{\zeta_0}(w)=\zeta_0+w(\Image\zeta_0/2)$.
We obtain
\[
 \rho_{\zeta_0}^{\ast}\varphi^{\ast}
 \bigl(\gbigp_G(\gminig),\ad(f)dz/z\bigr)
 =(\rho_{\zeta_0}^{\ast}E_0,0)\oplus
 \bigoplus_{\beta\in\Sp(v_0)}
 \bigl(\rho_{\zeta_0}^{\ast}E_{\beta},
 (\Image\zeta_0/2)\rho_{\zeta_0}^{\ast}\betatilde\id_{E_i}dw
 \bigr)
\]
By \cite[Corollary 2.6]{Decouple},
there exists $\epsilon>0$,
which is independent of $\zeta_0$,
such that 
\[
\Bigl[
 \rho_{\zeta_0}^{\ast}
 \varphi^{\ast}\theta,
 \bigl(
 \rho_{\zeta_0}^{\ast}
 \varphi^{\ast}\theta
 \bigr)^{\dagger}_{\rho_{\zeta_0}^{\ast}\varphi^{\ast}h}
 \Bigr]
 =
 O\Bigl(
 \exp\bigl(-\epsilon\Image \zeta_0\bigr)
 \Bigr)
\]
on $\{|w|<1/2\}$.
We obtain
\[
\Bigl[
 \varphi^{\ast}\theta,
 \bigl(
 \varphi^{\ast}\theta
 \bigr)^{\dagger}_{\varphi^{\ast}h}
 \Bigr]
 =
 O\Bigl(
 \exp\bigl(-\epsilon'\Image \zeta_0\bigr)
 \Bigr) 
\]
on
$\bigl\{
|\zeta-\zeta_0|<\Image(\zeta_0)/4
\bigr\}$.
We obtain 
\[
 \bigl|
[\theta,\rho_h(\theta)]
\bigr|_h=O(|z|^{\epsilon''-2}\,dz\,d\zbar).
\]
Then, we obtain (\ref{eq;25.5.10.4})
from Proposition \ref{prop;25.6.18.50}.
\hfill\qed

\vspace{.1in}
Let $f_{\phi}$ $(\phi\in\Pi^Q)$
be determined by
$\theta_{\phi}=f_{\phi}dz/z$.
Because $h$ and $h^{\can}$ are mutually bounded,
we obtain the following.
\begin{lem}
For any $\phi\in \Pi^{Q}$,
there exists $C\geq 1$ such that
$C^{-1}\leq
 |f_{\phi}|_h\cdot |z|^{c(\theta)-1}\leq C$.
\hfill\qed
\end{lem}

\subsubsection{The case $c(\theta)<1$}
\label{subsection;25.5.9.2}

Let us consider the case $c(\theta)<1$.
By \cite[Theorem 1]{s2} with Lemma \ref{lem;25.9.25.10},
we obtain
$|f|_{h}
=O\bigl(
(-\log|z|)^{-1}
\bigr)$.
For $\phi\in \Pi^Q$,
we obtain
$|f_{\phi}|_h
=O\bigl(
(-\log|z|)^{-1}
\bigr)$.

We set
$E=\gbigp_T(\gminig)$,
$E_{\phi}=\gbigp_T(\gminig_{\phi})$
$(\phi\in\Delta)$
and
$E_{\gminit}=\gbigp_T(\gminit)$.
We have the orthogonal decomposition
\begin{equation}
\label{eq;25.5.7.20}
E=
E_{\gminit}
\oplus
\bigoplus_{\phi\in\Delta}
E_{\phi}.
\end{equation}
We obtain the filtered bundle $\nbigp^h_{\ast}E$ on $(U,0)$.
(See \cite[\S2.5]{Mochizuki-KH-Higgs}.)
Because the decomposition (\ref{eq;25.5.7.20}) is orthogonal,
we obtain the decomposition of filtered bundles
\[
 \nbigp^h_{\ast}E
 =\bigoplus_{\phi\in\Delta} \nbigp^h_{\ast}E_{\phi}
 \oplus
 \nbigp^h_{\ast}E_{\gminit}.
\]

For $\phi=\sum \phi_i\alpha_i\in\Delta$,
we have the isomorphism
\begin{equation}
\label{eq;25.9.25.11}
 E_{\phi}\simeq \bigotimes_{i=1}^{\ell}E_{\alpha_i}^{\otimes\phi_i}.
\end{equation}
Because $h$ is induced by a harmonic reduction
$\gbigp^h_{T^c}\subset\gbigp_{T}$,
the isomorphism (\ref{eq;25.9.25.11})
is an isometry with respect to the metrics induced by $h$.
Hence, we have the isomorphism of filtered bundles
\[
 \nbigp^h_{\ast}E_{\phi}
 \simeq
 \bigotimes_{i=1}^{\ell}
 \bigl(
 \nbigp^h_{\ast}E_{\alpha_i}
 \bigr)^{\otimes\phi_i}.
\]

Because $|f|_h=O((-\log|z|)^{-1})$,
$f$ is a section of $\nbigp^h_0E$.
For $\phi\in\Pi^Q$,
let $a(\phi)\leq 0$ be the real number determined
as the parabolic degree of $f_{\phi}$:
\[
 a(\phi)=\min\{
 c\in\real\,\big|\,
 f_{\phi}\in
 \nbigp_cE_{\phi}
 \}.
\]
By the norm estimate in \cite{s1},
for $\phi\in \Pi^Q$,
we obtain
\begin{equation}
\label{eq;25.9.26.20}
 \log|f_{\phi}|_{h}
 +\log|z|^{a(\phi)}
 =O\Bigl(
\log(-\log|z|^2)
 \Bigr).
\end{equation}
Because
$\gminio(\theta)=\bigotimes_{i=1}^{\ell}f_{\alpha_i}^{\otimes\psi_i}
\otimes f_{-\psi}\cdot(dz/z)^{\tth+1}$,
we have
\[
 \sum_{i=1}^{\ell} \psi_ia(\alpha_i)+a(-\psi)
+\ord(\gminio(\theta))+(\tth+1)=0.
\]
The condition $c(\theta)<1$
is equivalent to
$m:=\ord(\gminio(\theta))+(\tth+1)>0$.
In particular,
$\sum \psi_ia(\alpha_i)+a(-\psi)<0$.

We set
$S=\bigl\{
 \phi\in \Pi^Q\,\big|\,
 a_{\phi}=0
 \bigr\}\subsetneq\Pi^Q$.
By changing a holomorphic trivialization
of $\gbigp_T\simeq U^{\ast}\times T$,
we may assume that
there exist 
$u_S\in\gminit_{\real}$
and $v_S=\sum_{\phi\in S}v_{\phi}\in\bigoplus_{\phi\in S}\gminig_{\phi}$
as in Proposition \ref{prop;25.5.8.1}
such that the following holds for some $t>0$:
\begin{itemize}
 \item $f_{\alpha_i}=tv_{\alpha_i}$
       $(\alpha_i \in S\setminus\{-\psi\})$.
 \item $f_{\alpha_i}=\gamma_ie_{\alpha_i}$ $(\alpha_i\not\in S)$
       for a nowhere vanishing holomorphic function
       $\gamma_i$ on $U$.
 \item (the case $-\psi\in S$)
       $f_{-\psi}=tz^mv_{-\psi}$.
 \item (the case $-\psi\not\in S$)
       $f_{-\psi}=\gamma_{-\psi}z^me_{-\psi}$
       for a nowhere vanishing holomorphic function $\gamma_{-\psi}$
       on $U$.
 \end{itemize}
We set $\beta:=\sum a(\alpha_i)\epsilon_i$.
\begin{lem}
We have 
$\beta(-\psi)
=-\sum a(\alpha_i)\psi_i
=m+a(-\psi)$.
In particular,
if $a(-\psi)=0$,
we have
$\beta(-\psi)=m$.
\hfill\qed
\end{lem}

Let $g:U^{\ast}\to T_{\real}$ be determined as follows:
\[
 g=\exp\Bigl(
 -\beta\log|z|
 -u_S\log(-\log|z|^{2t})
 \Bigr).
\]
We obtain the Hermitian metric
$h_0=g^{-1}Kg$ of $\gbigp_G$.
Let $\rho_{h_0}$ be the anti-$\cnum$-linear automorphism of
$\gbigp_G(\gminig)=\gbigp^{h_0}_{T_c}(\gminig)$ and $\rho$.
Let $R(h_0)$ denote the curvature of the Chern connection of
$(\gbigp_G,h_0)$.

\begin{prop}
$h$ and $h_0$ are mutually bounded.
Moreover,
there exists $\epsilon>0$ such that
the following holds with respect to $h_0$:
\begin{equation}
\label{eq;25.9.26.12}
       R(h_0)+[\theta,-\rho_{h_0}(\theta)]
       =O\bigl(|z|^{\epsilon-2}dz\,d\zbar\bigr).
\end{equation}
\end{prop}
\pf
By using Lemma \ref{lem;25.9.25.20},
we obtain
$\log|e_{\alpha_i}|_{h_0}
+\log|z|^{a(\alpha_i)}
=O\bigl(\log(-\log|z|^2)\bigr)$.
We obtain
\begin{equation}
\label{eq;25.9.26.10}
 \log|f_{\alpha_i}|_{h_0}
 +\log|z|^{a(\alpha_i)}
 =O\bigl(\log(-\log|z|^2)\bigr).
\end{equation}
We also obtain
\begin{equation}
\label{eq;25.9.26.11}
 \log|f_{-\psi}|_{h_{0}}
 +\log|z|^{a(-\psi)}
\sim\log|e_{-\psi}|_{h_{0}}
+\log|z|^{a(-\psi)+m}
=\log|e_{-\psi}|_{h_{0}}
+\log|z|^{-\sum \psi_ia(\alpha_i)}
=O\bigl(\log(-\log|z|^2)\bigr).
\end{equation}

We set
$\theta_0=\sum_{\phi\in S}f_{\phi}dz/z$.
By Proposition \ref{prop;25.9.25.30},
we have
$R(h_0)+[\theta_0,-\rho_{h_0}(\theta_0)]=0$.
By (\ref{eq;25.9.26.10}) and (\ref{eq;25.9.26.11}),
we obtain the following estimate with respect to $h_0$:
\[
 [\theta,-\rho_{h_0}(\theta)]
-[\theta_0,-\rho_{h_0}(\theta_0)]
=O\bigl(
 |z|^{\epsilon-2}\,dz\,d\zbar
 \bigr).
\]
We obtain the estimate (\ref{eq;25.9.26.12}).

By the estimate (\ref{eq;25.9.26.20}),
for any $\phi\in\Delta$,
we have
\[
 \log|e_{\phi}|_{h}
 +\log|z|^{\beta(\phi)}
=O\bigl(\log(-\log|z|)\bigr).
\]
By the estimates (\ref{eq;25.9.26.10}), (\ref{eq;25.9.26.11}),
we have
\[
 \log|e_{\phi}|_{h_0}
 +\log|z|^{\beta(\phi)}
=O\bigl(\log(-\log|z|)\bigr).
\]
We obtain
\[
 \log \Tr \Ad(s(h_0,h))
 =O\bigl(
 \log(-\log|z|)
 \bigr).
\]
By (\ref{eq;25.9.26.12}) and \cite[Lemma 3.1]{s1},
we also have
\[
 -\del_z\del_{\zbar}\Tr \Ad s(h_0,h)
=O\bigl(|z|^{\epsilon-2}\bigr).
\]
There exists an $L_2^p$-function $\gamma$
such that
$-\del_z\del_{\zbar}\gamma
= -\del_z\del_{\zbar}\Tr\Ad s(h_0,h)$.
Then,
$\Tr\Ad(s(h_0,h))-\gamma$
is harmonic on $U^{\ast}$,
and
$O\bigl(
\log(-\log|z|)
\bigr)$ holds.
It implies that
the maximum principle holds for 
$\Tr\Ad(s(h_0,h))-\gamma$
on $U$.
In particular, we obtain that
$\Tr\Ad(s(h_0,h))$ is bounded.
We obtain that $h$ and $h_0$ are mutually bounded.
\hfill\qed

\subsubsection{Auxiliary metrics}

Let us give a complement.
Let $\theta=f\,dz/z$
be a holomorphic section of
$\gbigp_T(\gminig)\otimes\Omega^1_{U^{\ast}}$.

\begin{condition}
Assume that $\gminio(\theta)$ is a section of
$\Upsilon(\gminig,\gminit)\otimes\Omega^1_{U}(\ast 0)$.
\hfill\qed
\end{condition}

We have the decomposition
$f=\sum_{i=1}^{\ell} f_{\alpha_i}+f_{-\psi}$,
where $f_{\phi}$ are sections of
$\gbigp_T(\gminig_{\phi})$.
There exists a holomorphic isomorphism
$\gbigp_T\simeq U^{\ast}\times T$
such that
$f_{\alpha_i}=e_{\alpha_i}$ $(i=1,\ldots,\ell)$
under the induced holomorphic isomorphisms
\begin{equation}
\label{eq;25.5.9.1}
 \gbigp_T(\gminig_{\phi})\simeq
 \gminig_{\phi}\otimes\nbigo_{U^{\ast}}
 \quad
 (\phi\in\Delta).
\end{equation}

By the construction in \S\ref{subsection;25.5.9.2},
we obtain the following lemma.
\begin{lem}
\label{lem;25.5.10.10}
Suppose that $\tth+1+\ord\gminio(\theta)>0$.
Let $\beta\in\gminit_{\real}$
such that the following holds.
\begin{itemize}
 \item $\alpha_i(\beta)\leq 0$ $(i=1,\ldots,\ell)$.
 \item $\psi(\beta)+
       (\tth+1+\ord\gminio(\theta))\geq 0$.
\end{itemize}
Then, there exists a Hermitian metric $h_{\beta}$ of
$\gbigp_G$ compatible with $\Ad\ttw$
satisfying the following conditions.
\begin{itemize}
 \item For any $\phi\in \Pi^Q=\{\alpha_i\}\cup\{-\psi\}$,
       the estimate
      $\log|f_{\phi}|+\phi(\beta)\log|z|
      =O\bigl(\log(-\log|z|)\bigr)$ holds.
 \item There exists $\epsilon>0$ such that
       $R(h_{\beta})+[\theta,-\rho_{h_{\beta}}(\theta)]=
       O\bigl(|z|^{\epsilon-2}\bigr)\,dz\,d\zbar$
       with respect to $h_{\beta}$.
\hfill\qed
\end{itemize} 
\end{lem}

\subsection{Classification}

Let $X$ be a compact Riemann surface
with a finite subset $D\subset X$.
Let $\gbigp_T$ be a holomorphic principal $T$-bundle on $X\setminus D$.
Let $\theta$ be a section of $\gbigp_T(\gminig_1)$.
Assume the following condition.
\begin{condition}
 $\gminio(\theta)$ is a section of
 $\Upsilon(\gminig,\gminit)\otimes K_X^{\tth+1}(\ast D)$.
 It is not constantly $0$.
\hfill\qed
\end{condition}

Let $D^{>0}\subset D$ denote the set of $P\in D$ such that
$\ord_P\gminio(\theta)+\tth+1>0$.
Let $\nbigs(\gbigp_T,\theta)$ denote
the set of tuples
$(\beta_P\in\gminit_{\real}\,|\,P\in D^{>0})$
satisfying the following conditions:
\begin{equation}
\label{eq;25.9.26.25}
\alpha_i(\beta_P)\leq 0,
\quad
\psi(\beta_P)+(\tth+1+\ord_P\gminio(\theta))\geq 0.
\end{equation}
If $D^{>0}=\emptyset$,
$\nbigs(\gbigp_T,\theta)$ is the set of one point.

Let $h$ be any harmonic metric of
$(\gbigp_G,\theta)$ compatible with $\Ad\ttw$.
By the norm estimate in \cite{s1},
for any $P\in D^{>0}$,
we obtain $\beta_P(h)\in\gminit_{\real}$
determined by 
\[
  \log |f_{\alpha_i}|_{h}
 +\alpha_i(\beta_P(h))\log|z_P|
 =O\bigl(
 \log(-\log|z_P|)
 \bigr)
 \quad(i=1,\ldots,\ell).
\]
\begin{lem}
$\beta_P(h)$ satisfies the condition {\rm(\ref{eq;25.9.26.25})}.
\end{lem}
\pf
Let $(X_P,z_P)$ be a holomorphic coordinate neighbourhood
around $P$.
We set $X_P^{\ast}=X_P\setminus\{P\}$.
Because $P\in D^{<0}$,
we have
$|f_{\phi|X_P^{\ast}}|_h=O\bigl((-\log|z_P|)^{-1}\bigr)$
for $\phi\in\Pi^Q$.
Hence, we obtain
$\alpha_i(\beta_P(h))\leq 0$
$(i=1,\ldots,\ell)$.
By setting $m=\ord\gminio(\theta)+\tth+1$,
we obtain
\[
 \log|f_{|X_P^{\ast}}|_h
 +\log|z|^{-m}
 +\log|z|^{-\psi(\beta_P(h))}
=O\bigl(\log(-\log|z|)\bigr).
\]
We obtain
$-m-\psi(\beta_P(h))\leq 0$.
\hfill\qed

\vspace{.1in}

Let $\Harm(\gbigp_G,\theta,\Ad\ttw)$
be the set of harmonic metrics of
$(\gbigp_G,\theta)$
compatible with $\Ad\ttw$.
By the previous lemma,
we obtain the map
\[
 \Psi:
 \Harm(\gbigp_G,\theta,\Ad\ttw)
 \lrarr
 \nbigs(\gbigp_T,\theta),
 \quad
 \Psi(h)=\bigl(\beta_P(h)\bigr)_{P\in D^{>0}}.
\]

The following theorem is a generalization of
\cite[Corollary 1.11]{Li-Mochizuki2}.
\begin{thm}
\label{thm;25.5.10.20}
The map $\Psi$ is a bijection.
In particular,
if $\ord_P\gminio(\theta)+\tth+1\leq 0$ for any $P\in D$,
then $(\gbigp_G,\theta)$
has a unique harmonic metric compatible with $\Ad\ttw$.
\end{thm}
\pf
If $\Psi(h_1)=\Psi(h_2)$,
we obtain $h_1=h_2$
by using Proposition \ref{prop;25.5.9.3}.

Let us study the surjectivity of $\Psi$.
Let $(\beta_P\,|\,P\in D^{>0})\in\nbigs(\gbigp_T,\theta)$.
Let $(X_P,z_P)$ be a holomorphic coordinate neighbourhood
around $P\in D$ such that $X_P\simeq \{|z_P|<2\}$.
For any $0<r<2$, we set $X_P(r)=\{|z_P|<r\}$.
By Proposition \ref{prop;25.6.19.1} and Lemma \ref{lem;25.5.10.10},
there exists a Hermitian metric $h_0$ of $\gbigp_G$
compatible with $\Ad\ttw$
satisfying the following conditions.
\begin{itemize}
 \item On $X\setminus \bigcup_{P\in D^{>0}}X_P(1/2)$,
       $h$ is a harmonic metric of
       $(\gbigp_G,\theta)$
       compatible with $\Ad\ttw$.
 \item If $P\in D\setminus D^{>0}$,
       on $X_P(1/3)\setminus\{P\}$,
       $h_0$ equals the canonical decoupled harmonic metric
       of $(\gbigp_G,\theta)$.
 \item If $P\in D^{>0}$,
       on $X_P\setminus\{P\}$, we have the following estimates:
\[
       \log |f_{\alpha_i}|_{h_{0}}
       +\alpha_i(\beta_P) \log|z_P|
       =O\bigl(
       \log(-\log|z_P|)
       \bigr)
       \quad(i=1,\ldots,\ell),
\]       
\[
       R(h_0)+[\theta,-\rho_{h_0}(\theta)]
       =O(|z|^{\epsilon-2}dz_P\,d\zbar_P).
\]       
\end{itemize}

Let $X_i$ $(i=1,2,\ldots)$
be a smooth exhaustive family of $X\setminus D$
such that
$\del X_i\subset \bigcup_{P\in D}X_P(1/4)$.
Let $h_i$ be harmonic metrics of
$(\gbigp_G,\theta)$ compatible with $\Ad\ttw$
such that
$h_{i|\del X_i}=h_{0|\del X_i}$.
By Proposition \ref{prop;25.6.20.30},
there exists $C>0$
such that the following holds on
$X\setminus\bigcup_{P\in D}X_P(1)$:
\[
 \bigl|
 \Ad(s(h_0,h_i))
 \bigr|_{h_0}
 \leq C.
\]
We set $F(h_0)=R(h_0)+[\theta,-\rho_{h_0}(\theta)]$.
Let $g_X$ be a conformal metric of $X$.
By \cite[Lemma 3.1]{s1},
we have
\[
 \Delta_{g_X}\log
 \Tr\bigl(
\Ad(s(h_0,h_i))
 \bigr)
 \leq
 \bigl|
 F(h_0)
 \bigr|_{h_0,g_X}
\]
on $X_i$.
There exists an $L_2^p$-function $\gamma_P$
on $X_P$
such that
$\Delta_{g_X}\log\gamma_P
=|F(h_0)|_{h_0,g_X}$.
Then,
$\Tr\Ad((s(h_0,h_i)))-\gamma_P$
are subharmonic functions on
$X_P\cap X_i$.
Note that $|\gamma_P|\leq C'$ for some $C'>0$.
We have $\Tr s(h_0,h_i)=\dim\gminig$ on $\del X_i\cap X_P$.
We also have
$\Tr(s(h_0,h_i))\leq C$ on $X_i\setminus \bigcup_{P\in D^{>0}}X_P(1)$.
Hence, there exists $C_P>0$, independently from $i$,
such that the following holds on $X_P\cap X_i$:
\[
\Tr(\Ad(s(h_0,h_i)))-\gamma_P
\leq C_P.
\]
There exists $C>0$,
independently from $i$,
such that
$\Tr \Ad(s(h_0,h_i))
\leq C$
on $X_i$.
By taking a subsequence,
we may assume that
the sequence $h_i$ is convergent
on any relatively compact open subsets of $X\setminus D$
in the $C^{\infty}$-sense.
Let $h_{\infty}$ denote the limit.
We obtain that
$\Tr s(h_0,h_{\infty})\leq C$.
Hence, $h_0$ and $h_{\infty}$ are mutually bounded.
We obtain that
$\Psi(h_{\infty})=(\beta_P\,|\,P\in D^{>0})$.
\hfill\qed

\subsubsection{Dirichlet problem}

Let $X$ be any Riemann surface.
Let $D$ be a discrete subset of $X$.
Let $Y\subset X$ be a relatively compact open subset
with smooth boundary $\del Y$
such that $D\cap \del Y=\emptyset$.

Let $\gbigp_T$ be a holomorphic principal $T$-bundle
on $X\setminus D$.
Let $\theta$ be a section of
$\gbigp_T(\gminig_1)\otimes K_{X\setminus D}$
such that
$\gminio(\theta)$ is a section of
$\Upsilon(\gminig,\gminit)\otimes K_X^{\tth+1}(\ast D)$
which is not constantly $0$.

We define $\nbigs(\gbigp_T,\theta,Y)$ as
the set of tuples
$\beta_P\in\gminit_{\real}$ $(P\in D^{>0}\cap Y)$
such that
\[
 \alpha_i(\beta_P)\leq 0,
 \quad\quad
 \psi(\beta_P)
 +(\tth+1+\ord_P\gminio(\theta))\geq 0.
\] 

Let $h_{\del Y}$ be a Hermitian metric of $\gbigp_G$
compatible with $\Ad\ttw$ on $\del Y$.
Let $\Harm(\gbigp_G,\theta,\Ad\ttw;h_{\del Y})$ denote the set of
harmonic metrics $h$ of
$(\gbigp_G,\theta)$ compatible with $\Ad\ttw$
on $Y$
such that $h_{|\del Y}=h_{\del Y}$.
The following theorem is similar to Theorem \ref{thm;25.5.10.20}.
\begin{thm}
The natural map
$\Harm(\gbigp_T,\theta,\Ad\ttw;h_{\del Y})
\to \nbigs(\gbigp_T,\theta,Y)$ 
is a bijection.
\hfill\qed
\end{thm}

\subsubsection{Case of $\cnum$}

Let us consider the case $X=\proj^1$
and $D=\{\infty\}$.
Let $(\gbigp_T,\theta)$ on $\cnum$
such that $\gminio(\theta)$ is not constantly $0$,
and meromorphic at $\infty$.
In this case, we have $D^{>0}=\emptyset$.
We obtain the following from Theorem \ref{thm;25.5.10.20}.
The following is a generalization of a result in \cite{Li}.
\begin{cor}
There exists a unique harmonic metric of
$(\gbigp_G,\theta)$  compatible with $\Ad\ttw$.
\hfill\qed
\end{cor}

\subsubsection{Homogeneous case on $\cnum^{\ast}$}

Let $e_{\phi}$ $(\phi\in\Delta)$ be a tuple of root vectors
as in Proposition \ref{prop;25.5.3.3}.
We consider the principal $T$-bundle
$\gbigp_T=\cnum^{\ast}\times T$ on $\cnum^{\ast}$
with the Higgs field
\[
 \theta=
 \sum_{i=1}^{\ell} e_{\alpha_i}\frac{dz}{z}
 +z^{\tth+1}e_{-\psi}\frac{dz}{z}.
\]

For any $a\in\cnum^{\ast}$,
we consider the map
$F_a:\proj^1\to\proj^1$
defined by
$F_a(z)=az$.
We also consider $g_a\in T$
such that
$g_ae_{\alpha_i}=ae_{\alpha_i}$
for $i=1,\ldots,\ell$.
Because $\gbigp_T=\cnum^{\ast}\times T$,
we have the natural isomorphism
$\kappa_a:F_a^{\ast}\gbigp_T\simeq \gbigp_T$.
Moreover,
by multiplying $g_a$,
we obtain 
\[
 \kappatilde_a=g_a\cdot \kappa_a:
F_a^{\ast}(\gbigp_T)\simeq \gbigp_T.
\]
By this isomorphism,
$F_a^{\ast}\theta$ is identified with
$a\theta$.
We set $S^1=\{a\in\cnum\,|\,|a|=1\}$.
We obtain the following proposition
by the uniqueness of the solutions.
\begin{lem}
Any $h\in\Harm(\gbigp_G,\theta,\Ad\ttw)$
is $S^1$-invariant.
\hfill\qed 
\end{lem}

\begin{cor}
\label{cor;26.3.8.1}
$S^1$-invariant harmonic metrics of 
$(\gbigp_G,\theta)$ compatible with $\Ad\ttw$
bijectively correspond to  
$\beta\in\gminit_{\real}$
satisfying 
$\alpha_i(\beta)\leq 0$ $(i=1,\ldots,\ell)$
and
$\psi(\beta)+\tth+1\geq 0$.
\hfill\qed
\end{cor}

\begin{rem}
The author obtained a classification like 
Corollary {\rm\ref{cor;26.3.8.1}}
in the case $\gminig=\minisl(n,\cnum)$
in {\rm\cite{Mochizuki-Toda}}.
Though it was accepted for publication after a peer review
in June {\rm 2016},
it has been left unpublished since then, unfortunately.
In the companion paper {\rm\cite{Mochizuki-TodaII}},
we also studied the associated meromorphic connections
and their Stokes structures.
\hfill\qed
\end{rem}

\begin{rem}
Guest, Its and Lin
{\rm\cite{GuestItsLin, GuestItsLin2, GuestItsLin3, GuestLin}}
developed the theory of $tt^{\ast}$-Toda equation
on $\cnum^{\ast}$ 
in the case $\gminig=\minisl(n,\cnum)$ 
which includes a similar classification,
on the basis of isomonodromic deformations.
Their method is completely different from ours.
Guest and Ho developed the theory of the $tt^{\ast}$-Toda equation
for general semisimple Lie algebra $\gminig$
in {\rm\cite{Guest-Ho-2019}}.
More recently,
in a conference in January {\rm 2026} held
at the University of Mannheim,
the author was informed that
they established a similar classification
for general $G$.
The author expects that
the study in this paper would be useful to clarify
a complementary aspect of
$tt^{\ast}$-Toda equations.
\hfill\qed
\end{rem}

\end{document}

%% file: notation.tex
\newcommand{\nbiga}{\mathcal{A}}
\newcommand{\nbigb}{\mathcal{B}}
\newcommand{\nbigc}{\mathcal{C}}

\newcommand{\nbigg}{\mathcal{G}}

\newcommand{\nbigk}{\mathcal{K}}

\newcommand{\nbigo}{\mathcal{O}}
\newcommand{\nbigp}{\mathcal{P}}

\newcommand{\nbigs}{\mathcal{S}}

\newcommand{\proj}{\mathbb{P}}
\newcommand{\seisuu}{{\mathbb Z}}

\newcommand{\cnum}{{\mathbb C}}
\newcommand{\real}{{\mathbb R}}

\newcommand{\hyperk}{\mathbb{K}}

\newcommand{\DD}{\mathbb{D}}
\newcommand{\EE}{\mathbb{E}}

\newcommand{\gbigp}{\mathfrak P}

\newcommand{\gminia}{\mathfrak a}

\newcommand{\gminig}{\mathfrak g}
\newcommand{\gminih}{\mathfrak h}

\newcommand{\gminik}{\mathfrak k}
\newcommand{\gminil}{\mathfrak l}

\newcommand{\gminio}{\mathfrak o}

\newcommand{\gminis}{\mathfrak s}
\newcommand{\gminit}{\mathfrak t}
\newcommand{\gminiu}{\mathfrak u}

\newcommand{\veca}{{\boldsymbol a}}
\newcommand{\vecb}{{\boldsymbol b}}

\newcommand{\vecm}{{\boldsymbol m}}

\newcommand{\vecW}{{\boldsymbol W}}

\newcommand{\vecU}{{\boldsymbol U}}
\newcommand{\vecV}{{\boldsymbol V}}

\newcommand{\lrarr}{\longrightarrow}

\newcommand{\pf}{{\bf Proof}\hspace{.1in}}
\newcommand{\qed}{\mbox{\rule{1.2mm}{3mm}}}

\def\Hom{\mathop{\rm Hom}\nolimits}

\def\End{\mathop{\rm End}\nolimits}

\def\Image{\mathop{\rm Im}\nolimits}

\def\Re{\mathop{\rm Re}\nolimits}

\def\GL{\mathop{\rm GL}\nolimits}
\def\SL{\mathop{\rm SL}\nolimits}

\def\rank{\mathop{\rm rank}\nolimits}
\def\Spec{\mathop{\rm Spec}\nolimits}

\def\Ker{\mathop{\rm Ker}\nolimits}

\def\ad{\mathop{\rm ad}\nolimits}

\def\Rep{\mathop{\rm Rep}\nolimits}
\def\ord{\mathop{\rm ord}\nolimits}

\def\Tr{\mathop{\rm Tr}\nolimits}

\def\can{\mathop{\rm can}\nolimits}

\def\id{\mathop{\rm id}\nolimits}

\def\Irr{\mathop{\rm Irr}\nolimits}

\newcommand{\del}{\partial}
\newcommand{\delbar}{\overline{\del}}

\newcommand{\barz}{\overline{z}}
\newcommand{\zbar}{\barz}

\newcommand{\Sp}{{\mathcal Sp}}

\newcommand{\lefttop}[1]{{}^{#1}\!}

\def\Harm{\mathop{\rm Harm}\nolimits}

\newcommand{\gl}{\gminig\gminil}
\newcommand{\minisl}{\gminis\gminil}

\newcommand{\Etilde}{\widetilde{E}}

\newcommand{\kappatilde}{\widetilde{\kappa}}

\newcommand{\htilde}{\widetilde{h}}

\newcommand{\vtilde}{\widetilde{v}}

\newcommand{\Psitilde}{\widetilde{\Psi}}

\newcommand{\Utilde}{\widetilde{U}}

\def\ord{\mathop{\rm ord}\nolimits}
\def\Herm{\mathop{\rm Herm}\nolimits}

\def\Ad{\mathop{\rm Ad}\nolimits}
\def\rs{\mathop{\rm rs}\nolimits}

\newcommand{\Xbar}{\overline{X}}

\newcommand{\II}{\mathbb{I}}

\newcommand{\ttv}{{\tt v}}

\newcommand{\ttw}{{\tt w}}
\newcommand{\tth}{{\tt h}}
\newcommand{\ttx}{{\tt x}}

\newcommand{\betatilde}{\widetilde{\beta}}
\newcommand{\vecVtilde}{\widetilde{\vecV}}

%% file: new_theorem.tex
\newtheorem{thm}{Theorem}[section]
\newtheorem{cor}[thm]{Corollary}

\newtheorem{rem}[thm]{Remark}
\newtheorem{lem}[thm]{Lemma}
\newtheorem{prop}[thm]{Proposition}
\newtheorem{df}[thm]{Definition}

\newtheorem{condition}[thm]{Condition}

\newtheorem{notation}[thm]{Notation}

%% file: 2-14.bbl
\begin{thebibliography}{99}

\bibitem{Baraglia} D. Baraglia,
	 {\em $G_2$ Geometry and integrable system},
	 thesis, arXiv:1002.1767v2, 2010.

\bibitem{b}
O. Biquard,
{\it Fibr\'es de Higgs et connexions int\'egrables:
le cas logarithmique (diviseur lisse)},
Ann. Sci. \'Ecole Norm. Sup. {\bf 30} (1997), 41--96.

\bibitem{Biquard-Boalch}
	O. Biquard and P. Boalch
	{\em Wild non-abelian Hodge theory on curves,}
	Compos. Math. {\bf 140} (2004), 179--204.

\bibitem{Biquard-Garcia-Prada-Mundet}
	O. Biquard,
 O. Garc\'{\i}a-Prada, and
	I. Mundet i Riera,
	{\em Parabolic Higgs bundles and representations of
	the fundamental group of a punctured surface into a real group},
	Adv. Math. {\bf 372} (2020), 107305, 70 pp.

\bibitem{Bradlow-Garcia-Prada-Mundet}
	S. Bradlow,
	O. Garc\'{\i}a-Prada, and
	I. Mundet i Riera,
	{\em Relative Hitchin-Kobayashi correspondences
	for principal pairs},
	Q. J. Math. {\bf 54} (2003), 171--208.
	
\bibitem{corlette}
K. Corlette,
{\em
Flat $G$-bundles with canonical metrics},
J. Differential Geom. {\bf 28} (1988), 361--382.
	
\bibitem{DMOS}
	P. Deligne,
	J. S. Milne,
	A. Ogus and
	K.-Y.Shih,
	{\em Hodge cycles, motives, and Shimura varieties},
	Lecture Notes in Math., {\bf 900}
	Springer-Verlag, Berlin-New York, 1982. ii+414 pp.

\bibitem{don2}
S. K. Donaldson,
{\em Twisted harmonic maps and the self-duality equations},
Proc. London Math. Soc. {\bf 55} (1987), 127--131.
	
\bibitem{Donaldson-boundary-value}
S. K. Donaldson,
{\em Boundary value problems for Yang-Mills fields},
J. Geom. Phys. {\bf 8} (1992), 89--122.

\bibitem{Garcia-Prada-Gonzalez-cyclic}
	 O. Garc\'{\i}a-Prada and
	 M. Gonz\'{a}lez,
	 {\em Cyclic Higgs bundles and the Toledo invariant},
	 arXiv:2403.00415
\bibitem{Guest-Ho-2019}
	M. A. Guest,
	N.-K. Ho,
	{\em Kostant, Steinberg, and the Stokes matrices of
	the $tt^{\ast}$-Toda equations},
	Selecta Math. (N.S.) {\bf 25} (2019), no. 3, Paper No. 50, 43 pp.

\bibitem{GuestItsLin} M. A. Guest, A. R. Its and C. S. Lin,
	 {\em Isomonodromy aspects of the $tt^*$ equations of
	 Cecotti and Vafa I.
	 Stokes data},
	 Int. Math. Res. Notices 2015 (2015), 11745--11784.

\bibitem{GuestItsLin2} M. A. Guest, A. R. Its and C. S. Lin,
	 Isomonodromy aspects of the $tt^{\ast}$ equations
	 of Cecotti and Vafa II:
	 Riemann-Hilbert problem
	 Comm. Math. Phys. {\bf 336} (2015), 337--380.
	 
\bibitem{GuestItsLin3}
	M. A. Guest, 
	A. R. Its, 
	C.-S. Lin,
	{\em Isomonodromy aspects of the $tt^{\ast}$ equations
	of Cecotti and Vafa III:
	Iwasawa factorization and asymptotics,}
	Comm. Math. Phys. {\bf 374} (2020), 923--973.


 \bibitem{GuestLin} M. A. Guest and C. S. Lin,
	 {\it Nonlinear PDE aspects of the $tt^*$ equations of Cecotti and Vafa},
	 J. reine angew. Math. {\bf 689} (2014) 1--32.
	 

 \bibitem{Helgason1978}
	 S. Helgason
	 {\em Differential geometry, Lie groups, and symmetric spaces},
	 Academic Press, Inc.
	 1978, xv+628 pp.
\bibitem{Hitchin-self-duality}
	 N. J. Hitchin,
	 {\em The self-duality equations on a Riemann surface},
	 Proc. London Math. Soc. (3) {\bf 55} (1987), 59--126. 
\bibitem{Hitchin-section}
	 N. J. Hitchin,
	 {\em Lie groups and Teichm\"{u}ller space},
	 Topology {\bf 31} (1992), 449--473. 
\bibitem{Humphreys}
J. E. Humphreys,
{\em Introduction to Lie algebras and representation theory}.
Second printing, revised
Grad. Texts in Math., {\bf 9},
Springer-Verlag, New York-Berlin, 1978. xii+171 pp.

\bibitem{JZ2}
J. Jost and K. Zuo,
{\em Harmonic maps of infinite energy and
rigidity results for representations of fundamental
   groups of quasiprojective varieties},
J. Differential Geom. {\bf 47} (1997), 469--503.
	
\bibitem{Knapp-Book}
	A. Knapp,
	{\em Lie groups beyond an introduction.}
	Progress in Mathematics, {\bf 140},
	Birkh\"{a}user Boston, Inc., Boston, MA, 1996.

\bibitem{Kostant-TDS}
	 B. Kostant,
	{\em The principal three-dimensional subgroup and
	the Betti numbers of a complex simple Lie group},
	Amer. J. Math. {\bf 81} (1959), 973--1032. 
 \bibitem{Li} Q. Li,
	 {\em On the uniqueness of vortex equations and its geometric applications},
	 J. Geom. Anal. {\bf 29} (2019), 105--120.
	
 \bibitem{Li-Mochizuki1}
	 Q. Li and T. Mochizuki,
	 {\em Complete solutions of Toda equations and
	 cyclic Higgs bundles over non-compact surfaces},
	 Int. Math. Res. Not. IMRN 2025, no. 7, Paper No. rnaf081, 38 pp.

 \bibitem{Li-Mochizuki2}
	 Q. Li and T. Mochizuki,
	 {\em Isolated singularities of Toda equations
	 and cyclic Higgs bundles},
	 arXiv:2010.06129,
	 to appear
	 in Advanced Studies in Pure mathematics,
	 the 13th MSJ-SI proceedings
	 ``Differential Geometry and Integrable Systems''.
 \bibitem{Li-Mochizuki3}
 	Q. Li and T. Mochizuki,
	 {\em Harmonic metrics of generically regular semisimple Higgs bundles
	 on noncompact Riemann surfaces},
	 Tunis. J. Math. {\bf 5} (2023), 663--711.
 \bibitem{Matsushima-Lie-algebra}
	 Y. Matsushima,
	 {\em Theory of Lie algebras},
	 (in Japanese),
	 Kyoritsu Shuppan,
	 1956.

\bibitem{mochi4}
T. Mochizuki,
 {\em Kobayashi-Hitchin correspondence
 for tame harmonic bundles and an application},
Ast\'erisque {\bf 309} (2006).
	
\bibitem{Mochizuki-wild}
T. Mochizuki,
{\em Wild harmonic bundles and 
 wild pure twistor $D$-modules},
Ast\'{e}risque {\bf 340}, (2011)

\bibitem{Mochizuki-Toda}
	T. Mochizuki,
	{\em Harmonic bundles and Toda lattices with opposite sign I},
	the first part of arXiv:1301.1718,
	to appear in Kokyuroku Bessatsu.
	
\bibitem{Mochizuki-TodaII}
T. Mochizuki,
{\em Harmonic bundles and Toda lattices with opposite sign II},
	Comm. Math. Phys. {\bf 328} (2014), 1159--1198. 

\bibitem{Decouple}
	T. Mochizuki,
	{\em Asymptotic behaviour of certain families of
	harmonic bundles on Riemann surfaces},
	J. Topol. {\bf 9} (2016), 1021--1073.
	
\bibitem{Mochizuki-KH-Higgs}
	T. Mochizuki,
	{\em Good wild harmonic bundles and good filtered Higgs bundles},
	SIGMA Symmetry Integrability Geom. Methods Appl. {\bf 17} (2021),
	Paper No. 068, 66 pp
 \bibitem{Mochizuki-Hitchin-metric}
	 T. Mochizuki,
	 {\em Asymptotic behaviour of the Hitchin metric
	 on the moduli space of Higgs bundles},
	 arXiv:2305.17638,
	 to appear in Geometric And Functional Analysis.
	
 \bibitem{Mochizuki-Szabo}
	 T. Mochizuki,
	 Sz. Szab\'{o},
	 {\em Asymptotic behaviour of large-scale solutions of
	 Hitchin's equations in higher rank},
	 Moduli {\bf 2} (2025), Paper No. e6, 44 pp.
	
\bibitem{Nishiyama-Ohta}
	K. Nishiyama and T. Ohta,
	{\em Algebraic groups and Orbits} (in Japanese),
	Suugaku Shobo, 2015.

 \bibitem{Nori} M. V. Nori,
	 {\em On the representations of the fundamental group},
	 Compositio Math. {\bf 33} (1976), 29--41.

 \bibitem{Sagman-Smillie}
	 N. Sagman and P. Smillie,
	 {\em
	 Unstable minimal surfaces in symmetric spaces of non-compact type},
	  arXiv:2208.04885

	 
 \bibitem{Sagman-Tosic}
 N. Sagman, O. To\v{s}i\'{c},
	 {\em On Hitchin's equations for cyclic $G$-Higgs bundles},
	 Adv. Math. {\bf 482} (2025), part A, Paper No. 110599, 50 pp.

\bibitem{Serre-Lie-algebra}
	J.-P. Serre,
	{\em Complex semisimple Lie algebras},
	Translated from the French by G. A. Jones.
	Reprint of the 1987 edition.
	Springer Monographs in Mathematics. Springer-Verlag, Berlin, 2001
	
\bibitem{s1}
C. Simpson,
{\it Constructing variations of Hodge structure
using Yang-Mills theory
and application to uniformization},
J. Amer. Math. Soc. {\bf 1} (1988), 867--918.

\bibitem{s2}
C. Simpson,
{\it Harmonic bundles on non-compact curves},
J. Amer. Math. Soc. {\bf 3} (1990), 713--770.

\bibitem{s5}
C. Simpson,
{\em
 Higgs bundles and local systems},
Publ. I.H.E.S.,
{\bf 75} (1992),
 5--95.

	
\bibitem{Tauvel-Yu}
	P. Tauvel, 
	R. W. T. Yu,
	{\em Lie algebras and algebraic groups},
	Springer Monogr. Math.
	Springer-Verlag, Berlin, 2005. xvi+653 pp.

\bibitem{Vinberg}
\`{E}. B. Vinberg,
{\em The Weyl group of a graded Lie algebra}
	Izv. Akad. Nauk SSSR Ser. Mat. {\bf 40} (1976), no. 3, 488--526, 709.
	(English translation: Math. USSR-Izv. {\bf 10} (1976), 463--495 (1977).)

\end{thebibliography}
